\def\no{\if01}
\def\iftwelvept{\no}

\def\ifusepdf{\no}
\def\ifpsfont{\no}

\iftwelvept
\documentclass[leqno,12pt]{amsart}
\else
\documentclass[leqno,10pt]{amsart}
\fi
\usepackage{amssymb}
\usepackage{amscd}
\usepackage{latexsym}
\usepackage{verbatim}
\usepackage[all]{xy}

\ifusepdf
\usepackage{hyperref}
\else\fi
\ifpsfont
\usepackage[T1]{fontenc}
\usepackage{times}
\else\fi

\iftwelvept
\else\fi

\theoremstyle{plain}
\newtheorem{Theorem}{Theorem}[section]

\newtheorem{Proposition}[Theorem]{Proposition}
\newtheorem{Lemma}[Theorem]{Lemma}
\newtheorem{Corollary}[Theorem]{Corollary}

\theoremstyle{definition}

\newtheorem{Definition}[Theorem]{Definition}
\newtheorem{Remark}[Theorem]{Remark}
\newtheorem{Example}[Theorem]{Example}
\newtheorem{Conjecture}{Conjecture}

\newcommand{\ZZ}{{\mathbb{Z}}}
\newcommand{\QQ}{{\mathbb{Q}}}

\newcommand{\NN}{{\mathbb{N}}}
\newcommand{\NNNN}{\operatorname{N}}

\newcommand{\YY}{{\mathcal{Y}}}

\newcommand{\Sh}{\operatorname{Sh}}

\newcommand{\DDD}{{\mathcal{D}}}
\newcommand{\uni}{\mathbf{1}}
\newcommand{\CCC}{{\mathcal{C}}}
\newcommand{\ICoh}{{\textup{ICoh}}}
\newcommand{\PR}{\operatorname{Pr}^{\textup{L}}}

\newcommand{\SSSS}{\mathbb{S}}
\newcommand{\Rep}{\operatorname{Rep}}

\newcommand{\OO}{{\mathcal{O}}}

\newcommand{\LLL}{\mathcal{L}}

\newcommand{\MMM}{\mathcal{M}}

\newcommand{\Vect}{\operatorname{Vect}}

\newcommand{\DM}{\mathsf{DM}}

\newcommand{\Hom}{\operatorname{Hom}}
\newcommand{\Ext}{\operatorname{Ext}}

\newcommand{\Ker}{\operatorname{Ker}}

\newcommand{\Comp}{\operatorname{Comp}}
\newcommand{\Spec}{\operatorname{Spec}}

\newcommand{\Perf}{\operatorname{Perf}}

\newcommand{\Cor}{\operatorname{Cor}}
\newcommand{\SP}{\operatorname{Sp}}

\newcommand{\Mod}{\operatorname{Mod}}

\newcommand{\SSS}{\mathcal{S}}

\newcommand{\colim}{\operatorname{colim}}
\newcommand{\Cat}{\textup{Cat}_{\infty}}

\newcommand{\Map}{\operatorname{Map}}

\newcommand{\Fun}{\operatorname{Fun}}
\newcommand{\Alg}{\operatorname{Alg}}

\newcommand{\End}{\operatorname{End}}

\newcommand{\Aff}{\operatorname{Aff}}

\newcommand{\GL}{\textup{GL}}
\newcommand{\Grp}{\operatorname{Grp}}

\newcommand{\FIN}{\operatorname{Fin}_\ast}

\newcommand{\wCat}{\widehat{\textup{Cat}}_{\infty}}

\newcommand{\CAlg}{\operatorname{CAlg}}

\newcommand{\calC}{\mathcal{C}}
\newcommand{\calD}{\mathcal{D}}
\newcommand{\QC}{\operatorname{QC}}

\newcommand{\hhh}{\operatorname{h}}
\newcommand{\Ind}{\operatorname{Ind}}

\newcommand{\Coh}{\operatorname{Coh}}

\newcommand{\Aut}{\operatorname{Aut}}

\newcommand{\XX}{\mathcal{X}}
\newcommand{\Sym}{\operatorname{Sym}}
\newcommand{\Proof}{{\sl Proof.}\quad}
\newcommand{\QED}{{\unskip\nobreak\hfil\penalty50\quad\null\nobreak\hfil
{$\Box$}\parfillskip0pt\finalhyphendemerits0\par\medskip}}

\begin{document}

\title{Tannaka Duality and Stable Infinity-Categories}

\author{Isamu Iwanari}

\address{Mathematical Institute, Tohoku University, Sendai, Miyagi, 980-8578 Japan}


\thanks{The author is partially
supported by Grant-in-Aid for Scientific Research, Japan Society for the Promotion of Science.}

\email{iwanari@math.tohoku.ac.jp}

\begin{abstract}
We introduce a notion of fine Tannakian infinity-categories and prove Tannakian characterization results for symmetric monoidal stable infinity-categories over a field of characteristic zero. It connects derived quotient stacks with symmetric monoidal stable infinity-categories which satisfy a certain simple axiom.
We also discuss several applications to examples.
\end{abstract}

\maketitle

\section{Introduction}
\label{introduction}

The theory of {\it Tannakian categories} from Grothendieck-Saavedra \cite{Sa}, Deligne-Milne \cite{DM}, and
Deligne \cite{D}, \cite{D2} says that the symmetric monoidal abelian categories of representations
of a pro-algebraic group
can be characterized
as a symmetric monoidal abelian category that satisfies some categorical
conditions.
This characterization is interesting in its own right.
Grothendieck's original motivation for Tannakian categories
was to construct a motivic Galois theory of Grothendieck numerical motives.
Moreover, Tannakian theory has many applications; notably, it allows one to obtain pro-algebraic groups
from various categories, which encode the data
of categories as their representations (e.g. Picard-Vessiot theory,
and Nori's fundamental group schemes).
Similarly,
the theory of {\it Galois categories} \cite{SGA1}
by Grothendieck characterizes 
Cartesian symmetric monoidal categories of representations of
pro-finite groups.
Let us slightly reformulate the category of representations.
If $G$ is a pro-algebraic group, then any representation
of $G$ corresponds to a quasi-coherent sheaf on the classifying stack
$BG$.
That is, the symmetric monoidal category of quasi-coherent sheaves on $BG$
may be viewed as that of representations of $G$.
With this in mind, we can say that
a Tannakian theory provides a correspondence
between geometric objects (e.g. $BG$) and symmetric monoidal
categories that satisfy some condition.

\vspace{2mm}

The main results of this paper may be best understood as
Tannakian results.
Let us shift our interest to the world of
higher category theory.
The purpose of this paper is to 
establish Tannakian results for
{\it symmetric monoidal stable $\infty$-categories} \cite{HA}
with coefficients
in a field of {\it characteristic zero}.
In a sense, stable $\infty$-categories can be considered a correct generalization of
triangulated categories in the realm of $\infty$-categories
 (cf. e.g. \cite{HTT}, \cite{HA}, \cite{Be}), and
 we focus on stable $\infty$-categories
 in this paper.

Our principal result is a Tannakian characterization.
We introduce the notion of {\it fine $\infty$-categories} (or {\it fine Tannakian $\infty$-categories}).
Let $k$ be a field of characteristic zero.
Let $\mathcal{C}^\otimes$ be a $k$-linear symmetric monoidal
stable idempotent complete
$\infty$-category with the underlying $\infty$-category $\mathcal{C}$.

\begin{Definition}
\label{defwedgefinite}
Let $C$ be an object of $\CCC$.
We say that $C$ is {\it wedge-finite} (or {\it exterior-finite})
if there is a natural  number $n\ge0$
such that $\wedge^{n+1}C\simeq 0$ and $\wedge^nC$ is invertible in $\CCC^\otimes$.
We call $n$ the dimension of $C$.
Here an object $C$ of $\CCC$ is said to be invertible
if there is an object $C'$ such that $C\otimes C'\simeq C'\otimes C\simeq \uni_\CCC$ where $\uni_\CCC$ is a unit object of $\CCC^\otimes$.
The $n$-fold wedge product $\wedge^n C$ is defined to be
the image
of the idempotent map
$\textup{Alt}^n=\frac{1}{n!}\Sigma_{\sigma\in \Sigma_n}\textup{sign}(\sigma)\sigma:C^{\otimes n}\to C^{\otimes n}$, i.e. $\Ker(1-\textup{Alt}^n)$,
in the homotopy category $\textup{h}(\CCC)$
that is an idempotent complete triangulated category.
The symmetric group $\Sigma_n$ acts on $C^{\otimes n}$
by permutation. By convention, a zero object is a $0$-dimensional
wedge-finite object.
\end{Definition}

\begin{Remark}
By definition, the notion of wedge-finiteness descends to the
level of the homotopy category. Thus, this condition can be checked
at the level of symmetric monoidal triangulated categories.
Any symmetric monoidal exact functor preserves wedge-finite objects.
If the endomorphism algebra $\End_{\textup{h}(\CCC)}(\uni_{\CCC})$
of a unit object in the homotopy category
is a field, the invertibility of $\wedge^nC$ in Definition~\ref{defwedgefinite}
is automatic
(see Proposition~\ref{useful}).
\end{Remark}

\begin{Definition}
\label{fine}
Let $\mathcal{C}^\otimes$ be a $k$-linear symmetric monoidal
stable presentable
$\infty$-category.
We say that $\CCC^\otimes$ is a {\it fine $\infty$-category} over $k$
(or {\it fine Tannakian $\infty$-category})
if
\begin{enumerate}
\renewcommand{\labelenumi}{(\roman{enumi})}
\item There is a small set $\{C_{\alpha}\}_{\alpha\in A}$
of wedge-finite objects
such that $\CCC^\otimes$ is generated by $\{C_{\alpha},C_\alpha^{\vee}\}_{\alpha\in A}$ as a
symmetric monoidal
stable presentable $\infty$-category (cf. Definition~\ref{generation}).
Here $C_\alpha^\vee$ denotes the dual of $C_\alpha$ (see Remark~\ref{wedgedual}).

\item A unit object
is compact (cf. \cite[5.3.5]{HTT}, Remark~\ref{homlevelcpt}). 
\end{enumerate}
We refer to $\{C_{\alpha}\}_{\alpha\in A}$ with the property (i)
as a set of wedge-finite (or exterior-finite) generators.
Here, {\it fine} indicates
{\it fin}iteness + {\it e}xterior-product.
If no confusion seems likely to arise,
we omit ``over $k$''.
\end{Definition}

\begin{Remark}
\label{wedgedual}
By Theorem~\ref{characterization},
every wedge-finite object is dualizable (see Remark~\ref{wedgedual2}).
\end{Remark}

Our characterization theorem is the following
statement
(cf. Theorem~\ref{algebraic}, Theorem~\ref{main1}):

\begin{Theorem}[Characterization theorem]
\label{intmain1}
Let $\mathcal{C}^\otimes$ be a $k$-linear symmetric monoidal
stable presentable
$\infty$-category.
The following statements are equivalent to one another:
\begin{enumerate}
\renewcommand{\labelenumi}{(\theenumi)}

\item $\CCC^{\otimes}$ is a fine $\infty$-category.

\item There exist a derived quotient stack $X=[\Spec A/G]$
where a pro-reductive group $G$ acts on an affine derived scheme $\Spec A$ with $A$
a commutative differential graded algebra, and a $k$-linear symmetric monoidal equivalence
$\CCC^\otimes\simeq\QC^\otimes(X)$.
Here $\QC^\otimes(X)$ denotes the symmetric monoidal stable
 $\infty$-category of quasi-coherent complexes on $X$ (see Section~\ref{QCC}).
\end{enumerate}
\end{Theorem}

A derived stack is a stack in the theory of derived algebraic geometry,
which is a generalization of classical
algebraic geometry \cite{DAG7}, \cite{HAG2}
that uses homotopy-theoretic ideas and techniques.
Here, we think of derived stacks of the form $[\Spec A/G]$ appearing in Theorem~\ref{intmain1}
as the generalization of classifying stacks of affine
group schemes as well as
a nice class of derived stacks.
This characterization makes it possible to
obtain a derived stack $X=[\Spec A/G]$  from
an abstract symmetric monoidal stable $\infty$-category.
We consider not only gerbe-like stacks,
but also the class of nice derived quotient stacks,
which allows access to the extensive power of derived algebraic geometry.
More importantly, 
our construction of a derived quotient stack (from a fine $\infty$-category
with a given set of wedge-finite generators)
is quite explicit, and  the associated stack has a specific form;
see Section~\ref{TC}.

Fine $\infty$-categories are defined by reasonably simple conditions,
which we can use to find examples in practice.
The recent fascinating development
of higher category theory has highlighted various examples of symmetric monoidal stable $\infty$-categories.
Among them, the next result proves that the following symmetric monoidal
$\infty$-categories are fine $\infty$-categories
(see Section~\ref{EX} for details):

\begin{Theorem}
\label{exlist}
The following symmetric monoidal
$\infty$-categories are examples of fine $\infty$-categories:
\begin{enumerate}
\renewcommand{\labelenumi}{(\roman{enumi})}

\item The unbounded derived $\infty$-category of representations of a
pro-algebraic group over a field of characteristic zero (cf. Proposition~\ref{alggp}, Remark~\ref{alggpR}),

\item The stable $\infty$-category of mixed motives generated by Kimura finite dimensional Chow motives (cf. Theorem~\ref{motuncond}),

\item The stable $\infty$-category of noncommutative mixed motives generated by
Kimura finite dimensional noncommutative motives (cf. Proposition~\ref{noncommutativemot}),

\item The stable $\infty$-category of Ind-coherent complexes on a topological
space (cf. Proposition~\ref{RRHT1}),

\item The unbounded derived $\infty$-category of quasi-coherent complexes on a
quasi-projective variety (cf. Theorem~\ref{dSerre}).
\end{enumerate}
\end{Theorem}

Examples (ii) and (iii) are of great interest in view of the
motivic Galois theory
of mixed motives.
A striking aspect of (ii) is that it reveals
an intimate relation
between the Kimura finiteness of motives (see \cite{Kim}, \cite{A}, \cite{A2}, \cite{Ivo}) and Theorem~\ref{intmain1} about
fine $\infty$-categories (see Section~\ref{exmot}).
Example (iv) is closely
related to rational homotopy theory (see Section~\ref{exRHT}).
We will prove that if $S$ is a topological space that
satisfies an appropriate property, the higher
rational homotopy groups and pro-algebraic completion
of the fundamental group
can be recovered in a Tannakian way from the stack associated with
the stable $\infty$-category of Ind-coherent complexes on $S$.
As an illustration, we briefly describe two examples of Galois categories: the category of finite
covers of a topological space and the category
of finite \'etale covers of a field.  Namely, the theory of Galois categories
simultaneously generalizes the fundamental groups of topological spaces
and the classical Galois theory:
\[
\textup{(classical Galois theory)}\leftarrow \textup{(Galois categories)}\to \textup{($\widehat{\pi}_1$ of topological spaces)}.
\]
Motivic Galois theory generalizes the classical Galois
theory to motives generated by general algebraic varieties.
Rational homotopy theory is the homotopy theory for rational homotopy types,
in which rational homotopy groups play a central role.
Hence, examples (ii) and (iv)
of fine $\infty$-categories could be considered
a ``higher'' generalization of the above, as follows:
\[
\textup{(motivic Galois theory)}\leftarrow \textup{(fine Tannakian $\infty$-categories)}\to \textup{(rational homotopy theory)}.
\]

Another important aspect of Theorem~\ref{intmain1} is the following:
A theorem of Schwede and Shipley \cite{SS} states roughly that
if $\mathcal{C}$ is
a stable presentable $\infty$-category (without any sort of monoidal structure)
that has a compact generator, then there is a 
(not necessarily commutative) ring spectrum $A$
such that $\mathcal{C}$ is equivalent to the
stable presentable $\infty$-category of modules over $A$.
Thus, one may think of Theorem~\ref{intmain1} as
a symmetric monoidal generalization of this theorem of Schwede and Shipley
(but the characteristic zero assumption is essential for our result).

\vspace{2mm}

The main difficulty in the proof of Theorem~\ref{intmain1}
arises from the fact that $\QC^\otimes([\Spec A/G])$
(or a given symmetric monoidal stable $\infty$-category)
does not have a Tannakian category or the like
as its full subcategory in general, and so
it is hard to rely on the classical Tannakian
theory and methods in our setting.
We use a new method of characterizing the derived $\infty$-category of
representations of a general linear group $\GL_d$ by a universal
property. This may be of independent interest, but is also a key ingredient
to the proof of
Theorem~\ref{intmain1} (cf. Theorem~\ref{characterization}):

\begin{Theorem}[A universal property]
\label{intchara}
Let $\mathcal{C}^\otimes$ be a $k$-linear symmetric monoidal
stable presentable
$\infty$-category.
Let $\mathcal{C}_{\wedge,d}$ be the full subcategory of $\mathcal{C}$
that consists of $d$-dimensional
wedge-finite (exterior-finite) objects in
$\mathcal{C}^\otimes$. Let $\mathcal{C}_{\wedge,d}^{\simeq}$ be
the largest Kan subcomplex (i.e. $\infty$-groupoid) of $\mathcal{C}_{\wedge,d}$,
obtained
by restricting to those
morphisms which are equivalences.
Then there exists a natural homotopy equivalence of spaces
\[
\Map_{\CAlg(\textup{Pr}_k^{\textup{L}})}(\QC^\otimes(B\GL_d), \mathcal{C}^\otimes) \to \mathcal{C}_{\wedge,d}^{\simeq}
\]
which carries $f:\QC^\otimes(B\GL_d) \to \CCC^\otimes$
to the image $f(K)$ of the standard representation $K$ of $\GL_d$.
That is, an object $C\in \mathcal{C}_{\wedge,d}$ corresponds to
a $k$-linear symmetric monoidal functor
$\QC^\otimes(B\GL_d)\to \mathcal{C}^\otimes$ that sends $K$ to $C$.
\end{Theorem}

The classical Tannakian theory tells us that
for a pro-algebraic group $G$ over $k$ and a $k$-algebra $R$, the groupoid
$\Map_{k-\textup{stacks}} (\Spec R, BG)$ of morphisms to $BG$
is equivalent to $\Map_{k}^\otimes (\textup{qcoh}^\otimes(BG), \textup{qcoh}^\otimes(\Spec R))$
of $k$-linear symmetric monoidal exact functors between symmetric monoidal
abelian categories of quasi-coherent sheaves.
That is, $f:\Spec R\to BG$
corresponds to $f^*:\textup{qcoh}^\otimes(BG)\to \textup{qcoh}^\otimes(\Spec R)$, cf. \cite{DM} for precise details.
Its analogue for derived $\infty$-categories of schemes and Deligne-Mumford
stacks is proved in \cite{FI}.
We now invite the reader's attention to the fact
that in the setting of our derived (Artin) stacks,
symmetric monoidal functors
do {\it not} correspond to morphisms of stacks.
There exists a symmetric monoidal functor
which is not the pullback functor of a morphism of stacks:
Let $B\mathbb{G}_m$ be the usual classifying stack of the algebraic torus
$\mathbb{G}_m$.
We have a symmetric monoidal equivalence
\[
\QC^\otimes(B\mathbb{G}_m) \to \QC^\otimes(B\mathbb{G}_m)
\]
which carries each character $\chi_n$ of weight $n$ of $\mathbb{G}_m$ to $\chi_n[2n]$. However, it does not arise
as the pullback functor of any morphism $B\mathbb{G}_m\to B\mathbb{G}_m$ (because it does not
preserve the heart of standard $t$-structure).
To analyze this exotic and new phenomenon\footnote{Note that such phenomena naturally appear in some generalizations of Tannaka duality. For example, it appears in Tannaka duality for (braided) monoidal categories, which involves a Drinfeld associator and twist.}, inspired by \cite{FI}
we introduce the geometric notion of {\it correspondences}
between derived stacks. A correspondence from $X$ to $Y$
can be viewed as a twisted morphism and is defined in a similar way to algebraic correspondences, thus capturing the phenomenon.
That is, we prove that correspondences (rather than
morphisms) corresponds to symmetric monoidal functors (see Section~\ref{correspondence}):

\begin{Theorem}[Symmetric monoidal functors versus correspondences]
\label{intmain}
Let $X$ and $Y$ be two quotient stacks
of the forms $[\Spec A/G]$ and $[\Spec B/H]$ respectively,
where $A,B\in \CAlg_k$, and $G$ and $H$ are pro-reductive groups
over $k$.
There is a natural equivalence of $\infty$-groupoids
\[
\Map_{\textup{Cor}_k}(X,Y)\to \Map_{\CAlg(\PR_k)}(\QC^\otimes(Y),\QC^\otimes(X)); f\mapsto f^*.
\]
Here, the left-hand side is the spaces of correspondences from
$X$ to $Y$ (defined in Section~\ref{correspondence}).
Moreover, the composition of symmetric monoidal functors
corresponds to a composition of correspondences.
\end{Theorem}

Recall that there are two aspects of Tannakian theory of
a given symmetric monoidal category $\CCC^\otimes$.
One is to think of $\CCC^\otimes$ as the category of sheaves
on a geometric object (or the representation category of a group object).
The other is to consider the group object that represents
the automorphism group of $p$
when $\CCC^\otimes$ is equipped with a ``fiber functor'' $p$.
Let us focus on the second aspect.
Suppose that a symmetric monoidal stable $\infty$-category $\CCC^\otimes$ is
equipped with
a symmetric monoidal functor $p:\CCC^\otimes\to \Mod_k^\otimes$ to a symmetric monoidal
stable $\infty$-category of $Hk$-module spectra.
In \cite{Tan}, we constructed a derived
affine group scheme
that represents the automorphism group $\Aut(p)$ of $p$. We refer to  \cite{Tan}, \cite{Bar} for details.
When $\CCC^\otimes$ is a fine $\infty$-category (and thus $\CCC^\otimes\simeq \QC^\otimes([\Spec A/G])$),
one can apply the construction of a based loop space for $[\Spec A/G]$,
under a suitable condition,
to obtain a derived affine group scheme $\mathsf{G}:=\Spec k\times_{[\Spec A/G]}\Spec k=\Omega_*[\Spec A/G]$ that represents the automorphism
group of $p$ (see Remark~\ref{Tandual} and~\ref{motGalcon}, Section~\ref{exRHT}). The fiber product $\Spec k\times_{[\Spec A/G]}\Spec k$
can naturally be regarded as a $G$-equivariant version of the bar construction.
This derived group scheme $\mathsf{G}$ is  the Tannaka dual
of $\CCC^\otimes$ with respect to $p$.

\vspace{4mm}

Now we view our results from the perspective of the motivic Galois
theory.
The results in this paper have applications to motivic Galois
theory of mixed motives.
To explain this, we briefly describe the motivic Galois theory
of {\it mixed Tate motives}.
Let $\mathbf{DM}_{gm}^\otimes$ be the $\QQ$-linear symmetric monoidal triangulated category
of mixed motives over a perfect field
$K$. Here, we adopt the category of Voevodsky's motives \cite{Vtri}.
The symmetric monoidal triangulated category $\mathbf{DTM}^\otimes_{gm}$ of mixed Tate motives 
is defined to be the triangulated subcategory of $\mathbf{DM}_{gm}$
generated by the motive of the one dimensional torus $\mathbb{G}_m$,
which is closed under taking retracts, tensor products and duals.
This category has a natural enhanced formulation in $\infty$-categories,
namely, a symmetric monoidal
stable presentable subcategory of mixed Tate motives
$\mathsf{DTM}^\otimes$ in the symmetric monoidal
stable presentable $\infty$-categories
$\mathsf{DM}^\otimes$.
Here, we switch to the formulation of presentable $\infty$-categories.
In \cite{LevT} and \cite{Spi},
it is essentially proved that
there exists an augmented commutative differential graded $\QQ$-algebra $A$ endowed with an
action of $\mathbb{G}_m$ and a $\QQ$-linear
symmetric monoidal equivalence $\mathsf{DTM}^\otimes\simeq \QC^\otimes([\Spec A/\mathbb{G}_m])$. This is a Tannakian theorem for mixed Tate motives.
The motivic Galois group of mixed Tate motives is constructed from $[\Spec A/\mathbb{G}_m]$ through the $\mathbb{G}_m$-equivariant bar construction.
Namely, in geometric terms, it is the truncation of
the derived affine group scheme
$\Spec \QQ\times_{[\Spec A/\mathbb{G}_m]}\Spec \QQ$, which is an ordinary
pro-algebraic group.
Ultimately, the conjectural motivic Galois theory (cf. \cite{A}) suggests that
this well-established Galois theory of mixed Tate motives $\mathsf{DTM}^\otimes$ should be generalized
to all mixed motives $\mathsf{DM}^\otimes$.
The main obstacle to generalization is
an extension of the above Tannakian theorem for mixed Tate motives.

Our Tannakian result can be applied to the stable presentable subcategory
generated by Kimura finite dimensional motives. 
The motives of abelian varieties and $\mathbb{G}_m$
(more generally, semi-abelian varieties) are Kimura finite,
and the mixed motives generated by the motives of abelian varieties
are the important class of {\it mixed elliptic motives} and
{\it its higher dimensional generalization}, i.e., mixed abelian motives.
Thus, by Theorem~\ref{intmain1} (see also Theorem~\ref{algebraic}),
for example, if $X$ is an abelian variety and
$\mathsf{DM}^\otimes_X$ denotes the smallest symmetric monoidal
stable presentable subcategory
in $\mathsf{DM}^\otimes$ containing the motives of $X$ and its dual,
we
unconditionally obtain a commutative differential graded algebra $B$ endowed
with an action of $\GL_d$ and a $\QQ$-linear symmetric monoidal equivalence
$\mathsf{DM}^\otimes_X\simeq \QC^\otimes([\Spec B/\GL_d])$.
It is worth mentioning that even for the case of mixed Tate motives,
the methods in this paper are new.
Moreover, as we will
explain in Section~\ref{exmot}, by using a realization functor of mixed motives,
one can apply a bar construction to obtain a motivic Galois group
of such a motive.
In a subsequent paper \cite{PM}, by combining the Tannakian results of this
paper with results on
Galois representations,
we prove a structure theorem of the motivic Galois group of such motives,
which was predicted by the conjectural perspective of mixed motives
of Beilinson and Deligne.
The results of this paper
can also be applied to a motivic generalization of
rational homotopy theory
\cite{MRHT}, the universal family of
a modular variety, and forth. Thus, one may expect more applications and
other directions.

\vspace{3mm}

Next we consider recent progress on Tannakian theory for
symmetric monoidal stable $\infty$-categories {\it endowed with $t$-structures}.
Lurie \cite[Section 4,5]{DAG8} establishes
a Tannakian theory of symmetric monoidal stable $\infty$-categories with
coefficients in a field of characteristic zero
which are endowed with $t$-structures and satisfy some conditions
(locally dimensional $\infty$-categories),
and in \cite{Wa} a version of Tannakian theory for stable
$\infty$-categories over ring spectra
equipped with $t$-structures and fiber functors is developed.
In \cite{DAG8}, Tannakian results for morphisms of stacks
were proved where the data of $t$-structures is essential.
As well as the motivation from motives,
Deligne's idea \cite{D}, \cite{D2} and Lurie's idea
on
beautiful internal characterizations of Tannakian (and super-Tannakian, locally dimensional)
categories
without fiber functors influence our work.
Meanwhile, as one can easily imagine,
there are substantial differences between the present paper and theories
taking account of $t$-structures.
First, if a symmetric stable $\infty$-category is endowed with $t$-structure,
its heart is a Tannakian category (or a suitable symmetric monoidal
abelian category) under
an appropriate condition on $t$-structure.
Thus, unlike the setting of this paper, one can rely on the classical theory of Tannakian category or a similar argument.
Second, since we do not assume $t$-structures,
Theorem~\ref{intmain1} is relatively easy to apply. For example,
this generalization is crucial for unconditional applications to
mixed motives
(cf. \cite{PM}, \cite{MRHT}, Section~\ref{EX}).
Third, as observed above, symmetric monoidal functors of fine $\infty$-categories correspond
not to morphisms of derived stacks but to correspondences.

As a byproduct of our Tannakian result,
we 
associate
a fine $\infty$-category
with a topological space,
and study the rational homotopy groups
and the pro-algebraic completion of the fundamental group
by means of the fine $\infty$-category 
and the associated stack (cf. Section~\ref{exRHT}),
although we do not obtain a reconstruction of a rational homotopy type
from the associated fine $\infty$-category.
There have been various works on rational homotopy
theory for non-nilpotent spaces, for example, see
Bousfield-Kan \cite{BK}, Brown-Szczarba \cite{BS},
G\'omez-Tato-Halperin-Tanr\'e \cite{GHT}, To\"en \cite{aff}, Pridham
\cite{P}, and Moriya \cite{Mori} (cf. Remark~\ref{classicRHT}).
These works proposed several different definitions of a rational
homotopy type for a non-nilpotent space.
It might be quite interesting
to study fine $\infty$-categories associated with topological spaces
from the viewpoint of these definitions of rational homotopy types.

\vspace{3mm}

This paper is organized as follows.
In Section 2, we recall/prepare basic definitions and results
about derived stacks, symmetric monoidal stable $\infty$-categories,
quasi-coherent complexes, and so on.
In Section 3, we discuss a universal characterization of
the derived $\infty$-category of representations of a general linear group
in terms of wedge-finite objects. Namely, we prove Theorem~\ref{intchara}.
In Section 4, we prove Theorem~\ref{intmain1} and its algebraic
version Theorem~\ref{algebraic}. Moreover,
we study an explicit presentation of the derived stack
associated to a fine $\infty$-category with a prescribed
wedge-finite generator.
In Section 5, we introduce correspondences between derived stacks
and prove Theorem~\ref{intmain}. The reader can skip this Section
for the first reading.
In Section 6, we present some examples of fine $\infty$-categories.
We discuss (i) the relation with the classical Tannakian categories,
(ii) unconditional application to stable $\infty$-category of mixed motives,
which opens up a nice relationship with Kimura finiteness of motives, (iii) derived
$\infty$-category of quasi-coherent sheaves on a quasi-projective variety,
(iv) a Tannakian theory of
coherent sheaves on a topological space in the context of the rational
homotopy theory.

\vspace{3mm}

{\it Convention and notation.}
Throughout this paper, we use the theory of {\it quasi-categories}.
A quasi-category is a simplicial set which
satisfies the weak Kan condition of Boardman-Vogt.
The theory of quasi-categories from the viewpoint of higher category theory
were extensively developed by Joyal and Lurie \cite{Jo}, \cite{HTT}, \cite{HA}.
Following \cite{HTT} we shall refer to quasi-categories
as {\it $\infty$-categories}.
Our main references are \cite{HTT}
 and \cite{HA}.
For the brief introduction to $\infty$-categories, we
refer to \cite[Chapter 1]{HTT}, \cite{Gro}, \cite[Section 2]{FI}.
For the quick survey on various approaches to $(\infty,1)$-categories
(e.g. simplicial categories, Segal categories, complete Segal spaces, etc.)
and their relations,
we refer to \cite{Be}.
As a set-theoretic foundation,
we employ the axiom of ZFC together with the axiom of Grothendieck universes
(i.e., every Grothendieck universe is an element of a larger universe).
We fix a sequence of universes
$(\NN\in) \mathbb{U}\in \mathbb{V}\in \mathbb{W}\in \ldots$ and refer to
sets belonging to $\mathbb{U}$ (resp. $\mathbb{V}$, $\mathbb{W}$) to
as small sets (resp. large sets, super-large sets).
But in the text
we avoid using the notation $\mathbb{U}$, $\mathbb{V}$, $\mathbb{W}$.
To an ordinary category, we can assign an $\infty$-category by taking
its nerve, and therefore
when we treat ordinary categories we often omit the nerve $\NNNN(-)$
and directly regard them as $\infty$-categories.
We often refer to a map $S\to T$ of $\infty$-categories
as a functor. We call a vertex in an $\infty$-category $S$
(resp. an edge) an object (resp. a morphism).
Here is a list of (some) of the conventions and notation that we will use:

\begin{itemize}

\item $\Delta$: the category of linearly ordered non-empty
finite sets (consisting of $[0], [1], \ldots, [n]=\{0,\ldots,n\}, \ldots$)

\item $\Delta^n$: the standard $n$-simplex

\item $\textup{N}$: the simplicial nerve functor (cf. \cite[1.1.5]{HTT})

\item $\mathcal{C}^{op}$: the opposite $\infty$-category of an $\infty$-category $\mathcal{C}$

\item Let $\mathcal{C}$ be an $\infty$-category and suppose that
we are given an object $c$. Then $\mathcal{C}_{c/}$ and $\mathcal{C}_{/c}$
denote the undercategory and the overcategory, respectively (cf. \cite[1.2.9]{HTT}).

\item $\CCC^\simeq$: the largest Kan subcomplex (contained) in an $\infty$-category
$\CCC$,
that is, the Kan complex obtained from $\CCC$
by restricting to those
morphisms (edges) which are equivalences.

\item $\operatorname{Cat}_\infty$: the $\infty$-category of small $\infty$-categories

\item $\wCat$: $\infty$-category of large $\infty$-categories

\item $\SSS$: $\infty$-category of small spaces. We denote by $\widehat{\SSS}$
the $\infty$-category of large spaces (cf. \cite[1.2.16]{HTT}).

\item $\textup{h}(\mathcal{C})$: homotopy category of an $\infty$-category (cf. \cite[1.2.3.1]{HTT})

\item $\Ind(\mathcal{C})$: $\infty$-category of Ind-objects
in an $\infty$-category $\mathcal{C}$
(see \cite[5.3.5.1]{HTT}, \cite[4.8.1.14]{HA} for the symmetric monoidal setting).

\item $\Fun(A,B)$: the function complex for simplicial sets $A$ and $B$

\item $\Fun_C(A,B)$: the simplicial subset of $\Fun(A,B)$ classifying
maps which are compatible with
given projections $A\to C$ and $B\to C$.

\item $\Map(A,B)$: the largest Kan subcomplex of $\Fun(A,B)$ when $A$ and $B$ are $\infty$-categories.

\item $\Map_{\mathcal{C}}(C,C')$: the mapping space from an object $C\in\mathcal{C}$ to $C'\in \mathcal{C}$ where $\mathcal{C}$ is an $\infty$-category.
We usually view it as an object in $\mathcal{S}$ (cf. \cite[1.2.2]{HTT}).

\end{itemize}

{\it Stable $\infty$-categories, symmetric monoidal $\infty$-categories and spectra}.
For the definitions of (symmetric) monoidal $\infty$-categories,
$\infty$-operads, and
their algebra objects, we shall refer to \cite{HA}.
A stable $\infty$-category is an $\infty$-category which satisfies
the conditions (i) there is a zero object, i.e., an object
which is both initial and final, (ii) every morphism has a fiber and a cofiber,
(iii) for any sequence $X\stackrel{f}{\to} Y \stackrel{g}{\to} Z$ of morphisms,
$X$ is a fiber of $g$ if and only if $Z$ is a cofiber of $f$ (see \cite[1.1.1.9]{HA}).
Our reference for stable $\infty$-categories is
\cite[Chapter 1]{HA}. We list some further notation.

\begin{itemize}

\item $\Mod_R$: $\infty$-category of $R$-module spectra
for a commutative ring spectrum $R$. When $R$ is the Eilenberg-MacLane
spectrum $Hk$ of an ordinary commutative ring $k$, we write $\Mod_k$ for
$\Mod_{R}$ (thus $\Mod_k$ is not the category of usual $k$-modules).
If $\mathcal{D}(k)$ denotes the stable $\infty$-category
obtained from the category
of (possibly unbounded) chain complexes of $k$-modules by inverting
quasi-isomorphisms, there is a canonical
equivalence $\Mod_k\simeq \mathcal{D}(k)$ (of symmetric monoidal
$\infty$-categories), see \cite[7.1.2, 7.1.2.13]{HA} for more details.

\item $\FIN$: the category of pointed finite sets $\langle 0\rangle=\{\ast\}, \langle 1\rangle=\{1, \ast\},\ldots,\langle n\rangle=\{1\ldots,n, \ast\},\ldots$. A morphism
is a map $f:\langle n\rangle\to \langle m \rangle$ such that
$f(\ast)=\ast$. Note that $f$ is not assumed to be order-preserving.

\item Let $\mathcal{M}^\otimes\to \mathcal{O}^\otimes$ be a fibration of $\infty$-operads. We denote by
$\Alg_{/\mathcal{O}^\otimes}(\mathcal{M}^\otimes)$ the $\infty$-category of algebra objects (cf. \cite[2.1.3.1]{HA}).  We often write $\Alg(\MMM^\otimes )$ or 
$\Alg(\MMM)$ for $\Alg_{/\mathcal{O}^\otimes}(\MMM^\otimes)$.
Suppose that $\mathcal{P}^\otimes\to \mathcal{O}^\otimes$ is a map of $\infty$-operads.
We write $\Alg_{\mathcal{P}^\otimes/\mathcal{O}^\otimes}(\mathcal{M}^\otimes)$
for the $\infty$-category of $\mathcal{P}$-algebra objects.

\item $\CAlg(\mathcal{M}^\otimes)$: $\infty$-category of commutative
algebra objects in a symmetric
monoidal $\infty$-category $\mathcal{M}^\otimes\to \NNNN(\FIN)$.
When the symmetric monoidal structure is clear,
we usually write $\CAlg(\mathcal{M})$ for $\CAlg(\mathcal{M}^\otimes)$.

\item $\CAlg_R$: $\infty$-category of commutative
algebra objects in the symmetric monoidal $\infty$-category $\Mod_R^\otimes$
where $R$ is a commutative ring spectrum. When $R$ is the sphere
spectrum $\mathbb{S}$, we set
$\CAlg=\CAlg_{\mathbb{S}}$.
When $R$ is the Eilenberg-MacLane spectrum $Hk$ with $k$ a commutative ring,
then we write $\CAlg_k$ for $\CAlg_R$.
If $k$ is a field of characteristic zero,
the $\infty$-category $\CAlg_k$ is
equivalent to the $\infty$-category
obtained from the model category of commutative
differential graded $k$-algebras by inverting quasi-isomorphisms
(cf. \cite[7.1.4.11]{HA}, \cite{Hi}). Therefore
we often refer to objects in $\CAlg_k$ as
commutative differential graded algebras.

\item $\Mod_A^\otimes(\mathcal{M}^\otimes)\to \NNNN(\FIN)$: symmetric monoidal
$\infty$-category of
$A$-module objects,
where $\mathcal{M}^\otimes$
is a symmetric monoidal $\infty$-category such that (1)
the underlying $\infty$-category admits a colimit for any simplicial diagram, and (2)
its tensor product functor $\mathcal{M}\times\mathcal{M}\to \mathcal{M}$
preserves
colimits of simplicial diagrams separately in each variable.
Here $A$ belongs to $\CAlg(\mathcal{M}^\otimes)$.
cf. \cite[3.3.3, 4.5.2]{HA}.

\end{itemize}

\begin{Definition}
\label{presentably}
Let $\CCC^\otimes$ be a symmetric monoidal $\infty$-category with
the underlying $\infty$-category $\CCC$.
We say that $\CCC^\otimes$
is a presentably symmetric monoidal $\infty$-category if the following
two conditions are satisfied:
\begin{itemize}
\item The underlying $\infty$-category $\CCC$ is presentable.
\item The tensor product functor $\CCC\times \CCC\to \CCC$ preserves
small colimits separately in each variable.
\end{itemize}
\end{Definition}

\begin{Definition}
\label{generation}
Let $\CCC$ be a stable presentable $\infty$-category.
Let $\{C_{\alpha}\}_{\alpha \in A}$ be a small set of objects
in $\CCC$.
We say that
$\{C_{\alpha}\}_{\alpha \in A}$ generates $\CCC$
as a stable presentable $\infty$-category
if $\CCC$ is the smallest stable subcategory which contains
$\{C_{\alpha}\}_{\alpha \in A}$ and
is closed under
small coproducts.

Suppose that $\CCC^\otimes$ is a presentably symmetric monoidal stable $\infty$-category (cf. Definition~\ref{presentably}).
We say that
$\{C_{\alpha}\}_{\alpha \in A}$ generates $\CCC^\otimes$
as a symmetric monoidal stable presentable $\infty$-category
if $\CCC$ is the smallest stable subcategory which contains
the unit object and $\{C_{\alpha}\}_{\alpha \in A}$ and
is closed under
small coproducts and tensor products.
(We remark that
any stable $\infty$-category which has small coproducts admits
all small colimits.)
\end{Definition}

\begin{Remark}
\label{homlevelcpt}
If each object $C_{\alpha}$
is compact and $\{C_{\alpha}\}_{\alpha\in A}$
generates $\CCC$ as a stable presentable $\infty$-category,
we say that the stable presentable $\infty$-category
$\CCC$ is compactly generated.
This notion is compatible with the notion of compactly generated
triangulated category. Namely,
the compactness of $C_{\alpha}$ in $\CCC$ and that in the triangulated
category $\textup{h}(\CCC)$ coincide, and
$\textup{h}(\CCC)$ is the smallest triangulated subcategory of $\textup{h}(\CCC)$
which contains $\{C_{\alpha}\}_{\alpha\in A}$ and is closed under small coproducts if and only if $\{C_{\alpha}\}_{\alpha\in A}$
generates $\CCC$ as a stable presentable $\infty$-category.
In addition, if each object $C_{\alpha}$
is compact, these conditions are equivalent to the following condition:
for any $C\in \CCC$,
the vanishing $\Hom_{\textup{h}(\CCC)}(C_{\alpha},C[r])=0$ for any pair $(\alpha,r)\in A\times \ZZ$
implies $C\simeq 0$. Our references are \cite[2.2.1]{SS}, \cite[1.4.4.3]{HA}.
\end{Remark}

\section{Preliminaries on stacks and quasi-coherent complexes}
\label{presection}

In this Section, we will recall
some definitions and prepare several results
concerning derived stacks, symmetric monoidal stable $\infty$-categories, and so on.

\subsection{Derived stacks}

Let $k$ be a field of characteristic zero.
First, we recall notions of flatness, \'etaleness, etc.
Let $\CAlg_k$ be the $\infty$-category of commutative ring spectra
over the Eilenberg-MacLane spectrum $Hk$.
We usually identify $\CAlg_k$ with the $\infty$-category
obtained from the category of commutative differential graded
$k$-algebras by inverting quasi-isomorphisms.
A morphism $\phi:A\to B$ in $\CAlg_k$
is said to be flat (resp. faithfully flat, resp. \'etale)
if the induced morphism of usual commutative algebras
$H_0(A)\to H_0(B)$ is flat (resp. faithfully flat, resp. \'etale) in the classical sense,
and the canonical morphism
\[
H_0(B)\otimes_{H_0(A)}H_n(A)\to H_n(B)
\]
is an isomorphism for any $n\in \ZZ$.
A morphism $\phi:A\to B$ is said to be formally perfect
if the cotangent complex $L_{B/A}$ is dualizable in $\Mod_B$
(see \cite[1.2.7.1]{HAG2} and \cite[7.3]{HA}
for the notions of formally perfectness
and cotangent complexes).

Next, write $\Aff_k:=\CAlg_k^{op}$.
We refer to $\Aff_k$ as the $\infty$-category of derived
affine schemes
over $k$. We denote by $\Spec R$ the object in $\Aff_k$
corresponding to $R$ in $\CAlg_k$.
We say that $\Spec B\to \Spec A$ is flat (resp. faithfully flat,
resp. \'etale, resp. \'etale surjective, resp. formally perfect)
if the corresponding morphism $A\to B$ is flat
(resp. faithfully flat,
resp. \'etale, resp. \'etale and faithfully flat, resp. formally perfect).

We say that a functor $F:\Aff_k^{op} \to \widehat{\SSS}$
is a flat (or fpqc) sheaf (resp. an \'etale sheaf) if
the following conditions are satisfied:
\begin{itemize}
\item For any finite coproduct $\sqcup_{i\in I}\Spec R_i$
in $\Aff_k$,
$F(\sqcup_{i\in I}\Spec R_i)\simeq \prod_{i\in I}F(A_i)$.

\item For any flat hypercovering (resp. any \'etale hypercovering)
$\Spec B^\bullet\to \Spec A$,
$F(A)\simeq \varprojlim_n F(B^n)$.

\end{itemize}
Here a flat hypercovering (resp. an \'etale hypercovering) of $\Spec A$ is
an augmented simplicial diagram of derived affine schemes
$\Spec B^\bullet \to \Spec A$
such that for any $n\ge0$, $\Spec B^{n+1}\to (\textup{cosk}_{n}\Spec B^\bullet)_{n+1}$
is faithfully flat
(resp. \'etale surjective) and $\Spec B^0\to \Spec A$ is faithfully flat
(resp. \'etale surjective), where the coskeleton is taken in $(\Aff_k)_{/\Spec A}$.
Strictly speaking, a sheaf in this sense should be
referred to as a hypercomplete sheaf (see e.g.
\cite[5.12]{DAG7}) but we usually omit ``hypercomplete''
in this paper.

Let $\Sh(\Aff_k)$ be the full subcategory of $\Fun(\CAlg_k,\widehat{\SSS})$ spanned by \'etale sheaves.
By the Yoneda Lemma, there is a fully faithful
functor $\Aff_k\to \Fun(\CAlg_k,\widehat{\SSS})$.
The essential image is contained in $\Sh(\Aff_k)$. By abuse of notation,
we often think of $\Aff_k$ as its essential image in $\Sh(\Aff_k)$.

\begin{Definition}
\label{defderivedstack}
Let $X:\CAlg_k\to \widehat{\SSS}$ be
a sheaf, that is, an object of $\Sh(\Aff_k)$.

\begin{enumerate}
\renewcommand{\labelenumi}{(\roman{enumi})}

\item The sheaf $X$ is said to be a derived stack
if there is a groupoid object
$X_\bullet:\NNNN(\Delta)^{op}\to \Aff_k$
(see e.g. \cite[6.1.2.7]{HTT}, \cite[1.3.1.6]{HAG2} for groupoid objects) such that $X$ is
equivalent to a colimit
of the composite $\NNNN(\Delta)^{op}\to \Aff_k\to \Sh(\Aff_k)$.
We refer to $X_\bullet$ as a presentation of $X$.

\item If there is a presentation $X_\bullet$ such that
$d_0:X_1\to X_0$ is formally perfect (equivalently $d_1:X_1\to X_0$ is
formally perfect),
we say that $X$ is a derived algebraic stack.
Here, $d_i:X_1\to X_0$ is determined by the inclusion
$d^i:[0]\to [1]$.

\item
A morphism $X\to Y$ of derived stacks is a morphism in $\Sh(\Aff_k)$.

\item
A morphism $X\to Y$ in $\Sh(\Aff_k)$ is said to be affine if
for any $\Spec R\to Y$, the fiber product $\Spec R\times_{Y}X$
belongs to $\Aff_k$.

\end{enumerate}

\end{Definition}

\begin{Remark}
The existence of a presentation
is a quite weak condition. In this paper,
we use this condition to prove that the $\infty$-category of
quasi-coherent complexes on a derived stack is presentable.

The notion of derived algebraic stacks is closely related to
the notion of geometric stacks in the setting of complicial
algebraic geometry introduced in To\"en-Vezzosi \cite{HAG2}
(in {\it loc. cit.} the theory is  developed by means of the model category theory).
Suppose that $X$ is a derived algebraic stack
and $X_\bullet$ is a presentation of $X$ such that
$d_0:X_1\to X_0$ is formally perfect.
This presentation is a $(-1)$-$\mathbf{P}_w$ (Segal) groupoid
in the sense of {\it loc. cit.}
(see \cite[1.3.4.1, 2.3.2]{HAG2} for the notion of
$\mathbf{P}_w$ groupoids). 
According to \cite[1.3.4.2]{HAG2},
$X$ is
a $0$-geometric stack for the HAG context appearing \cite[2.3.2]{HAG2}
(such a stack is called a weakly 0-geometric $D$-stack, see \cite[1.3.3.1, 2.3.2.2]{HAG2} for these notions).

We also remark that (i) a scheme is generally not a derived algebraic stack
(for instance, consider the diagonal),
and (ii) there is another approach to derived geometry, which is
based on the theory of ringed $\infty$-topoi
\cite{DAG7}.
\end{Remark}

We fix convention for algebraic groups and their representations.

\begin{Definition}
\begin{enumerate}
\renewcommand{\labelenumi}{(\roman{enumi})}
\item
An affine group scheme $G$ over a field $k$
is a group object in the category of usual affine schemes over $k$.
If $G$ is an affine group scheme and we put $G=\Spec B$,
then $B$ is a commutative Hopf $k$-algebra.

\item By an algebraic group over $k$, we mean an affine group scheme of finite type
over $k$.
Every affine group scheme over $k$ is
a pro-algebraic group over $k$, that is,
a directed limit of algebraic groups.

\item A representation of an affine group scheme $G=\Spec B$ over $k$
is an (left or right) action of $G$ on a $k$-vector space $V$, that is
determined by the rule assigning 
to each $k$-algebra $R$ and $g\in G(R)$
an
isomorphism $\phi_g:V\otimes_kR\stackrel{\sim}{\to}V\otimes_kR$ of
$R$-modules in a coherently functorial fashion. Equivalently,
a representation is a coaction $V\to V\otimes_kB$
of the commutative Hopf algebra $B$ on $V$.
As is well-known, every
representation is a filtered colimit of finite-dimensional
representations. See e.g. \cite{DM} for the basic facts
on affine group schemes.
\end{enumerate}
\end{Definition}

Let $G$ be a usual affine group scheme over a field $k$ of characteristic zero.
Then since $k$ is a field,
it gives rise to a group object $D_G:\NNNN(\Delta)^{op}\to \Aff_k$
given by $[n]\mapsto G^{\times n}$.
We usually consider the usual affine group scheme
$G$ to be a group object in $\Aff_k$.
We denote by $BG$ the colimit of this group object
in $\Sh(\Aff_k)$ and refer to $BG$ as the classifying stack of $G$.

A derived stack $X$ is said to be
a quotient stack by an action of
$G$ if there exists a presentation $X_\bullet:\NNNN(\Delta)^{op}\to \Aff_k$
of $X$ and a natural transformation
$X_\bullet \to D_G$ such that for any $[m]\to [n]$, the diagram
\[
\xymatrix{
X_\bullet([n]) \ar[r] \ar[d] & X_\bullet([m]) \ar[d] \\
D_G([n]) \ar[r] & D_G([m])
}
\]
is a pullback square.
Put $\Spec A=X_\bullet([0])$.
In this case, we often write $X=[\Spec A/G]$ for the quotient stack.
(Here, we consider the
natural transformation $X_\bullet\to D_G$ to be an action of $G$ on $\Spec A$.)

When $G=\Spec B$ is an algebraic group over $k$, the quotient stack
$[\Spec A/G]$ is a derived algebraic stack.
To see this, it will suffice to check that
the first projection $\Spec A\times G\to \Spec A$
is a formally perfect morphism.
By the base change formula $L_{A\otimes B/A}\simeq (A\otimes B)\otimes_{B}L_{B/k}$
of cotangent complexes,
it is enough to observe that the structure morphism
$G\to \Spec k$ is formally perfect.
Since $k$ is of characteristic zero,
it follows from Cartier's theorem \cite[Exp. $\textup{VI}_B$ 1.6.1]{SGA3} that $G$ is smooth of finite type
over $k$.
The cotangent complex $L_{B/k}$ is 
equivalent to the dualizable usual $B$-module $\Omega_{B/k}$
of K\"ahler differentials, which is
placed in degree zero.
Hence $[\Spec A/G]$ is a derived algebraic stack.
Moreover, by \cite[1.2.8.3]{HAG2}, $G\to \Spec k$
is formally i-smooth relative to the HA context appearing
\cite[2.3.4]{HAG2} (see \cite[1.2.8.1]{HAG2} for the
formally i-smoothness). Since this property is stable under
base changes, the first projection $\Spec A\times G\to \Spec A$
is also formally i-smooth.
A direct consequence of
this observation and \cite[1.3.4.2]{HAG2}
is that $[\Spec A/G]$ is a $0$-geometric stack for the HAG context
appearing in \cite[2.3.4]{HAG2}, i.e., a $0$-geometric $D$-stack.

\subsection{Symmetric monoidal structure}
\label{QCC}
First, we briefly recall
the notion of symmetric monoidal $\infty$-categories.
Let $\xi_{n,i}:\langle n\rangle \to \langle 1\rangle$
be the map in $\FIN$ such that $\xi_{n,i}(j)$ is $1$ if $j=i$ and is $\ast$ if $j\neq i$.
A symmetric monoidal $\infty$-category is
defined to
be a coCartesian fibration $p:\CCC^\otimes\to \NNNN(\FIN)$
such that
\[
(\xi_{n,1})_*\times \ldots \times (\xi_{n,n})_*:\CCC_n\to \CCC_1\times \ldots\times \CCC_1
\]
is an equivalence for each $n\ge0$.
Here $\CCC_n:=p^{-1}(\langle n\rangle)$.
By convention, $\CCC_0\simeq \Delta^0$.
We refer to $\CCC_1$ as the underlying $\infty$-category (but
we usually denote by $\CCC$ the underlying $\infty$-category).
For ease of notation, we usually write $\CCC^\otimes$
for $\CCC^\otimes\to \NNNN(\FIN)$.
For two symmetric monoidal $\infty$-categories $p:\CCC^\otimes\to \NNNN(\FIN)$
and $q:\DDD^\otimes\to\NNNN(\FIN)$,
a symmetric monoidal functor $\CCC^\otimes\to \DDD^\otimes$
is a map of coCartesian fibrations $\CCC^\otimes \to \DDD^\otimes$
over $\NNNN(\FIN)$ which carries $p$-coCartesian edges to
$q$-coCartesian edges.

We say that an object $C$ in the underlying
$\infty$-category of a symmetric monoidal $\infty$-category
$\mathcal{C}^\otimes$
is dualizable if there exists an object
$C^\vee$ and two morphisms $e:C\otimes C^\vee\to \mathbf{1}$ and
$c:\mathbf{1} \to C\otimes C^\vee$
with $\mathbf{1}$ a unit such that the compositions
\[
C \stackrel{\textup{id}_C\otimes c}{\longrightarrow} C\otimes C^\vee\otimes C \stackrel{e\otimes\textup{id}_C}{\longrightarrow} C\ \ \textup{and}\ \ C^\vee \stackrel{c\otimes\textup{id}_{C^\vee}}{\longrightarrow} C^\vee \otimes C\otimes C^\vee \stackrel{\textup{id}_{C^\vee}\otimes e}{\longrightarrow} C^\vee 
\]
are equivalent to the identity of $C$ and the identity of $C^\vee$
respectively.
The symmetric monoidal structure of $\mathcal{C}$ induces
that of the homotopy category
$\textup{h}(\mathcal{C})$.
If we consider $C$ to be an object
in $\textup{h}(\mathcal{C})$,
then $C$ is dualizable in $\mathcal{C}$ if and only if $C$
is dualizable in $\textup{h}(\mathcal{C})$.

Let $\Cat^{\textup{Sym}}$ denote the $\infty$-category of
symmetric monoidal small $\infty$-categories,
that is obtained from the simplicial category of
symmetric monoidal $\infty$-categories (regarded as coCartesian 
fibrations) whose morphisms
are symmetric monoidal functors.
Using the straightening functor \cite[3.2]{HTT} and \cite[4.2.4.4]{HTT}
we have a fully faithful functor
\[
\Cat^{\textup{Sym}}\to \Fun(\NNNN(\FIN),\Cat).
\]
The essential image is spanned by commutative monoid objects
(i.e., $E_\infty$-monoid objects).
If we equip $\Cat$ with the symmetric monoidal structure
given by Cartesian product, then
a commutative monoid object amounts to a commutative algebra object.
Thus we have a natural categorical
equivalence $\Cat^{\textup{Sym}}\simeq \CAlg(\Cat)$.
We often think of a symmetric monoidal small $\infty$-category
as an object in $\CAlg(\Cat)$.

Let $\PR$ be the subcategory of $\widehat{\textup{Cat}}_{\infty}$
which consists of presentable $\infty$-categories and
whose edges (i.e. morphisms) are colimit-preserving functors.
The $\infty$-category $\PR$ has a symmetric monoidal structure
(see \cite[4.8.1.15, 4.8.1.17]{HA}).
For two presentable $\infty$-categories $\mathcal{C}$
and $\mathcal{D}$, the tensor product $\mathcal{C}\otimes \mathcal{D}$
is given
by $\Fun^{\textup{R}}(\mathcal{C}^{op},\mathcal{D})$,
where $\Fun^{\textup{R}}(-,-)$ denotes the full subcategory of
$\Fun(-,-)$ spanned by limit-preserving functors.
According to \cite[4.8.1.17]{HA} and the proof, the tensor product
$\mathcal{C}\otimes \mathcal{D}$ satisfies the following universal property:
it admits a functor $\mathcal{C}\times \mathcal{D}\to \mathcal{C}\otimes \mathcal{D}$ such that the composition induces a fully faithful functor
\[
\Map(\mathcal{C}\otimes \mathcal{D},\mathcal{E})\to \Map(\mathcal{C}\times \mathcal{D},\mathcal{E})
\]
whose essential image is spanned by functors which preserve
(small) colimits separately in each variable.
The $\infty$-category $\SSS$ of (small) spaces
is a unit object in $\PR$.
By \cite[4.8.1.19]{HA} and \cite[5.5.3.18]{HTT},
the tensor product $\PR\times \PR\to \PR$ preserves colimits
separately in each variable.

A presentably symmetric monoidal $\infty$-category $\mathcal{C}^\otimes$
(cf. Definition~\ref{presentably})
can be viewed as a commutative algebra object in the symmetric monoidal
$\infty$-category $(\PR)^{\otimes}$
in the same way that
a symmetric monoidal small $\infty$-category can
be viewed as an object in $\CAlg(\Cat)$.
Namely, $\mathcal{C}^\otimes$
belongs to $\CAlg(\PR)$.
A morphism in $\CAlg(\PR)$ corresponds to a symmetric monoidal functor
which preserves (small) colimits.
Let $R$ be a commutative ring spectrum and $\Mod_{R}^\otimes$
the symmetric monoidal (stable) $\infty$-category of $R$-module spectra.
Since $\Mod_R^\otimes$ lies in $\CAlg(\PR)$,
we can consider the symmetric monoidal $\infty$-category
$\Mod_{\Mod_R^\otimes}^\otimes(\PR)$ of $\Mod_R^\otimes$-module objects.
We write $\PR_R$ for $\Mod_{\Mod_R^\otimes}(\PR)$.
We shall refer to an object in $\PR_R$ as an $R$-linear
presentable $\infty$-category
and refer to $\PR_R$ as the $\infty$-category of
$R$-linear
presentable $\infty$-categories.
Similarly, we shall refer to an object in $\CAlg(\PR_R)$ as an $R$-linear
symmetric monoidal presentable $\infty$-category
and refer to $\CAlg(\PR_R)$ as the $\infty$-category of
$R$-linear
symmetric monoidal presentable $\infty$-categories.
A morphism in $\CAlg(\PR_R)$ will be referred to as
an $R$-linear symmetric monoidal functor.
Consider the case where $R$ is the sphere spectrum
$\mathbb{S}$. By \cite[4.8.2.18]{HA},
the forgetful functor
$\PR_{\mathbb{S}}\to \PR$ can be regarded as a fully faithful
embedding whose essential image consists of stable presentable $\infty$-categories (recall that $\mathbb{S}$ denotes the sphere spectrum). 
In particular, any $R$-linear presentable $\infty$-category is stable.
Let $\SP$ denote the stable presentable $\infty$-category of
spectra. We denote by $\otimes$ the smash product.
The left adjoint of $\PR_{\mathbb{S}}\to \PR$ is given
by $\PR\to \PR_{\mathbb{S}}$ which carries $\mathcal{C}$ to $\mathcal{C}\otimes \SP\simeq \varprojlim\Fun^{\textup{R}}(\mathcal{C}^{op},\mathcal{S}_*)$
where $\mathcal{S}_*=\mathcal{S}_{\Delta^0/}$ is the $\infty$-category of pointed spaces and the limit of the
sequence of the loop space functor
$\Omega_*:\mathcal{S}_*\to \mathcal{S}_*$ is taken in $\PR$.
If $R$ is the Eilenberg-MacLane spectrum $Hk$ for some ordinary
commutative ring $k$, then we write $\PR_k$ for $\PR_{Hk}$.
In that case, we use the term ``$k$-linear
presentable $\infty$-category'' instead of ``$Hk$-linear
presentable $\infty$-category''.
Recall that
the homotopy category $\textup{h}(\CCC)$
of a stable $\infty$-category $\CCC$
is a triangulated category (see \cite{HA}).
In particular, it is an additive category.
When $\CCC$ is a $k$-linear
presentable $\infty$-category, the additive
category $\textup{h}(\CCC)$ is $k$-linear; every hom
set $\Hom_{\textup{h}(\CCC)}(C,D)$
has the structure of a $k$-vector space, and the composition
$\Hom_{\textup{h}(\CCC)}(D,E)\times \Hom_{\textup{h}(\CCC)}(C,D)\to \Hom_{\textup{h}(\CCC)}(C,E)$ is $k$-bilinear.
The functor $\Mod_k\times \CCC\to \CCC$ given by
$\Mod_k^\otimes$-module
structure induces an action of $k=\Hom_{\textup{h}(\Mod_k)}(\uni_k,\uni_k)$ on $\Hom_{\textup{h}(\CCC)}(C,D)$, where $\uni_k$ is a unit in $\Mod_k$, and $C$ and $D$ belong to $\CCC$.
It gives rise to the structure of a $k$-vector space
\[
k\times \Hom_{\textup{h}(\CCC)}(C,D)= \Hom_{\textup{h}(\Mod_k)}(\uni_k,\uni_k)\times \Hom_{\textup{h}(\CCC)}(C,D)\to \Hom_{\textup{h}(\CCC)}(C,D)
\]
where the right map is determined by
$\Mod_k\times \CCC\to \CCC$. We easily see that the composition is $k$-bilinear.

\subsection{Quasi-coherent complexes}
\label{QCC}

Let $X:\Aff_k^{op}\to \widehat{\SSS}$ be a sheaf.
We will define the stable $\infty$-category of
quasi-coherent complexes on $X$ (cf. \cite[2.7, 2.7.9]{DAG8}).
The construction $\Spec A \mapsto \Mod_A^\otimes$
gives rise to a functor 
$\CAlg_{\mathbb{S}}\to \CAlg(\widehat{\textup{Cat}}_\infty)$.
By $\CAlg_k\simeq (\CAlg_{\mathbb{S}})_{k/}$,
it gives rise to $\QC^\otimes:\Aff_k^{op}=\CAlg_k\to \CAlg(\widehat{\textup{Cat}}_\infty)_{\Mod_k^\otimes/}$.
Since $\CAlg(\widehat{\textup{Cat}}_\infty)_{\Mod_k^\otimes/}$
admits large limits,
there is a right Kan extension of
$\QC^\otimes:\CAlg_k\to \CAlg(\widehat{\textup{Cat}}_\infty)_{\Mod_k^\otimes/}$:
\[
q:\Fun(\CAlg_k,\widehat{\SSS})^{op}\to \CAlg(\widehat{\textup{Cat}}_\infty)_{\Mod_k^\otimes/}
\]
which preserves (large) limits (see \cite[5.1.5.5]{HTT}
for the existence of a Kan extension).
By the flat hypercomplete descent property of modules
\cite[6.13]{DAG7},
the functor $q$ factors through the \'etale sheafification $\Fun(\CAlg_k,\widehat{\SSS})\to \Sh(\Aff_k)$ (moreover, it factors through the flat sheafification).
Consequently,
we obtain the induced functor
\[
\QC^\otimes:\Sh(\Aff_k)^{op}\to \CAlg(\widehat{\textup{Cat}}_\infty)_{\Mod_k^\otimes/}
\]
which carries $X$ to $\Mod_k^\otimes\to \QC^\otimes(X)$.
We define a symmetric monoidal stable $\infty$-category of quasi-coherent complexes on $X$ to be $\QC^\otimes(X)$.
Note that $\QC^\otimes(\Spec A)\simeq \Mod_A^\otimes$.
We also remark that the flat sheafification does not change $\QC(-)$.

Suppose that
$X$ is a derived stack over $k$ (cf. Definition~\ref{defderivedstack}).
Let $X_\bullet$ be a presentation of $X$.
Put $\Spec R^n=X_\bullet([n])$.
Let $X^\circ$ denote the colimit of $X_\bullet$ in $\Fun(\CAlg_k,\widehat{\SSS})$. Taking account of the right Kan extension $q$ and the sheafification
$X$ of $X^\circ$,
we have
\[
\varprojlim_{[n]\in\Delta} \Mod^\otimes_{R^n} \simeq q(X^\circ)\simeq \QC^\otimes(X).
\]
Thus, one may consider $\QC^\otimes(X)$ to be the limit
$\varprojlim \Mod^\otimes_{R^n}$.
In particular,
$\QC^\otimes(X)$ is a limit of the cosimplicial
diagram of presentably symmetric monoidal stable
$\infty$-categories.
Since a small limit in $\CAlg(\PR_{k})$ commutes with that in $\CAlg(\widehat{\textup{Cat}}_\infty)$,
$\QC^\otimes(X)$ belongs to
$\CAlg(\PR_k)$ for any derived stack $X$.
In particular, $\QC(X)$ is presentable.
We write $\OO_X$ for a unit object of the symmetric monoidal
$\infty$-category $\QC^\otimes(X)$.
For a morphism $f:X\to Y$ of derived stacks (that is, a morphism
as objects in $\Sh(\Aff_k)$), $\QC^\otimes$ induces
a morphism $f^*:\QC^\otimes (Y)\to \QC^\otimes (X)$ in $\CAlg(\PR_k)$
(note that any morphism $A\to B$ in $\CAlg_k$
induces the base change functor $\Mod_A^\otimes\to \Mod_B^\otimes$
that belongs to $\CAlg(\PR_k)$).
By the adjoint functor theorem, there is a right adjoint
$f_*:\QC(X)\to \QC(Y)$ of $f^*$. We shall refer to $f^*$ and $f_*$
as
the pullback functor and the pushforward functor, respectively.

\begin{Example}
\label{classifyingstackqc}
Let $G=\Spec B$ be an affine group scheme over a field $k$ of characteristic zero.
The associated group object $D_G:\NNNN(\Delta)^{op}\to \Aff_k$
can be regarded as a cosimplicial diagram $[n]\mapsto B^{\otimes n}$
of (usual) commutative $k$-algebras. Let $BG$ be the classifying stack.
Note that $BG$ is a colimit of $\NNNN(\Delta)^{op}\stackrel{D_G}{\to} \Aff_k\to \Sh(\Aff_k)$.
Thus, $\QC^\otimes(BG)$ can naturally
be identified with the limit $\varprojlim_{[n]\in \Delta} \Mod^\otimes_{B^{\otimes n}}$.
\end{Example}

\vspace{2mm}

{\it From model categories to $\infty$-categories.}
Here, we recall a version of Dwyer-Kan
localization in the context of $\infty$-categories
by which we can obtain
$\infty$-categories
from model categories (see \cite[1.3.4, 4.1.3]{HA}, \cite{Hin}).
Let $\mathbb{M}$ be a combinatorial model category (cf. \cite{HTT})
and $\mathbb{M}^c$ the full subcategory which consists of cofibrant
objects. Then there are an $\infty$-category $\NNNN_W(\mathbb{M}^c)$
and a functor $\xi:\textup{N}(\mathbb{M}^c)\to \NNNN_W(\mathbb{M}^c)$
such that for any $\infty$-category $\CCC$ the composition
induces a fully faithful functor
\[
\Map(\NNNN_W(\mathbb{M}^c),\CCC)\to \Map(\textup{N}(\mathbb{M}^c),\CCC)
\]
whose essential image consists of those functors $F:\textup{N}(\mathbb{M}^c)\to\CCC$ which carry weak equivalences in $\textup{N}(\mathbb{M}^c)$
to equivalences in $\CCC$.
By the Yoneda lemma, $\textup{N}(\mathbb{M}^c)\to \NNNN_W(\mathbb{M}^c)$
is unique up to a contractible space of choices.
We shall refer to $\NNNN_W(\mathbb{M}^c)$
as the $\infty$-category obtained from $\mathbb{M}$ (or $\mathbb{M}^c$)
by inverting weak equivalences.
An explicit construction of $\NNNN_W(\mathbb{M}^c)$
is given by the hammock localization (cf. \cite{DK}).
More precisely, one model
of $\NNNN_W(\mathbb{M}^c)$ is the simplicial nerve
of (a fibrant replacement of) the hammock localization of $\mathbb{M}^c$.
The homotopy category of $\NNNN_W(\mathbb{M}^c)$ coincides
with the homotopy category of the model category $\mathbb{M}$.
The $\infty$-category
$\NNNN_W(\mathbb{M}^c)$ is presentable. If $\mathbb{M}$ is
a stable model category,  $\NNNN_W(\mathbb{M}^c)$ is 
stable (cf. \cite{Tan}).
If $\mathbb{M}$ is a symmetric monoidal model
category, there are a
symmetric monoidal $\infty$-category $\NNNN_W^\otimes(\mathbb{M}^c)$
which belongs to $\CAlg(\PR)$ and a symmetric monoidal
colimit-preserving functor $\tilde{\xi}:\textup{N}^\otimes(\mathbb{M}^c)\to \NNNN_W^\otimes(\mathbb{M}^c)$ whose underlying functor is equivalent to $\xi$.
There is a universal property: for any symmetric
monoidal $\infty$-category $\CCC^\otimes$
the composition induces a fully faithful functor
$\Map_{\CAlg(\wCat)}(\NNNN_W^\otimes(\mathbb{M}^c),\CCC^\otimes)\to \Map_{\CAlg(\wCat)}(\textup{N}^\otimes(\mathbb{M}^c),\CCC^\otimes)$
whose essential image consists of those
functors $F:\textup{N}^\otimes(\mathbb{M}^c)\to\CCC^\otimes$ which carry weak equivalences in $\textup{N}(\mathbb{M}^c)$
to equivalences in $\CCC$.

Let us consider the model category
of chain complexes of representations.
Let $G$ be a pro-reductive group over a field $k$ of characteristic zero.
Let $\Vect(G)$ be the (symmetric monoidal) Grothendieck abelian category
of (not necessarily finite dimensional) representations
of $G$, that is, $k$-vector spaces equipped with
actions of $G$.
Let $\Comp(\Vect(G))$ be the symmetric monoidal
category
of (possibly unbounded) chain complexes
of objects in $\Vect(G)$.
Let $\mathcal{G}_{G}$ be the set of finite coproducts
of irreducible representations of $G$.
Let $\mathcal{H}=\{0\}$. Then
by the semi-simplicity of representations of $G$,
the pair $(\mathcal{G}_{G},\mathcal{H})$
is a flat descent structure in the sense of \cite{CD1}.
Consequently, there exists a combinatorial
symmetric monoidal model structure on
$\Comp(\Vect(G))$ such that (i) weak equivalences
are exactly quasi-isomorphisms, and (ii) coproducts
of objects in $\mathcal{G}$ are cofibrant \cite{CD1}.
Let $\DDD^{\otimes}(BG)$
denote the symmetric monoidal presentable $\infty$-category
obtained from the full subcategory $\Comp(\Vect(G))^c$
of cofibrant objects
by inverting weak equivalences.
Let $\Comp(k)\to \Comp(\Vect(G))$
be the natural left adjoint symmetric monoidal
functor which carries a unit $k$ to the trivial representation.
Inverting weak equivalences yields a symmetric monoidal
colimit-preserving functor $\Mod_k^\otimes\to \DDD^\otimes(BG)$.
Hence $\DDD^\otimes(BG)$ belongs to $\CAlg(\PR_k)\simeq \CAlg(\PR)_{\Mod_k^\otimes/}$.
There exists a natural equivalence
$\DDD^\otimes(BG)\simeq \QC^\otimes(BG)$; see e.g. \cite[Lemma 5.9]{PM}
(in {\it loc. cit.}, the reductive algebraic case is treated but that works
mutatis-mutandis in the case of pro-reductive groups).
We often write $\Rep^\otimes(G)$ for $\QC^\otimes(BG)$.

\vspace{3mm}

{\it Associated quotient stacks.}
Let $G=\Spec B$ be an affine group scheme over $k$.
Let $BG$ be the classifying stack of $G$.
Let $A$ be an object of $\CAlg(\QC(BG))$.
We will construct a quotient stack which we will denote by $[\Spec A/G]$.

Let $D_G:\NNNN(\Delta)^{op}\to \Aff_k$
be the group object corresponding to $G$.
For ease of notation, put $G_n:=D_G([n])\simeq G^{\times n}=\Spec B^{\otimes n}$.
Let $\pi:\Spec k\to BG$ be the canonical projection.
We often abuse notation by writing $A$
for the image of $A$ in $\CAlg_k$
under the pullback functor $\pi^*:\CAlg(\QC(BG))\to \CAlg(\QC(\Spec k))\simeq \CAlg_k$.
Let $\CAlg_{G_\bullet}\to \NNNN(\Delta)$
be a coCartesian fibration corresponding to
$[n]\mapsto \CAlg_{B^{\otimes n}}$ via the unstraightening functor (cf. \cite[3.2]{HTT}).
The limit $\varprojlim_{[n]\in \Delta}\CAlg_{B^{\otimes n}}$
can be identified with the full subcategory of 
$\Fun_{\NNNN(\Delta)}(\NNNN(\Delta),\CAlg_{G_\bullet})$
spanned by those functors which send all edges to
coCartesian edges.
The homomorphism $G\to \Spec k$ of group schemes
to the trivial group scheme $\Spec k$
gives rise to a map of coCartesian fibrations
\[
\xymatrix{
\CAlg_k\times \NNNN(\Delta) \ar[rr] \ar[rd]_{\textup{pr}_2} & &  \CAlg_{G_\bullet} \ar[dl]\\
  &  \NNNN(\Delta) & 
}
\]
such that each fiber
$\CAlg_k\to \CAlg_{B^{\otimes n}}$ is given by $A\mapsto A\otimes B^{\otimes n}$.
This map carries coCartesian edges to coCartesian edges. 
By the relative adjoint functor theorem \cite[7.3.2.7]{HA},
there is a right adjoint functor
$c:\CAlg_{G_\bullet}\to \CAlg_k\times \NNNN(\Delta)$
relative to $\NNNN(\Delta)$. (On each fiber,
$c$ induces $\CAlg_{B^{\otimes n}}\to \CAlg_{k}$
determined by the composition with $k\to B^{\otimes n}$,
that is the right adjoint of the base change $(-)\otimes_k B^{\otimes n}:\CAlg_k \to \CAlg_{B^{\otimes n}}$.)
The composition with
\[
\CAlg_{G_\bullet}\stackrel{c}{\to} \CAlg_k\times \NNNN(\Delta)\stackrel{\textup{pr}_1}{\to} \CAlg_k\hookrightarrow \Sh(\Aff_k)^{op}
\]
induces
\[
\Fun_{\NNNN(\Delta)}(\NNNN(\Delta),\CAlg_{G_\bullet})\stackrel{\eta}{\to} \Fun(\NNNN(\Delta),\CAlg_k)\to
\Fun(\NNNN(\Delta),\Sh(\Aff_k)^{op})\stackrel{\textup{lim}}{\to} \Sh(\Aff_k)^{op}.
\]
The right functor carries $f:\NNNN(\Delta)\to \Sh(\Aff_k)^{op}$ to its limit.
Note that it carries
an initial object in $\varprojlim_{[n]\in \Delta}\CAlg_{B^{\otimes n}}$
to $BG$.
We then obtain
$\varprojlim_{[n]\in \Delta}(\CAlg_{B^{\otimes n}})^{op}\to \Sh(\Aff_k)_{/BG}$.
By the construction, any $X\in \varprojlim_{[n]\in \Delta}(\CAlg_{B^{\otimes n}})^{op}\simeq \CAlg(\QC(BG))^{op}$ gives rise to $\eta(X)^{op}:\NNNN(\Delta)^{op}\to \Aff_k$
and a natural transformation $\eta(X)^{op}\to D_G$. 
We easily see that $\eta(X)^{op}\to D_G$ satisfies
the axiom of quotient stacks, i.e., the canonical morphism $\eta(X)^{op}([n])\to \eta(X)^{op}([m])\times_{D_G([m])}D_G([n])$ is an equivalence
for any map $[m]\to [n]$ of $\Delta$.
Now suppose that $X$ corresponds to $A\in \CAlg(\QC(BG))$.
We define $[\Spec A/G]$ to be a colimit of $\NNNN(\Delta)^{op}\stackrel{\eta(X)^{op}}{\to} \Aff_k\hookrightarrow \Sh(\Aff_k)$.
We shall refer to $[\Spec A/G]$ as the quotient stack associated to
$A\in \CAlg(\QC(BG))$.
The natural transformation $\eta(X)^{op}\to D_G$ induces a canonical morphism
$p:[\Spec A/G]\to BG$.
Informally, the quotient stacks and their presentations are depicted as follows:
\[
\xymatrix{
   \cdots \Spec A \times G\times G \ar[d] \ar@<5pt>[r] \ar@<-5pt>[r] \ar[r] &  \Spec A\times_k G \ar[d] \ar@<3pt>[r] \ar@<-3pt>[r]  \ar@<2.5pt>[l] \ar@<-2.5pt>[l] &  \Spec A \ar[l] \ar[d] \ar[r] &  [\Spec A/G] \ar[d]^p   \\
    \cdots G\times G \ar@<5pt>[r] \ar@<-5pt>[r] \ar[r] &   G   \ar@<2.5pt>[l] \ar@<-2.5pt>[l]   \ar@<3pt>[r] \ar@<-3pt>[r]      &   \Spec k \ar[l] \ar[r]^{\pi} & BG.
}
\]
The upper horizontal diagram is 
the augmented simplicial diagram induced by
$\eta(X)^{op}$,
and the lower horizontal diagram is the
augmented simplicial diagram induced by $D_G$.
Every square diagram is a pullback square.

In subsequent sections,
the following propositions will be useful:

\begin{Proposition}
\label{PMA}
Let $\CCC^\otimes$ and $\DDD^\otimes$ be 
objects in $\CAlg(\PR_\mathbb{S})$,
i.e., stable presentably symmetric monoidal
$\infty$-categories (cf. Definition~\ref{presentably}).
Let $F:\DDD^\otimes\to \CCC^\otimes$
be a symmetric monoidal functor which preserves small colimits.
Let $G:\CCC^\otimes\to \DDD^\otimes$ be a lax symmetric monoidal right adjoint functor of $F$ (which exists by the relative version of adjoint functor theorem \cite{HA}). Let $\uni_{\CCC}$ be a unit object of $\CCC$ (thus $\uni_\CCC\in \CAlg(\CCC)$)
and $B:=G(\uni_\CCC)\in \CAlg(\DDD)$.
Consider the composite of symmetric monoidal colimit-preserving
functors
\[
F':\Mod_B^\otimes(\DDD)\stackrel{F}{\to} \Mod_{F(B)}^\otimes(\CCC)\to \Mod_{\uni_\CCC}^\otimes(\CCC)\simeq \CCC^\otimes
\]
where the right functor is determined by the counit map $F\circ G(\uni_\CCC)\to \uni_\CCC$. 
Suppose that
\begin{enumerate}
\renewcommand{\labelenumi}{(\theenumi)}
\item There is a small set $\{I_\lambda\}_{\lambda\in \Lambda}$
of compact and dualizable objects of $\DDD$
which generates $\DDD$ as a stable presentable $\infty$-category,

\item Each $F(I_{\lambda})$ is compact, and
$\{F(I_\lambda)\}_{\lambda\in \Lambda}$ generates $\CCC$ as a stable
presentable $\infty$-category.

\end{enumerate}
Then $F'$ is an equivalence.
\end{Proposition}

If $F$ satisfies (1) and (2) in Proposition~\ref{PMA} we say that
$F$ is perfect.
Let $G':\calC^\otimes \to \Mod_B^\otimes(\calD^\otimes)$ be a lax symmetric monoidal functor which is a right adjoint
functor of $F'$. The existence of the right adjoint functor
follows from the relative version of adjoint functor theorem (see \cite[7.3.2.7]{HA}).
Therefore we have a diagram 
\[
\xymatrix{
\calD^\otimes \ar@<2pt>[r]^{F} \ar@<2pt>[d]^{R}  & \calC^\otimes \ar@<2pt>[l]^{G} \ar@<2pt>[ld]^{G'} \\
\Mod^\otimes_B(\calD^\otimes) \ar@<2pt>[u]^{U} \ar@<2pt>[ru]^(0.4){F'}.
}
\]
where $U$ is the forgetful functor
and $R$ assigns to $M\in \calD$ a free left $B$-module $B\otimes M$.
All functors are exact.
The composite $F'\circ R:\calD^\otimes \to \calC^\otimes$ is equivalent to
$F$ as symmetric monoidal functors.

\begin{Lemma}
\label{coptgen}
Suppose that
$\{I_\lambda\}_{\lambda\in \Lambda}$
is a small set of compact objects which generates
$\calD$ as a stable presentable $\infty$-category.
Then $\{R(I_\lambda)\}_{\lambda\in \Lambda}$ is a set of
compact objects which generates $\Mod_B(\calD^\otimes)$
as a stable presentable $\infty$-category.
\end{Lemma}

\Proof
We first show that $R(I_\lambda)$ is compact.
Let $\varinjlim N_i$ be a filtered colimit in $\Mod_B(\calD^\otimes)$.
Note that by \cite[4.2.3.5]{HA}
$U$ preserves colimits and thus $\varinjlim U(N_i)\simeq U(\varinjlim N_i)$.
Then we have natural equivalences
\begin{eqnarray*}
\Map_{\Mod_B(\calD^\otimes)}(R(I_\lambda),\varinjlim N_i) &\simeq& \Map_{\calD}(I_\lambda,U(\varinjlim N_i)) \\
&\simeq& \Map_{\calD}(I_\lambda,\varinjlim U(N_i)) \\
&\simeq& \varinjlim \Map_{\calD}(I_\lambda,U(N_i)) \\
&\simeq& \varinjlim \Map_{\Mod_B(\calD^\otimes)}(R(I_\lambda),N_i)
\end{eqnarray*}
in $\SSS$. Notice that the third equivalence
follows from the compactness of $I_\lambda$.
By these equivalences, we conclude that $R(I_\lambda)$ is compact.
It remains to prove that if 
$\Ext^n_{\Mod_B(\calD^\otimes)}(R(I_\lambda),N)=0$ for any $\lambda\in \Lambda$
and any integer $n\in \ZZ$, then
$N$ is zero.
Since
\[
\Ext^n_{\Mod_B(\calD^\otimes)}(R(I_\lambda),N)\simeq \Ext^n_\calD(I_\lambda,U(N))\simeq \Ext^n_\calD(I_\lambda,U(N))=0,
\]
our claim follows from the fact that $\{I_\lambda\}_{\lambda\in \Lambda}$ is a compact generator and $\Mod_B(\calD^\otimes)\to \calD$ is conservative.
\QED

{\it Proof of Proposition~\ref{PMA}.}
If $F'$ is fully faithful, $F'$ is also essentially
surjective.
In fact, if $F'$ is fully faithful, the essential image of $F'$ is
the smallest stable subcategory of $\calC$ which has colimits
and contains $F(I_\lambda)$ for all $\lambda\in \Lambda$.
By the condition (2),
the essential image of $F'$ coincides with
$\calC$.
Hence 
we will prove that $F'$ is fully faithful.
For this purpose, since $F'$ is an exact functor
between stable $\infty$-categories
$\Mod_B(\calD^\otimes)$ and $\calC$,
it will suffice to show, by \cite[Lemma 5.8]{Tan}, that
$F'$ induces a fully faithful functor between their homotopy categories.
We will prove that $F'$ induces a bijection
\[
\alpha:\Hom_{\Mod_{B}(\calD^\otimes)}(R(I_\lambda),R(\Sigma^n I_\mu))\to \Hom_{\calC}(F'(R(I_\lambda)),F'(R(\Sigma^nI_\mu)))
\]
where $\Hom(-,-)$ indicates $\pi_0(\Map(-,-))$ and $n$ is an integer.
Note that by adjunction, we have natural bijections
\begin{eqnarray*}
\Hom_{\Mod_{B}(\calD^\otimes)}(R(I_\lambda),R(\Sigma^n I_\mu))&\simeq& \Hom_\calD(I_\lambda,U(R(\Sigma^n I_\mu))) \\
&\simeq& \Hom_\calD(I_\lambda\otimes (\Sigma^n I_\mu)^{\vee}, G(\mathbf{1}_\calC)).
\end{eqnarray*}
Here $(\Sigma^n I_\mu)^\vee$ is the dual of $\Sigma^n I_\mu$.
On the other hand, we have natural bijections
\begin{eqnarray*}
\Hom_{\calC}(F'(R(I_\lambda)),F'(R(\Sigma^n I_\mu)))&\simeq& \Hom_\calC(F(I_\lambda),F(\Sigma^n I_\mu)) \\
&\simeq& \Hom_\calC(F(I_\lambda\otimes (\Sigma^nI_\mu)^{\vee}), \mathbf{1}_\calC).
\end{eqnarray*}
Also, by adjunction there is a bijection
\[
\beta:\Hom_\calD(I_\lambda\otimes (\Sigma^nI_\mu)^{\vee}, G(\mathbf{1}_\calC))
\to \Hom_\calC(F(I_\lambda\otimes (\Sigma^nI_\mu)^{\vee}), \mathbf{1}_\calC)
\]
which carries $f:I_\lambda\otimes (\Sigma^nI_\mu)^{\vee}\to G(\mathbf{1}_\calC)$
to $F(I_\lambda\otimes (\Sigma^nI_\mu)^{\vee})\stackrel{F(f)}{\longrightarrow} F(G(\mathbf{1}_\calC))\to \mathbf{1}_\calC$ where the second morphism is the counit map.
Therefore, it is enough to identify $\alpha$ with $\beta$
through the natural bijections.
Since $F'$ is symmetric monoidal,
by replacing $I_\lambda$ and $\Sigma^nI_\mu$ by $I_\lambda\otimes (\Sigma^nI_\mu)^{\vee}$
and $\mathbf{1}_\calD$ respectively, we may and will assume that
$\Sigma^nI_\mu=\mathbf{1}_\calD$.
According to the definition, $\alpha$ carries
$f:R(I_\lambda)=G(\mathbf{1}_{\calC})\otimes I_\lambda\to R(\mathbf{1}_{\calD})=B$
to
\[
\mathbf{1}_{\calC}\otimes_{F(G(\mathbf{1}_{\calC}))}F\circ U(f):\mathbf{1}_{\calC}\otimes_{F(G(\mathbf{1}_{\calC}))}F(G(\mathbf{1}_{\calC})\otimes I_\lambda)\to \mathbf{1}_{\calC}\otimes_{F(G(\mathbf{1}_{\calC}))} F(G(\mathbf{1}_{\calC})\otimes \mathbf{1}_{\calC})\simeq \mathbf{1}_{\calC}.
\]
Unwinding the definitions shows that $\beta$ sends
$f:R(I_\lambda)\to R(\mathbf{1}_{\calD})$
to the composite
\[
F(I_\lambda)\to F\circ U\circ R(I_\lambda)\to F\circ U\circ R(\mathbf{1}_{\calD})=F(G(\mathbf{1}_{\calC}))\to \mathbf{1}_{\calC}=F(\mathbf{1}_{\calD})
\]
where the first functor is induced by
the unit $\textup{id}\to U\circ R$,
the second functor is $F\circ U(f)$,
and the third functor is induced by
the counit $F\circ G\to \textup{id}$.
Note that $F(I_\lambda)\to F\circ U\circ R(I_\lambda)$ can be identified
with $\mathbf{1}_{\calC}\otimes_{F(G(\mathbf{1}_{\calC}))}F(G(\mathbf{1}_{\calC})\otimes I_\lambda)\to F(G(\mathbf{1}_{\calC})\otimes I_\lambda)$
induced by the unit
$\mathbf{1}_{\calC}\to F(G(\mathbf{1}_{\calC}))$
of $F(G(\mathbf{1}_{\calC}))\in \CAlg(\calC^\otimes)$.
Now the desired identification with $\beta$
follows from the fact that $F(I_\lambda)\to F\circ U\circ R(I_\lambda)\to \mathbf{1}_{\calC}\otimes_{F(G(\mathbf{1}_{\calC}))}F(G(\mathbf{1}_{\calC})\otimes I_\lambda)\simeq F(I_{\lambda})$ is the identity (note that
$\mathbf{1}_{\calC}\to F(G(\mathbf{1}_{\calC}))\to \mathbf{1}_{\calC}$
is the identity).

Next we then apply the bijection $\alpha$ to conclude that $F'$
is fully faithful.
Since $F'$ preserves colimits (in particular, it is exact), we see that
if $N, M\in \Mod_B(\calD^\otimes)$ belongs to
the smallest stable subcategory $\mathcal{E}$
which contains $\{R(I_\lambda)\}_{\lambda\in \Lambda}$, then $F'$ induces a bijection
\[
\Hom_{\Mod_{B}(\calD^\otimes)}(N,M)\to \Hom_\calC(F'(N),F'(M)).
\]
There is a categorical equivalence
$\Ind(\mathcal{E})\simeq \Mod_B(\calD^\otimes)$ which follows
from Lemma~\ref{coptgen} and \cite[5.3.5.11]{HTT}.
Again by \cite[5.3.5.11]{HTT} and
the fact that
$F'(E)$ is compact for any $E\in \mathcal{E}$ (by condition (2)),
a left Kan extension $\Ind(\mathcal{E})\to \calC$
induced by $F':\mathcal{E}\to \calC$
(cf. \cite[5.3.5.10]{HTT}, \cite[4.8.1.14]{HA})
is fully faithful.
This implies that $F'$ is fully faithful.
\QED

\begin{Proposition}
\label{affinecor}
Let $G$ be an affine group scheme over a field $k$ of characteristic zero.
Let $A$ be an object of $\CAlg(\QC(BG))$,
and
let $[\Spec A/G]$ be the quotient stack associated to $A$.
Then there exists a canonical equivalence of symmetric monoidal $\infty$-categories
\[
\Mod_{A}^\otimes(\QC(BG))\to \QC^\otimes([\Spec A/G]).
\]
\end{Proposition}

\Proof
We use
the notation in the discussion before Proposition~\ref{PMA}.
Let $\eta(X)^{op}([n])=\Spec C^n \simeq \Spec A\otimes B^{\otimes n}$
and remember $D_G([n])\simeq G^{\times n}=\Spec B^{\otimes n}$.
We may view $C^n$ as a commutative
algebra over $B^{\otimes n}$,
 which is the image of $A\in \CAlg(\QC(BG))\simeq \varprojlim_{[n]\in \Delta}\CAlg_{B^{\otimes n}}$
 under the pullback functor $\CAlg(\QC(BG))\to \CAlg(\QC(D_G([n])))=\CAlg_{B^{\otimes n}}$.
By the definition of $[\Spec A/G]$,
$\QC^\otimes([\Spec A/G])$ can naturally be identified with
the limit $\varprojlim_{[n]\in \Delta}\Mod^\otimes_{C^n}$.
On the other hand, there are canonical equivalences
$\Mod_{A}^\otimes(\QC(BG))\simeq \varprojlim_{[n]\in \Delta}\Mod_{C^n}^\otimes(\Mod_{B^{\otimes n}})\simeq \varprojlim_{[n]\in \Delta}\Mod_{C^n}^\otimes$.
Thus we have $\Mod_{A}^\otimes(\QC(BG))\simeq \QC^\otimes([\Spec A/G])$.
\QED

\section{A universal characterization of representations of general linear groups}
\label{charsection}

Throughout this Section, $k$ is a field of characteristic zero.
Let $\CCC^\otimes$ be
a $k$-linear symmetric monoidal presentable $\infty$-category.
Let $\mathcal{C}_{\wedge,d}$ denote the full subcategory of $d$-dimensional wedge-finite objects in $\CCC$, and let $\mathcal{C}_{\wedge,d}^\simeq$ be the
largest Kan subcomplex.
Let $\Rep^\otimes(\GL_d)=\QC^\otimes(B\GL_d)$.
The main purpose of this Section is to prove the following result:

\begin{Theorem}
\label{characterization}
Let $\mathcal{C}^\otimes$ be a $k$-linear symmetric monoidal
(stable) presentable
$\infty$-category, i.e., an object in $\CAlg(\PR_k)$.
Then there exists a natural homotopy equivalence of spaces
\[
\Map_{\CAlg(\textup{Pr}_k^{\textup{L}})}(\Rep^\otimes(\GL_d), \mathcal{C}^\otimes) \to \mathcal{C}_{\wedge,d}^{\simeq}
\]
which carries $f:\Rep^\otimes(\GL_d) \to \CCC^\otimes$
to the image $f(K)$ of the standard representation $K$ of $\GL_d$.
That is, an object $C\in \mathcal{C}_{\wedge,d}$ corresponds to
a $k$-linear symmetric monoidal functor
$\Rep^\otimes(\GL_d)\to \mathcal{C}^\otimes$ that sends $K$ to $C$.
\end{Theorem}

\begin{Remark}
\label{wedgedual2}
By Theorem~\ref{characterization},
every wedge-finite object is the image of the standard representation
of $\GL_d$ for some $d\ge0$ under a symmetric monoidal functor.
The standard representation $K$ is dualizable in $\Rep^\otimes(\GL_d)$,
and any symmetric monoidal exact functor preserves dualizable objects.
Hence every wedge-finite object is dualizable.
\end{Remark}

\begin{Remark}
We use the assumption that the field $k$ is of characteristic zero
in an essential way.
\end{Remark}

\vspace{3mm}

We define a category $B\Sigma$ as follows:
Objects of $B\Sigma$ are finite sets, that is, $\bar{0}, \bar{1},\ldots , \bar{n}=\{1,\ldots,n\},\ldots$.
By convention $\bar{0}$ is the empty set.
A morphism
in $B\Sigma$ is a bijective map $\bar{n} \to \bar{n}$.
Namely, $\Hom_{B\Sigma}(\bar{n}, \bar{n})$ is isomorphic
to the symmetric group $\Sigma_n$ for $n\ge 0$, where $\Sigma_0$ is the group
consisting of one element.
If $n\neq m$, $\Hom_{B\Sigma}(\bar{n}, \bar{m})$
is the empty set.
Thus $B\Sigma$ is the coproduct $\sqcup_{n\ge0}B\Sigma_n$ (in $\Cat$)
where $B\Sigma_n$ is the category consisting of one object
$\bar{n}$ (regarded as a formal symbol)
such that $\Hom_{B\Sigma_n}(\bar{n}, \bar{n})=\Sigma_n$.
Let $\Vect_k$ be the category of $k$-vector spaces.
Here we denote by $\Fun(B\Sigma^{op},\Vect_k)$ 
the functor category.
It is a Grothendieck abelian category; it is presentable
(cf. \cite[5.5.3.6]{HTT}) and monomorphisms are closed under
filtered colimits.

Given an abelian category $\mathcal{A}$, we write $\Comp(\mathcal{A})$
for the category of chain complexes of objects in $\mathcal{A}$.
The category $\Comp(\Fun(B\Sigma^{op},\Vect_k))$
is isomorphic
to the functor category $\Fun(B\Sigma^{op},\Comp(k))$.
Here for ease of notation we write $\Comp(k)$ for $\Comp(\Vect_k)$.
An object $E:B\Sigma^{op}\to \Comp(k)$
corresponds to a symmetric sequence in the sense of \cite[Section 6]{Hov2}, that is,
\[
(E_0,E_1,\ldots,E_n,\ldots)
\]
where each chain complex $E_n=E(\bar{n})$
is endowed with a right $\Sigma_n$-action.
Recall from \cite{Hov2}
the symmetric monoidal structure on $\Fun(B\Sigma^{op},\Comp(k))$
defined as follows:
the tensor product $E\otimes F$ for $E,F\in \Fun(B\Sigma^{op},\Comp(k))$ is given by
\[
\bar{l} \mapsto \coprod_{A\sqcup B=\bar{l}, A\cap B=\phi}E(A)\otimes F(B)
\]
which is $\Sigma_l$-equivariantly isomorphic to
\[
\bar{l} \mapsto 
\coprod_{n+m=l}E(\bar{n})\otimes F(\bar{m}) \otimes_{k[\Sigma_n]\otimes_kk[\Sigma_m]}k[\Sigma_l]
\]
on which $\Sigma_l$ acts by the right multiplication.
Here for a finite group $G$, $k[G]$ denotes the group algebra,
and $E(\bar{n})\otimes F(\bar{m})$
is considered to be a right $k[\Sigma_n]\otimes_kk[\Sigma_m]$-module, and 
$k[\Sigma_l]$ is considered to be a left $k[\Sigma_n\times\Sigma_m]$-module
through the natural inclusion $\Sigma_n\times \Sigma_m\subset \Sigma_l$.
For any $a\ge 0$, we define a symmetric sequence
$I^a=(I^a_n)_{n\ge0}$
by
$I_a^a=k[\Sigma_a]$ equipped with the right multiplication of $\Sigma_a$,
and $I_n^a=0$ for $n\neq a$.
Then for any $a,b\ge0$, the tensor product
$I^a\otimes I^b$ is $I^{a+b}$,
and the commutative constraint on $I^{a+b}_{a+b}=k[\Sigma_{a+b}]$
is defined by the left action of the permutation
$(1,\ldots,a,a+1,\ldots,a+b)\mapsto (a+1,\ldots,a+b,1,\ldots, a)$.

By using the machinery in \cite{CD1},
we equip $\Fun(B\Sigma^{op},\Comp(k))$
with a combinatorial symmetric monoidal model structure.
The class of weak equivalences consists of termwise quasi-isomorphisms.
Let $\mathcal{G}$ be the set of finite coproducts of objects
in $\Fun(B\Sigma^{op},\Vect_k)$
which have
the form $(E_n)_{n\ge0}$ such that there is a non-negative integer
$i$ such that $E_i$ is an irreducible
$k$-linear $\Sigma_n$-representation, and
$E_n=0$ if $n \neq i$.
Set $\mathcal{H}={0}$.
Then by the representation theory
of symmetric groups in characteristic zero
and its semi-simplicity,
we see that
the pair $(\mathcal{G},\mathcal{H})$ is a flat descent structure
in the sense of \cite{CD1}.
According to \cite[Theorem 2.5, Proposition 3.2]{CD1},
there is a proper combinatorial symmetric monoidal model structure on 
$\Fun(B\Sigma^{op},\Comp(k))$ in which weak equivalences
are termwise quasi-isomorphisms (we do not recall the cofibrations
and fibrations, see \cite{CD1}).

Let $\DDD^\otimes(B\Sigma,k):=\NNNN_W(\Fun(B\Sigma^{op},\Comp(k))^c)$
be the symmetric monoidal presentable stable $\infty$-category
obtained from the full subcategory
$\Fun(B\Sigma^{op},\Comp(k))^c$ of cofibrant objects
by inverting weak equivalences.
The $\infty$-category $\DDD(B\Sigma,k)$
is equivalent to $\NNNN_W(\Fun(B\Sigma^{op},\Comp(k)))$.
By \cite[1.3.4.25]{HA} and $B\Sigma\simeq \sqcup_{n\ge0} B\Sigma_n$, there exist equivalences of $\infty$-categories
\[
\DDD(B\Sigma,k)\simeq \Fun(B\Sigma^{op},\Mod_k)\simeq \prod_{n\ge0}\Fun(B\Sigma_n^{op},\Mod_k),
\]
where $\Fun(-,-)$ denotes the function complex.
Here we abuse notation by indicating with $B\Sigma$
the nerve of $B\Sigma$.

Let us consider the functor category $\Fun(B\Sigma_n^{op},\Comp(k))$,
which we often identify with the category of chain complexes
of $k$-linear representations, that is, $k$-vector spaces
endowed with right actions of $\Sigma_n$.
As in the case of $\Fun(B\Sigma^{op},\Comp(k))$, $\Fun(B\Sigma_n^{op},\Comp(k))$ admits a combinatorial model structure
in which weak equivalences are exactly quasi-isomorphisms.
Let $\DDD(B\Sigma_n,k)$
be the presentable $\infty$-category
obtained from
$\Fun(B\Sigma_n^{op},\Comp(k))^c$
by inverting weak equivalences. The $\infty$-category
$\DDD(B\Sigma_n,k)$ can be identified with
$\Fun(B\Sigma_n^{op},\Mod_k)$. The homotopy category of $\DDD(B\Sigma_n,k)$ is
the (unbounded) derived category of $k$-linear representations
of $\Sigma_n$.
Since $\DDD(B\Sigma,k)$ is the product of $\{\DDD(B\Sigma_n,k)\}_{n\ge0}$, we often write $(E_n)_{n\ge0}$ with $E_n\in \DDD(B\Sigma_n,k)$
for an object in $\DDD(B\Sigma,k)$.
If we regard $E_n \in \DDD(B\Sigma_n,k)$ as an object in 
$\DDD(B\Sigma,k)$ in the obvious way,
the coproduct $\oplus_{n\ge0}E_n$ in $\DDD(B\Sigma,k)$
is $(E_n)_{n\ge0}$ since
\[
\Map(\oplus_{n\ge0}E_n,(F_n)_{n\ge0})\simeq \prod_{n\ge0}\Map(E_n,(F_n)_{n\ge0})\simeq \prod_{n\ge0}\Map(E_n,F_n)=\Map((E_n)_{n\ge0},(F_n)_{n\ge0})
\]
where we omit the subscript in each $\Map(-,-)$.

Next we construct a natural symmetric monoidal functor $\Mod_k^\otimes\to \DDD^\otimes(B\Sigma,k)$. For this,
let us consider a symmetric monoidal functor
$p:\Comp(k)\to \Fun(B\Sigma^{op},\Comp(k))$
given by $V$ to $p(V)=(V, 0,0,\ldots)$.
It is a left adjoint functor; the right adjoint is
determined by evaluation at the $0$-th term
$(E_0,E_1,\ldots)\mapsto E_0$.
There is a combinatorial symmetric monoidal
model structure on $\Comp(k)$
in which (i) weak equivalences are quasi-isomorphisms,
(ii) cofibrations are degreewise monomorphisms, and (iii)
fibrations are degreewise epimorphisms (cf. \cite[7.1.2.8, 7.1.2.11]{HA},
\cite[2.3.11]{Ho1}).
We remark that
generating cofibrations in $\Comp(k)$ (described in the proof of \cite[7.1.2.8]{HA} or \cite[2.3.3]{Ho1}) map to cofibrations in $\Fun(B\Sigma^{op},\Comp(k))$.
Hence $p$ preserves cofibrations and weak equivalences.
In particular, $p$ is a left Quillen adjoint functor.
Inverting weak equivalences of full subcategories of cofibrant objects
gives a symmetric monoidal colimit-preserving functor
\[
\NNNN_W(\Comp(k)^c)\longrightarrow \NNNN_W(\Fun(B\Sigma^{op},\Comp(k))^c)= \DDD(B\Sigma,k).
\]
Thanks to \cite[7.1.2.13]{HA}, there is a natural
symmetric monoidal equivalence
$\NNNN_W(\Comp(k)^c)\simeq \Mod_k$.
Hence we obtain a symmetric monoidal colimit-preserving functor 
$\Mod_k^\otimes\to \DDD^\otimes(B\Sigma,k)$.
Thus $\Mod_k^\otimes\to \DDD^\otimes(B\Sigma,k)$
belongs to $\CAlg(\PR_k)\simeq \CAlg(\PR)_{\Mod_k^\otimes/}$.

\begin{Lemma}
\label{linearleftkan}
Let $S$ be a small $\infty$-category. Let $R$ be a commutative ring
spectrum.
Let
\[
S\to \Fun(S^{op},\mathcal{S})\to \Fun(S^{op},\mathcal{S}_*)\to
\Fun(S^{op},\SP)\to \Fun(S^{op},\Mod_R)
\]
be the sequence of functors; the first functor is the Yoneda embedding,
the other functors are determined by the composition with
$\SSS\to \SSS_* \stackrel{\Sigma^\infty}{\to} \SP \stackrel{\otimes_{\mathbb{S}}R}{\to} \Mod_R$ where $\SSS\to \SSS_*$ carries $A$ to $A\sqcup \Delta^0$.
Let $\CCC$ be an $R$-linear presentable $\infty$-category,
that is, an object in $\PR_R$.
Then $\Fun(S^{op},\Mod_R)\simeq \Fun(S^{op},\SSS)\otimes \Mod_R$, and
the composition with the composite $S\to \Fun(S^{op},\Mod_R)$ induces
a homotopy equivalence
\[
\Map_{\PR_R}(\Fun(S^{op},\Mod_R),\mathcal{C})\to \Map(S,\CCC).
\]
\end{Lemma}

\Proof
By the left Kan extension (cf. \cite[5.1.5.6]{HTT}),
the Yoneda embedding induces
\[
\Map_{\PR}(\Fun(S^{op},\SSS),\CCC)\simeq \Map(S,\CCC)
\]
for any $\CCC\in \PR$.
Consider the adjoint pair $\PR\rightleftarrows \PR_R$
where the right adjoint $\PR_R\to \PR$ is the forgetful functor,
and the left adjoint is given by the base change $(-)\otimes \Mod_R$.
Taking account of this adjoint pair $\PR\rightleftarrows \PR_R$ we have
\[
\Map_{\PR}(\Fun(S^{op},\SSS),\CCC)\simeq \Map_{\PR_R}(\Fun(S^{op},\SSS)\otimes\Mod_R,\CCC)
\]
for any $\CCC\in \PR_R$.
Next we show that $\Fun(S^{op},\SSS)\otimes\DDD\simeq \Fun(S^{op},\DDD)$
for any $\DDD\in \PR$.
By definition,
\[
\Fun(S^{op}, \SSS)\otimes \DDD\simeq \Fun^{\textup{R}}(\DDD^{op},\Fun(S^{op},\SSS))\simeq \Fun'(\DDD^{op}\times S^{op},\SSS)
\]
where $\Fun'(\DDD^{op}\times S^{op},\SSS)$ denotes
the full subcategory of $\Fun(\DDD^{op}\times S^{op},\SSS)$ spanned by
those functors which preserve small limits in the variable $\DDD^{op}$.
There exist equivalences
\[
\Fun'(\DDD^{op}\times S^{op},\SSS)\simeq \Fun(S^{op},\Fun^{\textup{R}}(\DDD^{op},\SSS))\simeq \Fun(S^{op},\DDD\otimes \SSS)\simeq \Fun(S^{op},\DDD).
\]
Thus $S\to \Fun(S^{op},\Mod_R)$ induces the desired equivalence.
\QED

\begin{Lemma}
\label{coproduct}
There is a natural equivalence
\[
\Fun(B\Sigma^{op},\Mod_R)\simeq \bigoplus_{n\ge0} \Fun(B\Sigma_n^{op},\Mod_R)
\]
in $\PR_R$.
Here the coproduct $\bigoplus_{n\ge0}$ of the right-hand side is taken in $\PR_R$.
\end{Lemma}

\Proof
Invoking Lemma~\ref{linearleftkan},
we have
\begin{eqnarray*}
\Map_{\PR_R}(\Fun(B\Sigma^{op},\Mod_R),\CCC) &\simeq& \Map(B\Sigma,\CCC) \\
&\simeq& \prod_{n\ge0}\Map(B\Sigma_n,\CCC) \\
&\simeq& \prod_{n\ge0}\Map_{\PR_R}(\Fun(B\Sigma_n^{op},\Mod_R),\CCC) \\
&\simeq& \Map_{\PR_R}(\bigoplus_{n\ge0}\Fun(B\Sigma_n^{op},\Mod_R),\CCC)
\end{eqnarray*}
for any $\CCC\in \PR_R$. This proves our assertion.
\QED

\begin{Remark}
In $\widehat{\textup{Cat}}_\infty$, the $\infty$-category
$\Fun(B\Sigma^{op},\Mod_R)$
is not a coproduct of $\{\Fun(B\Sigma_n^{op},\Mod_R)\}_{n\ge0}$
but a product
$\prod_{n\ge0}\Fun(B\Sigma_n^{op},\Mod_R)$.
\end{Remark}

\begin{Proposition}
\label{absolute1}
Suppose that $\CCC^\otimes$ belongs to $\CAlg(\PR_k)$.
There exists a natural equivalence
\[
\Map_{\CAlg(\PR_k)}(\DDD^\otimes(B\Sigma,k),\CCC^\otimes)\simeq \CCC^\simeq.
\]
\end{Proposition}

To prove Proposition~\ref{absolute1},
we first recall the notion of free commutative algebra objects in (general)
symmetric monoidal $\infty$-categories (cf. \cite[3.1]{HA}).
Let $\CCC^\otimes$ be a symmetric monoidal $\infty$-category
and $\CAlg(\CCC)$ the $\infty$-category of commutative algebra objects.
We denote by $\theta:\CAlg(\CCC)\to \CCC$ the forgetful functor.
For $C\in \CCC$, $A\in \CAlg(\CCC)$ and $\phi:C\to \theta(A)$,
we say that $\phi$ makes $A$ a free commutative algebra object generated by $C$
if
$\Map_{\CAlg(\CCC)}(A,B)\to \Map_{\CCC}(C,\theta(B))$,
informally given by $f\mapsto \theta(f)\circ \phi$,
is a homotopy equivalence.
If we suppose that $\CCC$ admits countable colimits
and the tensor product preserves countable colimits separately in each
variable,
then $\theta$ has a left adjoint $\textup{Free}_{\CCC}:\CCC\to \CAlg(\CCC)$,
so that $(\textup{Free}_\CCC(C), C\to \textup{Free}_\CCC(C))$ is a free commutative
algebra object generated by $C$ where $C\to \textup{Free}_\CCC(C)$ is
the unit map
determined by the adjoint pair.

Consider the free commutative symmetric monoidal $\infty$-category
$\textup{Free}(\Delta^0)$ generated by the ``trivial'' category $\Delta^0$.
More precisely, $\textup{Free}(\Delta^0)$
is the image of $\Delta^0$ under the left adjoint functor $\textup{Free}$ in
\[
\textup{Free}:\Cat \rightleftarrows \CAlg(\Cat):\theta=\textup{forget}.
\]
The free algebra object $\textup{Free}(\Delta^0)$ has a more explicit form $B\Sigma$.
We define a (strict) symmetric monoidal structure on $B\Sigma$.
The tensor product $\otimes:B\Sigma\times B\Sigma\to B\Sigma$
is given by $\bar{n}\otimes \bar{m}:=\overline{n+m}$.
A pair of maps
$\phi:\bar{n}\to \bar{n}$ and $\psi:\bar{m}\to \bar{m}$
induces the map
$\phi\otimes \psi:\overline{n+m}\to \overline{n+m}$
determined by the permutations of $\{1,\ldots,n\}$ and $\{n+1,\ldots,n+m\}$
given by $\phi$ and $\psi$ respectively.
The commutative constraint
$\bar{n}\otimes \bar{m}=\overline{n+m} \to \overline{n+m}=\bar{m}\otimes \bar{n}$ is given by the left multiplication by the permutation $(1,\ldots,n,n+1,\ldots,n+m)\mapsto (n+1,\ldots,n+m,1,\ldots,n)$.
The unit object is $\bar{0}$.

\begin{Proposition}
\label{freebasic}
Let $v:\Delta^0\to B\Sigma$ be a functor determined by the value $\bar{1}$.
Then a pair $(B\Sigma,v:\Delta^0\to B\Sigma)$
is a free commutative algebra object in $\Cat$, generated by $\Delta^0$.
In particular, there exists a symmetric monoidal equivalence
$\textup{Free}(\Delta^0)\simeq B\Sigma$.
\end{Proposition}

Before giving the proof of Proposition~\ref{freebasic},
we recall the symmetric powers $\Sym^*(-)$.
Our main reference is \cite[3.1.3]{HA}
but we prefer to use a simpler formulation
which is described in its former version \cite[Section 3]{DAG3}
(note that \cite{DAG3} is not the newest version on arXiv but the version 3).
Let $\FIN^{\sim}$ be the subcategory of $\FIN$ such that
(i) objects in $\FIN^{\sim}$ are same with $\FIN$,
and (ii) a morphism of $\FIN$ lies in $\FIN^\sim$ if and only if
it is an isomorphism.
Notice that $\NNNN(\FIN^\sim)\simeq B\Sigma$.
For an $\infty$-category $\CCC$ we refer to a
 functor
 $\NNNN(\FIN^\sim)\to \CCC$
 as a symmetric sequence in $\CCC$.
Roughly speaking, a symmetric sequence in $\CCC$
consists of data $\{C_n\}_{n\ge0}$ where each $C_n$
is endowed with the left action of $\Sigma_n$.
As constructed in \cite[Section 3]{DAG3} for any symmetric monoidal
$\infty$-category
there is a functor $\textup{PSym}:\CCC\to \Fun(\NNNN(\FIN^\sim),\CCC)$
which sends $C$ to $\{C^{\otimes n}\}_{n\ge0}$
such that each $C^{\otimes n}$ is equipped with the permutation action
of $\Sigma_n$.
Suppose that $\CCC$ has countable colimits.
We define $\Sym^*:\CCC\to \CCC$
to be the composite
\[
\CCC\stackrel{\textup{PSym}}{\longrightarrow} \Fun(\NNNN(\FIN^\sim),\CCC)\to \CCC
\]
where the right functor carries $\NNNN(\FIN^\sim)\to \CCC$ to the colimit.
If $\FIN^\sim(n)$ is the full subcategory of $\FIN^\sim$
spanned by $\langle n\rangle$, then we define $\Sym^n$ to be the composite
\[
\CCC\stackrel{\textup{PSym}}{\to} \Fun(\NNNN(\FIN^\sim),\CCC)\to \Fun(\NNNN(\FIN^\sim(n)),\CCC) \to \CCC
\]
where the middle functor is induced by the restriction and the right functor
carries diagrams to colimits.
By definition, $\Sym^nC$ is a colimit of the permutation
action of $\Sigma_n$ on $C^{\otimes n}$, and
$\Sym^*C$ is the coproduct $\sqcup_{n\ge0}\Sym^nC$.

\vspace{2mm}

\Proof
We apply \cite[3.12]{DAG3} (cf. \cite[3.1.3.13 (ii)]{HA}) to our situation:
$B\Sigma$ is 
a free commutative algebra object generated by
$\Delta^0$
if and only if
the composite
$\Sym^*(\Delta^0)\stackrel{\Sym(v)}{\to} \Sym^*(B\Sigma)\to B\Sigma$
is a categorical equivalence.
Here $\Sym^*(B\Sigma)\simeq \sqcup_{n\ge0}\Sym^n(B\Sigma)\to B\Sigma$ is induced by the evaluation of
the natural transformation $\textup{PSym}(B\Sigma)\to B\Sigma_{\NNNN(\FIN^\sim)}$
from $\textup{PSym}(B\Sigma):\NNNN(\FIN^\sim) \to \Cat$ to the constant functor $B\Sigma_{\NNNN(\FIN^\sim)}:\NNNN(\FIN^\sim)\to \Cat$ taking the value $B\Sigma$ (see \cite[3.10]{DAG3}).
In concrete terms, for each $n\ge 0$
the 
evaluation at $\langle n\rangle $ induces
the $n$-fold tensor product
$B\Sigma^{\times n}\to B\Sigma$
which factors through the projection
$B\Sigma^{\times n}\to \Sym^n(B\Sigma)$.
To prove that the composite is an equivalence, note first
that
\[
B\Sigma^{\times n}=\sqcup_{(r_1,\ldots,r_n)}B\Sigma_{r_1}\times \ldots \times B\Sigma_{r_n}
\]
since the cartesian product in $\Cat$ preserves colimits separately
in each variable.
Hence $B\Sigma_{1}\times \ldots \times B\Sigma_{1}$
is a direct summand of $B\Sigma^{\times n}$
which is compatible with the permutation (left) action of $\Sigma_n$.
Note that the action of $\Sigma_n$ on
$B\Sigma_{1}\times \ldots \times B\Sigma_{1}$ is trivial
since $B\Sigma_1$ is contractible.
We have the following diagram:

\[
\xymatrix{
(\Delta^0)^{\times n} \ar[r]^(0.5){\sim} \ar[d] & B\Sigma_1^{\times n} \ar[r]^f \ar[d] & B\Sigma_n \\
\Sym^n(\Delta^0)=B\Sigma_n \ar[r]^(0.6){\sim} & B\Sigma_1^{\times n}/\Sigma_n. \ar[ur]_g & 
}
\]
The vertical functors are natural projections.
The functor $f$ is induced by the $n$-fold
tensor product $B\Sigma^{\times n}\to B\Sigma$.
By the commutative constraint of the symmetric monoidal
structure of $B\Sigma$, $f$ factors through
the projection $B\Sigma_1^{\times n}\to B\Sigma_1^{\times n}/\Sigma_n$,
which gives rise to $g$.
Here we consider $B\Sigma_1^{\times n}/\Sigma_n$ as
a direct summand of $\Sym^n(B\Sigma)$, and $g$ is
$B\Sigma_1^{\times n}/\Sigma_n
 \hookrightarrow \Sym^n(B\Sigma)\to B\Sigma$.
The lower horizontal functor is induced by $\Sym^*(v)$.
It will suffice to show that $g$ is a categorical equivalence.
The functor $g$ is determined by $f$.
More precisely,
we think of $f$ as a morphism in $\Fun(B\Sigma_n,\Cat)$, i.e.,
a natural transformation from the constant functor $B\Sigma_n\to \Cat$ taking the value $B\Sigma_1^{\times n}$ to the constant functor $B\Sigma_n \to \Cat$
taking the value $B\Sigma_n$.
Note that for any group $G$
there is an adjoint pair
\[
\alpha:\Fun(BG,\Cat)\rightleftarrows \Cat:\delta
\]
where the right adjoint $\delta:\Cat\to \Fun(BG,\Cat)$ is the diagonal
embedding by the
composition with $BG\to \Delta^0$.
The left adjoint carries $BG\to \Cat$ to its colimit.
Through this adjoint pair the morphism $f$ in $\Fun(B\Sigma_n,\Cat)$
corresponds to $g:B\Sigma_n\to B\Sigma_n$.
In concrete terms, the data of a functor $h:BG\to B\Sigma_n$ amounts to a left action of $G=\Hom_{BG}(*_{BG},*_{BG})$ on $\Hom_{B\Sigma_n}(*_{B\Sigma_n},*_{B\Sigma_n})$ in the obvious way,
where $*_{BG}$ and $*_{B\Sigma_n}$ denote unique
objects in $BG$ and $B\Sigma_n$ respectively
(keep in mind the case $G=\Sigma_n$).
A left action $G=\Hom_{BG}(*_{BG},*_{BG})$ on $\Hom_{B\Sigma_n}(*_{B\Sigma_n},*_{B\Sigma_n})$ corresponds to a natural transformation from
the constant functor $BG\to \Cat$ taking the value $\Delta^0$ to the
constant functor taking value $B\Sigma_n$. It relates
$g$ with $f$. The identity functor $B\Sigma_n\to B\Sigma_n$ corresponds to
the natural left multiplication $\Sigma_n$ on $\Sigma_n=\Hom_{B\Sigma_n}(*_{B\Sigma_n},*_{B\Sigma_n})$.
Therefore it is enough to prove that
$f$ corresponds to the natural left multiplication $\Sigma_n$ on $\Sigma_n=\Hom_{B\Sigma_n}(*_{B\Sigma_n},*_{B\Sigma_n})$.
Recall that $f$ is induced by the $n$-fold tensor product of $B\Sigma$.
By the definition of the commutative constraint of $B\Sigma$,
the (trivial) permutation action of $\Sigma_n$ on $(B\Sigma_1)^{\times n}$
gives rise to the left multiplication of $\Sigma_n$ on
$\Sigma_n=\Hom_{B\Sigma_n}(*_{B\Sigma_n},*_{B\Sigma_n})$
(consider the natural transformations given by the commutative constraint
\[
\xymatrix{
B\Sigma_1^{\times n}\ar[rr]^{\textup{trivial action of }\sigma\in \Sigma_n} \ar[dr]_f & & B\Sigma_1^{\times n} \ar[dl]^f \\
 & B\Sigma_n & 
}
\]
which give rise to the action of $\Sigma_n$ on
$\Hom_{B\Sigma_n}(*_{B\Sigma_n},*_{B\Sigma_n})$).
Hence we conclude that $g$ is the identity.
\QED

Consider the presentable $\infty$-category $\Fun(B\Sigma^{op},\SSS)$.
According to \cite[4.8.1.10, 4.8.1.12]{HA},
$\Fun(B\Sigma^{op},\SSS)$ inherits from $B\Sigma$ a symmetric monoidal
structure with the following properties:
\begin{itemize}
\item The Yoneda embedding $B\Sigma\hookrightarrow \Fun(B\Sigma^{op},\SSS)$
is extended to a symmetric monoidal functor.

\item The tensor product $\otimes:\Fun(B\Sigma^{op},\SSS)\times \Fun(B\Sigma^{op},\SSS)\to \Fun(B\Sigma^{op},\SSS)$ preserves small colimits separately in each variable.
\end{itemize}
Hence $\Fun(B\Sigma^{op},\SSS)$ belongs to $\CAlg(\PR)$,
and let us consider the coproduct
\[
\Fun(B\Sigma^{op},\SSS)\otimes \Mod_R^\otimes
\]
in $\CAlg(\PR)$ for a commutative ring spectrum $R$.
Namely, 
$\Fun(B\Sigma^{op},\SSS)\otimes \Mod_R^\otimes$ lies in $\CAlg(\PR_R)$.

\begin{Proposition}
\label{freeleftKan}
The sequence of functors $\Delta^0\stackrel{v}{\to} B\Sigma \hookrightarrow \Fun(B\Sigma^{op}, \SSS)\to \Fun(B\Sigma^{op},\SSS)\otimes \Mod_R$
induces a homotopy equivalence
\[
\Map_{\CAlg(\PR_R)}(\Fun(B\Sigma^{op},\SSS)\otimes \Mod^\otimes_R, \CCC^\otimes) \stackrel{\sim}{\longrightarrow} \Map(\Delta^0,\CCC)=\CCC^\simeq
\]
for any $\CCC^\otimes \in \CAlg(\PR_R)$.
\end{Proposition}

\Proof
We note the three points:
\begin{itemize}
\item $B\Sigma\simeq \textup{Free}(\Delta^0)$ by Proposition~\ref{freebasic},

\item $\Map_{\CAlg(\PR)}(\Fun(B\Sigma^{op},\SSS),\CCC^\otimes)\simeq \Map_{\CAlg(\widehat{\textup{Cat}}_\infty)}(B\Sigma,\CCC^\otimes)$ by \cite[4.8.1.10]{HA},

\item $\Map_{\CAlg(\PR_R)}(\Fun(B\Sigma^{op},\SSS)\otimes \Mod_R^\otimes,\CCC^\otimes)\simeq \Map_{\CAlg(\PR)}(\Fun(B\Sigma^{op},\SSS),\CCC^\otimes)$ by the adjoint pair $(-)\otimes \Mod_R^\otimes:\CAlg(\PR) \rightleftarrows \CAlg(\PR_R):\textup{forget}$.

\end{itemize}
Our claim follows.
\QED

Suppose that $R$ is the Eilenberg-Maclane spectrum $Hk$
of the field $k$ of characteristic zero.
We relate $\Fun(B\Sigma^{op},\SSS)\otimes \Mod^\otimes_R$
with $\DDD(B\Sigma,k)$.

\begin{Proposition}
\label{super}
There exists an
equivalence
$\Fun(B\Sigma^{op},\SSS)\otimes \Mod^\otimes_k\simeq \DDD^\otimes(B\Sigma,k)$
in $\CAlg(\PR_k)$.
\end{Proposition}

\Proof
Let
$B\Sigma \to \Fun(B\Sigma^{op},\SSS)\otimes\Mod_k^\otimes$
be the symmetric monoidal functor
given by
\[
B\Sigma \hookrightarrow \Fun(B\Sigma^{op},\SSS) \to \Fun(B\Sigma^{op},\SSS_*) \to \Fun(B\Sigma^{op},\SP) \to \Fun(B\Sigma^{op},\SP)\otimes \Mod_k
\]
where the first functor is the Yoneda embedding, the subsequent
functors are given by compositions with
$\SSS\to \SSS_*\stackrel{\Sigma^{\infty}}{\to}\SP=\Mod_{\mathbb{S}}\to \Mod_k$.
If we identify $\Fun(B\Sigma^{op},\SP)\otimes \Mod_k$ with
$\Fun(B\Sigma^{op},\Mod_k)$ in light of Lemma~\ref{linearleftkan},
the image of $\bar{r}\in B\Sigma$ that lies in
$\Fun(B\Sigma^{op},\Mod_k)\simeq \prod_{n\ge0}\Fun(B\Sigma_n^{op},\Mod_k)$
is $J^r:=(J^r_n)_{n \ge 0}$ such that $J^r_n \in \Fun(B\Sigma_n^{op},\Mod_k)$,
$J^r_r=\oplus_{g\in \Sigma_r} k\cdot g=k[\Sigma_r] \in \Mod_k$ equipped with
the right multiplication of $\Sigma_r$,
and $J^r_n=0$ for $n\neq r$.
Here we regard $J^r_n$ as an object in $\Mod_k$ endowed with right action
of $\Sigma_n$ (arising from the functoriality of $B\Sigma_n^{op}\to \Mod_k$).

By Proposition~\ref{freeleftKan},
the object $I^1$ in $\DDD(B\Sigma,k)$
induces a morphism $\phi:\Fun(B\Sigma^{op},\SSS)\otimes\Mod_k^\otimes\to \DDD^\otimes(B\Sigma,k)$ in $\CAlg(\PR_k)$.
To prove that it is a symmetric monoidal equivalence,
it will suffice to show that $\phi$ induces
a categorical equivalence of underlying $\infty$-categories.
Since $\DDD(B\Sigma,k)$ and $\Fun(B\Sigma^{op},\Mod_k)$ are stable,
and $\Fun(B\Sigma^{op},\Mod_k)\to \DDD(B\Sigma,k)$ is exact,
it follows that it is enough to prove that $\phi$
induces an equivalence between their homotopy categories (see e.g. \cite[Lemma 5.8]{Tan}).
Since $J^1$ maps to $I^1$, we see that $(J^1)^{\otimes r}=J^r$ maps to
$(I^1)^{\otimes r}=I^r$.
Thus the colimit-preserving functor
\[
\phi:\Fun(B\Sigma^{op},\Mod_k)\simeq \prod_{n\ge0}\Fun(B\Sigma_n^{op},\Mod_k) \to \DDD(B\Sigma,k)\simeq \prod_{n\ge0} \Fun(B\Sigma_n^{op},\Mod_k)
\]
is determined by the product of each restriction
$\phi_n:\Fun(B\Sigma_n^{op},\Mod_k)\to \Fun(B\Sigma_n^{op},\Mod_k)$.
Here $\Fun(B\Sigma_n^{op},\Mod_k)$
is considered to be the full subcategory of $\prod_{n\ge0}\Fun(B\Sigma_n^{op},\Mod_k)$ spanned by $(E_i)_{i\ge0}$
such that $E_i=0$ for $n\neq i$.
Thus it will suffice to prove that $\phi_n$ induces an equivalence 
of homotopy categories.
To this end, consider
the map
\[
\theta:\Hom_{\textup{h}(\Fun(B\Sigma_n^{op},\Mod_k))}(k[\Sigma_n],k[\Sigma_n])\to \Hom_{\textup{h}(\Fun(B\Sigma_n^{op},\Mod_k))}(k[\Sigma_n],k[\Sigma_n])
\]
induced by
$\textup{h}(\phi_n):\textup{h}(\Fun(B\Sigma_n^{op},\Mod_k)) \to \textup{h}(\Fun(B\Sigma_n^{op},\Mod_k))$.
Recall that $\textup{h}(\Fun(B\Sigma_n^{op},\Mod_k))$ is the (unbounded) derived category of $k$-linear representations of $\Sigma_n$.
Note that the category of $k$-linear representations of $\Sigma_n$ is
semi-simple, and every irreducible representation of $\Sigma_n$
is isomorphic to a direct summand of $k[\Sigma_n]$.
Therefore, to show that 
the exact functor $\textup{h}(\phi_n)$ of triangulated categories is an
equivalence, we are reduced to proving
that $\theta$ is a bijective map.
Observe that
$\Hom_{\textup{h}(\Fun(B\Sigma_n^{op},\Mod_k))}(k[\Sigma_n],k[\Sigma_n])$
can be identified with the set of homomorphisms $k[\Sigma_n]\to k[\Sigma_n]$
as right $k[\Sigma_n]$-modules. Thus it is isomorphic to
$k[\Sigma_n]$, and
we can view $\theta$ as a $k$-linear morphism
$\xi:k[\Sigma_n]\to k[\Sigma_n]$.
By the construction of $\phi$,
$\textup{h}(\phi_n)$ commutes with the natural functor
$B\Sigma_n\to \textup{h}(\Fun(B\Sigma_n^{op},\Mod_k))$.
Hence the $k$-linear map $\xi:k[\Sigma_n]\to k[\Sigma_n]$ preserves
$\Sigma_n\subset k[\Sigma_n]$. It follows that
$\theta$ is a bijective map.
\QED

{\it Proof of Proposition~\ref{absolute1}.}
It follows from
Proposition~\ref{freeleftKan} and~\ref{super}.
\QED

\vspace{5mm}


Let $K$ be the standard representation of $\GL_d$,
that is, $k^{\oplus d}$ endowed with the natural action of $\GL_d$.
Applying Proposition~\ref{absolute1} to $K$
we obtain a morphism in $\CAlg(\PR_k)$:
\[
u:\DDD^\otimes(B\Sigma,k)\to \DDD^\otimes(B\GL_d)\simeq \Rep^\otimes(\GL_d)= \QC^\otimes(B\GL_d)
\]
which carries $I^1$ to $K$ placed in degree zero.
Since $I^{n}=(I^1)^{\otimes n}$, we have $u(I^n)=K^{\otimes n}$.
Moreover, we have

\begin{Proposition}
\label{vanish}
Suppose that $W$ is a representation of $\Sigma_n$ viewed
as an object in $\Fun(B\Sigma_n^{op},\Mod_k)\subset \DDD(B\Sigma,k)$.
Then $u(W)\simeq W\otimes_{k[\Sigma_n]}K^{\otimes n}$.
\end{Proposition}

\Proof
Note first that $W$ can be described
as a coproduct
of retracts in $k[\Sigma_n]$.
Thus we may and will assume that $W$ is a retract of $k[\Sigma_n]$.
Since $W$ is a retract,
it is a filtered colimit of a 
linearly ordered sequence consisting of the idempotent maps
(the standard heart consisting of part of (co)homological degree zero is closed under formulation of filtered colimits).
Since $u$ preserves small
colimits, $u(W)$ is a filtered colimit of the linearly ordered sequence
of idempotent maps between
$u(I^n)=K^{\otimes n}\simeq K[\Sigma_n]\otimes_{K[\Sigma_n]}K^{\otimes n}$.
The standard heart of $\DDD(B\GL_d)$ is also closed under filtered 
colimits. Thus, we conclude that
$u(W)\simeq W\otimes_{k[\Sigma_n]}K^{\otimes n}$.
\QED

Before proceeding further we need the representation theory
of symmetric groups.
Let us recall that every representation of a symmetric group
can be constructed by means of Young diagrams (see e.g. \cite[4.1]{Ful1}, \cite[Section 7]{Ful2}):
Let $\lambda$ be a Young diagram having $n$ boxes. Then after choosing 
a standard Young tableau whose underlying Young diagram
is $\lambda$, we can associate to it an idempotent map
between $k[\Sigma_n]$ called the Young symmetrizer. Its retract $V_\lambda$
(the image of the
idempotent map) is an irreducible representation of $\Sigma_n$.
The isomorphism class of $V_\lambda$ (as a representation) does not depend
on the choice of a Young tableau.
Any irreducible representation of a symmetric group is obtained in this way
for a unique Young diagram. The tensor product of irreducible representations
can be described by the Littlewood-Richardson rule.

Let us consider any $k$-linear representation
of $\Sigma_n$ for $n\ge 0$ as an object in $\Fun(B\Sigma_n^{op},\Mod_k)\subset \DDD(B\Sigma,k)$.
Let $T$ be the set consisting of objects $W$
in $\DDD(B\Sigma,k)$
such that $W$ is of the form $V[r]$
such that $[r]$ indicates the shift for $r\in \ZZ$,
$V$ is an irreducible representation of some $\Sigma_n$
associated to Young diagrams having more than $d$ rows.

\begin{Lemma}
\label{preuniversal2}
There is an object $\DDD^\otimes(B\Sigma,k)_{T}$ and a morphism $\DDD^\otimes(B\Sigma,k)\to \DDD^\otimes(B\Sigma,k)_{T}$
in $\CAlg(\PR_k)$
such that the composition induces a homotopy equivalence of spaces
\[
\Map_{\CAlg(\PR_k)}(\DDD^\otimes(B\Sigma,k)_{T},\mathcal{C}^\otimes)\to \Map_{\CAlg(\PR_k)}^{\wedge}(\DDD^\otimes(B\Sigma,k),\CCC^\otimes)
\]
for any $\CCC^\otimes\in \CAlg(\PR_k)$.

Here the
space on the right-hand side denotes the full subcategory
of $\Map_{\CAlg(\PR_k)}(\DDD^\otimes(B\Sigma,k),\CCC^\otimes)$,
spanned by those functors which carry the $(d+1)$-fold wedge 
product $\wedge^{d+1}(I^1)$ to
zero, where $d$ is the integer that appears in the definition of $T$.
\end{Lemma}

\Proof
Notice first
that for any $\Mod_{k}^\otimes \to \mathcal{E}^\otimes\in \CAlg(\PR)_{\Mod_k^\otimes/}\simeq \CAlg(\PR_k)$
\[
\Map_{\CAlg(\PR_k)}(\mathcal{E}^\otimes,\mathcal{C}^\otimes)\simeq \Map_{\CAlg(\PR)}(\mathcal{E}^\otimes,\mathcal{C}^\otimes)\times_{\Map_{\CAlg(\PR)}(\Mod_k^\otimes,\mathcal{C}^\otimes)}\{s\}
\]
where the right-hand side denotes the fiber product (in $\SSS$),
and $s:\Mod_k^\otimes \to \CCC^\otimes$ is the structure functor of $\CCC^\otimes \in \CAlg(\PR)_{\Mod_k^\otimes/}$.
Therefore we may replace $\CAlg(\PR_k)$ by $\CAlg(\PR)$
in the statement.
We apply symmetric monoidal localizations \cite[4.1.3.4]{HA} to
$T':=\{W\to 0\}_{W\in T}$.
For this, we need to show that for any $W\in T$ and any $C\in \DDD(B\Sigma,k)$, $W\otimes C$ is a coproduct of objects in $T$
(it follows that $W\otimes C\to 0$ belongs to a strongly saturated class
generated by the small set $T'$; cf. \cite[5.5.4.5]{HTT}).
We deduce it from Littlewood-Richardson rule or its special case: Pieri rule
(see
\cite[Section 5]{Ful2}).
For this purpose,
we may assume that $W$ is an irreducible representation $V_{\lambda}$
associated to a Young diagram $\lambda$ having $m$ rows with $m>d$, and
$C$ is an irreducible representation $V_{\mu}$
associated to a Young diagram $\mu$.
Let $\alpha=(1,\ldots,1)$ be the Young diagram
corresponding to the partition $m=1+\ldots +1$ of $m$,
 that is, $\alpha$ has $m$ boxes in one column.
Let $\lambda-\alpha$ be the Young diagram
obtained from $\lambda$ by removing $m$ boxes from the left end column.
Then by Littlewood-Richardson rule we see that
the decomposition in $\DDD(B\Sigma,k)$
\[
V_\alpha\otimes V_{\lambda-\alpha} \simeq \oplus_{\nu}V_\nu
\]
where the right-hand side is a coproduct
of those $V_\nu$ such that Young diagram $\nu$
is obtained from $\lambda-\alpha$ by adding $m$
boxed, with no two in the same row.
Hence $V_{\lambda}$ is a retract of $V_\alpha\otimes V_{\lambda-\alpha}$.
Thus it is enough to
prove that $V_{\alpha}\otimes V_{\lambda-\alpha}\otimes V_\mu$
is decomposed into a coproduct of the representations $V_{\beta}$ such that
$\beta$ has more than $d$ rows.
For this, we may replace $V_{\lambda-\alpha}\otimes V_{\mu}$
by $V_{\mu}$. Then again by Littlewood-Richardson rule
we see that $V_{\alpha}\otimes V_{\mu}$ is decomposed into 
$\oplus_{\beta}V_\beta$ where $\beta$
run over the set of
Young diagrams obtained from $\mu$ by adding $m$
boxed, with no two in the same row.
In particular,
$\beta$ has at least $m$ rows.
Consequently, we can apply symmetric monoidal localization \cite[4.1.3.4]{HA}
with respect to $T'$; inverting $T'$ we obtain
$\DDD^\otimes(B\Sigma,k)\to \DDD^\otimes(B\Sigma,k)_T:=\DDD^\otimes(B\Sigma,k)[T^{'-1}]$ which induces a homotopy equivalence
\[
\Map_{\CAlg(\PR)}(\DDD^\otimes(B\Sigma,k)_T,\mathcal{C}^\otimes)\to \Map_{\CAlg(\PR)}^{T}(\DDD^\otimes(B\Sigma,k),\CCC^\otimes)
\]
where the superscript $T$ indicates the full subcategory consisting of
those functors which carry all objects in $T$ to zero.
Finally, we prove that
any morphism $F:\DDD^\otimes(B\Sigma,k)\to \CCC^\otimes$
in $\CAlg(\PR)$
sends all objects in $T$ to zero
if and only if it sends $\wedge^{d+1}(I^1)$ to zero.
The ``only if'' direction is obvious
since the $(d+1)$-fold wedge product is obtained from $k[\Sigma_{d+1}]\simeq (I^1)^{\otimes d+1}$ by using the Young
symmetrizer arising from 
the Young diagram having $d+1$ boxes in one column.
Suppose that $F$ sends $\wedge^{d+1}(I^1)$ to zero.
As observed above, if the Young diagram $\lambda$ has $m$ rows with
$m>d$, then
$V_{\lambda}\in \DDD(B\Sigma,k)$ is 
a retract of a tensor product of $\wedge^{m}(I^1)$ and another
object. Therefore $V_{\lambda}$ maps to zero.
\QED

\begin{Remark}
The underlying functor
$\DDD(B\Sigma,k)\to \DDD(B\Sigma,k)_T$
is a localization (cf. \cite[5.2.7.2]{HTT}),
i.e., a left adjoint functor which has a fully faithful
right adjoint functor whose essential image consists of $T'$-local objects.
It sends $C$ to a $T'$-local object $C_T$
such that the unit map $C\to C_T$ is a $T'$-equivalence
(cf. \cite[5.2.7, 5.5.4.1, 5.5.4.15]{HTT}).
Suppose that $C$ is $\oplus_{i\in I}M_i$  of a coproduct of those
$M_i$ such that $M_i$ is of the form $N[r]$
where $N$ is an irreducible representation of some $\Sigma_m$
and $r\in \ZZ$.
Then $C_T$ is isomorphic to the retract of $\oplus_{i\in I}M_i$
obtained by removing retracts belonging to $T$.
\end{Remark}

For an irreducible representation $V_\lambda$
of $\Sigma_n$ associated to a Young diagram $\lambda$,
$V_{\lambda}\otimes_{k[\Sigma_n]}K^{\otimes n}$ is zero if and only if
the number of rows of $\lambda$ is bigger than $d$.
By Proposition~\ref{vanish},
we see that $u(W)\simeq 0$ for any $W\in T$.
Hence invoking Lemma~\ref{preuniversal2}
we obtain a morphism $u_T:\DDD^\otimes(B\Sigma,k)_T\to \DDD^\otimes(B\GL_d)$ induced by $u:\DDD^\otimes(B\Sigma,k)\to \DDD^\otimes(B\GL_d)$.
Let $\DDD(B\GL_d)_{\textup{eff}}$
be the stable subcategory which contains the standard representation $K$ and the unit
and is closed under tensor products and
coproducts.
The stable presentable full subcategory
$\DDD(B\GL_d)_{\textup{eff}}$ inherits a symmetric monoidal
structure from $\DDD^\otimes(B\GL_d)$.

\begin{Proposition}
\label{step1}
The functor
$u_T:\DDD(B\Sigma,k)_T\to \DDD(B\GL_d)$
is a fully faithful functor whose essential image is
$\DDD(B\GL_d)_{\textup{eff}}$.
In particular, $\DDD^\otimes(B\Sigma,k)_T\simeq \DDD^\otimes(B\GL_d)_{\textup{eff}}$.
\end{Proposition}

Before the proof, let us recall the consequences from Schur-Weyl duality.
Let $V_\lambda$ be the irreducible representation of $\Sigma_n$
associated to a Young diagram $\lambda$ having $n$ boxes.
Then if $\lambda$ has at most $d$ rows, $V_\lambda\otimes_{k[\Sigma_n]}K^{\otimes n}$ is a nonzero irreducible
representation of $\GL_d$. If $\lambda$ has $m$ rows with $m>d$,
then $V_\lambda\otimes_{k[\Sigma_n]}K^{\otimes n}\simeq 0$.
One can obtain any irreducible representation of $\GL_d$
which is a retract of the power $K^{\otimes n}$ in this way for a unique Young diagram.

\Proof
We first prove that
$u_T:\DDD^\otimes(B\Sigma,k)_T \to \DDD^\otimes(B\GL_d)_{\textup{eff}}$ is essentially surjective.
Note that by the semi-simplicity
any object in $\DDD(B\GL_d)_{\textup{eff}}$
is isomorphic to a coproduct $\oplus_{i\in I}P_i$ such that 
$P_i$ is (up to shift)
equivalent to an irreducible representation of $\GL_d$
which is a retract of $K^{\otimes n}$ for some $n\ge0$.
For any nonzero irreducible representation $W$ of $\GL_d$
contained in $K^{\otimes n}$,
there is a unique irreducible representaion of
$V$ of $\Sigma_n$, up to isomorphisms,
 such that $V\otimes_{K[\Sigma_n]}K^{\otimes n}\simeq W$.
Thus Proposition~\ref{vanish} implies that $u_T$ is
essentially surjective.
Next we will prove that $u_T:\DDD^\otimes(B\Sigma,k)_T \to \DDD^\otimes(B\GL_d)_{\textup{eff}}$ is fully faithful.
Let $C$ and $D$ be objects in $\DDD(B\Sigma,k)$.
We may and will assume that
$C$ lies in $\Fun(B\Sigma_{n}^{op},\Mod_k)$
and $D$ lies in $\Fun(B\Sigma_{m}^{op},\Mod_k)$.
Suppose that $n\neq m$.
Then $\Map_{\DDD(B\Sigma,k)}(C,D)$ is a contractible space.
On the other hand, if $P[r], Q[s]\in \DDD(B\GL_d)$
such that $r$ and $s$ are integers, and
$P$ and $Q$ are retracts in
$K^{\otimes n}$ and $K^{\otimes m}$ respectively, then $\Map_{\DDD(B\GL_d)}(P[r],Q[s])$ is a contractible space (for weight reasons).
It follows that
$\Delta^0\simeq \Map_{\DDD(B\Sigma,k)}(C,D)\to \Map_{\DDD(B\GL_d)}(u_T(C),u_T(D))\simeq \Delta^0$
is a homotopy equivalence.
Finally, consider the case of $n=m$.
To prove 
$\Map_{\DDD(B\Sigma,k)}(C,D)\to \Map_{\DDD(B\GL_d)}(u_T(C),u_T(D))$
is a homotopy equivalence, using
decompositions and shifts we are reduced to the case when
$C$ and $E$ are irreducible representations of $\Sigma_n$,
and $D=E[r]$ for some $r\in \ZZ$.
When $C\simeq E$ and $r\ge0$, we have
a natural homotopy equivalence
$k[r] \simeq \Map_{\DDD(B\Sigma,k)}(C,D)\to \Map_{\DDD(B\GL_d)}(u_T(C),u_T(D))\simeq k[r]$.
Here for a space $S$, by $S\simeq k[r]$ we means that $\pi_r(S)\simeq k$
and $\pi_l(S)$ is trivial for $l\neq r$ (i.e., an Eilenberg-MacLane space).
When either $C$ is not equivalent to $E$ or $r<0$,
both $\Map_{\DDD(B\Sigma,k)}(C,D)$
and $\Map_{\DDD(B\GL_d)}(u_T(C),u_T(D))$
are contractible.
This proves that $u_T$ is fully faithful.
\QED

Let $\CCC^\otimes\in \CAlg(\PR_k)$ and let $C$ be an object in $\CCC$.
Then there is a categorical construction
which makes $C$ an invertible object, i.e.,
there is an object $C^\vee$ such that $C\otimes C^\vee$ is a unit of
$\CCC$. Namely, we say that
$\CCC^\otimes \to \CCC^{\otimes}[C^{-1}]$ in $\CAlg(\PR_k)$
is the inversion of $C$ if it induces a homotopy equivalence
\[
\Map_{\CAlg(\PR_k)}(\CCC^{\otimes}[C^{-1}],\mathcal{E}^\otimes)\to \Map^C_{\CAlg(\PR_k)}(\CCC^{\otimes},\mathcal{E}^\otimes)
\]
for any $\mathcal{E}^\otimes\in \CAlg(\PR_k)$, where the superscript $C$
in the right-hand side indicates that we consider only those functors which carry $C$
to an invertible object in $\mathcal{E}^\otimes$.
By a result of Robalo \cite[Proposition 2.9]{Rob}, for any such
$\CCC^\otimes$ and $C\in \CCC$, such an inversion exists.

\begin{Proposition}
\label{step2}
Let $U:=\wedge^{d} K$ be the $d$-fold wedge product of the standard representation.
We denote by $\DDD^\otimes(B\GL_d)_{\textup{eff}}\to \DDD^\otimes(B\GL_d)_{\textup{eff}}[U^{-1}]$ the inversion of $U$.
Then the natural inclusion
$\DDD^\otimes(B\GL_d)_{\textup{eff}}\hookrightarrow \DDD^\otimes(B\GL_d)$ induces an equivalence 
\[
\DDD^\otimes(B\GL_d)_{\textup{eff}}[U^{-1}]\to 
\DDD^\otimes(B\GL_d).
\]
\end{Proposition}

\Proof
Let $\DDD_c(B\GL_d)_{\textup{eff}}$ be the 
stable subcategory of $\DDD(B\GL_d)_{\textup{eff}}$
spanned by compact objects.
Namely, $\DDD_c(B\GL_d)_{\textup{eff}}$
consists of those objects which are
of the form
$\oplus_{i\in I}N_i[r_i]$ where each $r_i$ is an integer,
and each $N_i$ is an
irreducible representation
which belongs to $\DDD_c(B\GL_d)_{\textup{eff}}$.
The small stable $\infty$-category $\DDD_c(B\GL_d)_{\textup{eff}}$
inherits
a symmetric monoidal structure in the natural way.
By \cite[Proposition 2.1, Proposition 2.2]{Rob}, there is the ``small version''
of the inversion of $U$; there exist
a small symmetric monoidal $\infty$-category $\DDD_c^\otimes(B\GL_d)_{\textup{eff}}[U^{-1}]$
and 
a symmetric monoidal
functor $\DDD_c^\otimes(B\GL_d)_{\textup{eff}}\to \DDD_c^\otimes(B\GL_d)_{\textup{eff}}[U^{-1}]$
which induces
a homotopy equivalence
\[
\Map_{\CAlg(\Cat)}(\DDD_c^\otimes(B\GL_d)_{\textup{eff}}[U^{-1}],\mathcal{E}^\otimes) \to \Map_{\CAlg(\Cat)}^U(\DDD_c^\otimes(B\GL_d)_{\textup{eff}},\mathcal{E}^\otimes)
\]
for any $\mathcal{E}^\otimes \in \CAlg(\Cat)$, where the superscript $U$
on the right-hand side indicates the full subcategory consisting
of those functors which carry $U$
to an invertible object in $\mathcal{E}^\otimes$.
Then since $U$ is a symmetric object in the sense of \cite{Rob},
it follows from \cite[Proposition 2.19, Corollary 2.22]{Rob}
that the underlying $\infty$-category
$\DDD_c(B\GL_d)_{\textup{eff}}[U^{-1}]$ is equivalent to a colimit of the
linearly ordered sequence
\[
\DDD_c(B\GL_d)_{\textup{eff}}\stackrel{\otimes U}{\to} \DDD_c(B\GL_d)_{\textup{eff}} \stackrel{\otimes U}{\to} \DDD_c(B\GL_d)_{\textup{eff}}\stackrel{\otimes U}{\to} \ldots
\]
in $\CAlg(\Cat)$ (in \cite[Corollary 2.22]{Rob}, the presentable situation
is treated but the proof is also applicable to this case).
In particular, $\DDD_c(B\GL_d)_{\textup{eff}}[U^{-1}]$
is a stable $\infty$-category since 
the filtered colimit of stable $\infty$-categories
in $\Cat$ is also a stable $\infty$-category \cite[1.1.4.6]{HA}.
Since $(-)\otimes U:\DDD_c(B\GL_d)_{\textup{eff}}\to \DDD_c(B\GL_d)_{\textup{eff}}$ is a fully faithful exact functor
and $(-)\otimes U:\DDD_c(B\GL_d)\to \DDD_c(B\GL_d)$
is an equivalence,
the colimit can be identified with the essential image of the natural functor induced by the inclusion
$\DDD_c(B\GL_d)_{\textup{eff}}\hookrightarrow \DDD_c(B\GL_d)$:
\begin{eqnarray*}
\colim (\DDD_c(B\GL_d)_{\textup{eff}} \stackrel{\otimes U}{\to} \ldots) &\to& 
\colim (\DDD_c(B\GL_d)\stackrel{\otimes U}{\to} \DDD_c(B\GL_d) \stackrel{\otimes U}{\to}\ldots) \\
&\simeq&  \DDD_c(B\GL_d).
\end{eqnarray*}
Since every object in $\DDD_c(B\GL_d)$
has the form $(U^\vee)^{\otimes m} \otimes W$
such that $m \in \NN$, and $W$ belongs to $\DDD_c(B\GL_d)_{\textup{eff}}$,
we thus conclude that
the colimit is $\DDD_c(B\GL_d)$.
Hence we deduce that
the natural symmetric monoidal
functor $\DDD_c^\otimes(B\GL_d)_{\textup{eff}}[U^{-1}]\to 
\DDD_c^\otimes(B\GL_d)$ is an equivalence.
Note that since $(-)\otimes U:\DDD_c(B\GL_d)_{\textup{eff}} \to\DDD_c(B\GL_d)_{\textup{eff}}$ preserves finite colimits, then for any
symmetric monoidal stable
$\infty$-category $\mathcal{E}^\otimes$
a symmetric monoidal functor
$\DDD_c(B\GL_d)_{\textup{eff}}[U^{-1}] \to \mathcal{E}$
preserves finite colimits
if and only if composite
$\DDD_c(B\GL_d)_{\textup{eff}}\to \DDD_c(B\GL_d)_{\textup{eff}}[U^{-1}] \to \mathcal{E}$
preserves finite colimits.
Hence we have a fully faithful functor 
\[
\alpha:\Map_{\CAlg(\widehat{\textup{Cat}}_{\infty})}^{\textup{ex}}(\DDD_c^\otimes(B\GL_d)_{\textup{eff}}[U^{-1}],\mathcal{E}^\otimes)\to \Map_{\CAlg(\widehat{\textup{Cat}}_\infty)}^{\textup{ex}}(\DDD_c^\otimes(B\GL_d)_{\textup{eff}},\mathcal{E}^\otimes)
\]
where by ``$\textup{ex}$'' indicates the full subcategory spanned by
exact functors, i.e., functors which preserve finite colimits.
The essential image consists of those functors $F:\DDD_c^\otimes(B\GL_d)_{\textup{eff}}\to \mathcal{E}^\otimes$
which carry $U$ to an invertible object.
Since $\DDD(B\GL_d)$ is compactly generated,
the symmetric monoidal
Ind-category $\Ind(\DDD_c^\otimes(B\GL_d))$
(cf. \cite[5.3.6.8]{HTT}, \cite[4.8.1.14, 4.8.1.15]{HA})
is equivalent to $\DDD^\otimes(B\GL_d)$.
Similarly, 
$\Ind(\DDD_c^\otimes(B\GL_d)_{\textup{eff}})$
is equivalent to $\DDD^\otimes(B\GL_d)_{\textup{eff}}$.
The left Kan extension
$\Ind(\DDD_c^\otimes(B\GL_d)_{\textup{eff}}) \simeq \DDD^\otimes(B\GL_d)_{\textup{eff}}\to \mathcal{E}^\otimes$ (cf. \cite[4.8.1.14]{HA})
preserves small colimits if and only if 
the composite $\DDD_c^\otimes(B\GL_d)_{\textup{eff}} \to \mathcal{E}^\otimes$
preserves finite colimits (see \cite[the proof of 1.1.3.6]{HA}).
Thus we have a homotopy equivalence
\[
\beta:\Map_{\CAlg(\PR)}(\DDD^\otimes(B\GL_d)_{\textup{eff}},\mathcal{E}^\otimes)
\simeq \Map_{\CAlg(\widehat{\textup{Cat}}_\infty)}^{\textup{ex}}(\DDD_c^\otimes(B\GL_d)_{\textup{eff}},\mathcal{E}^\otimes)
\]
for any $\mathcal{E}^\otimes\in \CAlg(\PR_{\mathbb{S}})$.
Similarly, we have
\[
\gamma:\Map_{\CAlg(\PR)}(\DDD^\otimes(B\GL_d),\mathcal{E}^\otimes)
\simeq \Map_{\CAlg(\widehat{\textup{Cat}}_\infty)}^{\textup{ex}}(\DDD_c^\otimes(B\GL_d),\mathcal{E}^\otimes).
\]
Combining these $\alpha, \beta$ and $\gamma$, we obtain a homotopy equivalence
\[
\Map_{\CAlg(\PR)}(\DDD^\otimes(B\GL_d),\mathcal{E}^\otimes)\to \Map^U_{\CAlg(\PR)}(\DDD^\otimes(B\GL_d)_{\textup{eff}},\mathcal{E}^\otimes)
\]
induced by
$\DDD^\otimes(B\GL_d)_{\textup{eff}}\to \DDD^\otimes(B\GL_d)$.
Note that by \cite[Corollary 2.23]{Rob},
$\DDD^\otimes(B\GL_d)_{\textup{eff}}[U^{-1}]$
is stable.
Hence
$\DDD^\otimes(B\GL_d)_{\textup{eff}}[U^{-1}] \simeq \DDD^\otimes(B\GL_d)$.
\QED

\begin{Proposition}
\label{step3}
We adopt the notation in the paragraph preceding Lemma~\ref{preuniversal2}. In particular, $T=\{V_{\lambda}[r]\}_{\lambda,r\in \ZZ}$ where $\lambda$ run over the set of Young diagrams having more than $d$ rows.
Let $L$ be the $d$-fold wedge product of $I^1$ in $\DDD(B\Sigma,k)$.
Then there exists a natural equivalence
\[
\DDD^\otimes(B\Sigma,k)_T[L^{-1}]\simeq \DDD^\otimes(B\GL_d).
\]
\end{Proposition}

\Proof
Combine Proposition~\ref{step1} and~\ref{step2}.
\QED

\begin{Remark}
If we identify $\DDD^\otimes(B\Sigma,k)$
with $\Fun(\textup{Free}(\Delta^0)^{op},\SSS)\otimes\Mod_k^\otimes$
by Proposition~\ref{super}, we have
$(\Fun(\textup{Free}(\Delta^0)^{op},\SSS)\otimes\Mod_k^\otimes)_T[L^{-1}]\simeq \DDD^\otimes(B\GL_d)$.
\end{Remark}

{\it Proof of Theorem~\ref{characterization}}.
Consider the sequence of functors
\[
\Delta^0\to \textup{Free}(\Delta^0) \to  \Fun(\textup{Free}(\Delta^0)^{op},\Mod_k) \simeq \DDD(B\Sigma,k)\stackrel{s}{\to}
\DDD(B\Sigma,k)_T\stackrel{t}{\to} \DDD(B\Sigma,k)_T[L^{-1}].
\]
The left functor is induced by the adjoint pair
$\textup{Free}:\Cat\rightleftarrows \CAlg(\Cat):\textup{forget}$,
$\textup{Free}(\Delta^0) \to  \Fun(\textup{Free}(\Delta^0)^{op},\Mod_k)$
is the ``natural'' functor,
and the middle equivalence follows from Proposition~\ref{freebasic} and ~\ref{super}.
The functors $s$ and $t$ are left adjoint functors arising from
the localization and the inversion respectively.
The composition with this sequence gives rise to
\[
\alpha:\Map^\otimes_{\CAlg(\PR_k)}(\DDD(B\Sigma,k)_T[L^{-1}],\CCC^\otimes)\to \Map(\Delta^0,\CCC)=\CCC^{\simeq}
\]
for any $\CCC^\otimes\in \CAlg(\PR_k)$.
Combining Proposition~\ref{absolute1}, Lemma~\ref{preuniversal2},
and universal properties,
we deduce that $\alpha$ is fully faithful and its essential
image is $\CCC^{\simeq}_{\wedge,d}$.
Note that through the equivalence $\DDD(B\Sigma,k)_T\simeq \DDD(B\GL_d)_{\textup{eff}}$,
$I^1$ corresponds to $K$, and $L$ corresponds to $U$.
Finally, according to Proposition~\ref{step3},
$\DDD^\otimes(B\Sigma,k)_T[L^{-1}]\simeq \DDD^\otimes(B\GL_d)=\textup{Rep}^\otimes(\GL_d)$.
Therefore, our assertion follows.
\QED

\section{Tannakian characterization}
\label{TC}

\subsection{}
In this Section we prove Theorem~\ref{intmain1}; see 
Theorem~\ref{algebraic} and Theorem~\ref{main1}.
We also describe an explicit presentation (construction)
of $A$ and $G$ in $[\Spec A/G]$ in Theorem~\ref{algebraic} and~\ref{main1}.
We begin by treating its algebraic version,
that is, the case when a fine $\infty$-category 
admits a single wedge-finite generator.
Throughout this Section,
$k$ is a field of characteristic zero.

\begin{Theorem}
\label{algebraic}
Let $\CCC^\otimes$ be a $k$-linear
symmetric monoidal presentable
$\infty$-category. That is, $\CCC^\otimes$ belongs to $\CAlg(\PR_k)$.
Then the following conditions are equivalent:
\begin{enumerate}
\renewcommand{\labelenumi}{(\theenumi)}

\item There exists a wedge-finite object $C$ such that
$\CCC^\otimes$ is generated by $\{C,C^\vee\}$
as a symmetric monoidal stable presentable $\infty$-category.
A unit object $\uni_{\CCC}$ is a compact object.

\item There exist a quotient
stack $[\Spec A/G]$ 
and an equivalence $\CCC^\otimes \simeq \QC^\otimes([\Spec A/G])$
in $\CAlg(\PR_k)$,
where a reductive algebraic group $G$
over $k$ acts on $\Spec A$ with $A\in \CAlg_k$.

\item There exist a quotient stack $[\Spec A/\GL_d]$ and
an equivalence $\CCC^\otimes \simeq \QC^\otimes ([\Spec A/\GL_d])$
in $\CAlg(\PR_k)$,
where the 
general linear group $\GL_d$ for some $d\ge0$
acts on $\Spec A$ with $A\in \CAlg_k$.
\end{enumerate}
\end{Theorem}

\begin{Remark}
The conditions in Theorem~\ref{algebraic} are equivalent to one more
condition, see Corollary~\ref{nonred}.
\end{Remark}

\Proof
The implication from (3) to (2) is obvious.

We will prove that (2) implies (1).
Let $V$ be a finite dimensional
faithful representation of $G$.
If we think of $V$ and $V^\vee$
as objects in $\QC(BG)$, then
$\QC^\otimes(BG)$ is generated by $V$ and $V^\vee$ as a symmetric monoidal
stable presentable $\infty$-category.
Let $\QC^\otimes(BG)\to \QC^\otimes([\Spec A/G])\simeq \Mod_A^\otimes(\QC(BG))$
be the symmetric monoidal functor (informally)
given by $M\mapsto A\otimes M$ (cf. Proposition~\ref{affinecor}).
Since $V$ is wedge-finite in $\QC^\otimes(BG)$,
$A\otimes V$ is wedge-finite.
Observe that $\QC^\otimes([\Spec A/G])$
is generated by $A\otimes V$ and $A\otimes V^\vee$
as a symmetric monoidal stable presentable $\infty$-category.
For the present, we assume that $A\otimes V^{\otimes n}$ and $A\otimes (V^\vee)^{\otimes n}$ are compact.
We will prove that for any
$N\in \QC([\Spec A/G])$, the condition
\[
\Hom_{\textup{h}(\QC([\Spec A/G]))}(A\otimes V^{\otimes n}, N[r])=0\ \ \textup{and}\ \ \Hom_{\textup{h}(\QC([\Spec A/G]))}(A\otimes (V^\vee)^{\otimes n}, N[r])=0
\]
for any $n\ge0$ and any $r\in \ZZ$ implies $N\simeq 0$ (cf. Remark~\ref{homlevelcpt}).
Consider the adjoint pair $A\otimes (-):\QC(BG)\rightleftarrows \QC^\otimes([\Spec A/G]):U$ where $U$ is the forgetful functor.
The vanishing
\[
\Hom_{\textup{h}(\QC(BG))}(V^{\otimes n}, U(N[r]))=0\ \ \textup{and}\ \ 
\Hom_{\textup{h}(\QC(BG))}((V^\vee)^{\otimes n}, U(N[r]))=0
\]
for any $n\ge0$ and $r\in \ZZ$ implies $U(N)=0$.
Using the adjoint pair we conclude that $\QC^\otimes([\Spec A/G])$
is generated by $A\otimes V^{\otimes n}$ and $A\otimes (V^\vee)^{\otimes n}$
($n\ge 0$) as a stable presentable $\infty$-category.
Thus $\QC^\otimes([\Spec A/G])$
is generated by $A\otimes V$ and $A\otimes V^\vee$
as a symmetric monoidal stable presentable $\infty$-category.
Now we will observe that $A\otimes V^{\otimes n}$ and $A\otimes (V^\vee)^{\otimes n}$ are compact.
Taking account of this adjoint pair and the fact that
(i)
a unit in $\QC(BG)$ is compact, (ii) $U$ preserves colimits,
we see that a unit in $\QC([\Spec A/G])$
is compact. Consequently, every dualizable object is compact.
It follows that $A\otimes V^{\otimes n}$ and $A\otimes (V^\vee)^{\otimes n}$ are compact.
Hence (2) implies (1).

Finally, we will prove that (3) follows from (1).
Suppose that there is a $d$-dimensional wedge-finite object $C$
such that
$\CCC^\otimes$ is generated by $C$ and $C^\vee$.
By Theorem~\ref{characterization} there
is a morphism $F:\QC^\otimes(B\GL_d)\to \CCC^\otimes$ in $\CAlg(\PR_k)$
which carries the standard representation of $\GL_d$
to $C$. It is unique up to a contractible space of choices.
We apply Proposition~\ref{PMA} to $F$.
To this end, let us verify the existence of a small set of
compact and dualizable objects generating $\QC(B\GL_d)$ as a stable presentable
$\infty$-category; $\{V^{\otimes n}, (V^\vee)^{\otimes n}\}_{n\ge0}$ generates $\QC(B\GL_d)$ as a
stable presentable $\infty$-category. 
Also, $F(V^{\otimes n})$
and $F((V^\vee)^{\otimes n})$ are compact
(notice that
the compactness of the unit implies that
every dualizable object is compact). If $G$ denotes the right adjoint of $F$
and $\uni_{\CCC}$ denotes a unit of $\CCC$,
we let $A=G(\uni_{\CCC})$.
Then since $\uni_{\CCC}$ belongs to $\CAlg(\CCC)$, $G$ is a lax 
symmetric monoidal functor (by relative adjoint functor theorem \cite[7.3.2.7]{HA})
which induces $G:\CAlg(\CCC)\to \CAlg(\QC(B\GL_d))$. Therefore,
$A$ belongs to $\CAlg(\QC(B\GL_d))$.
According to Proposition~\ref{PMA},
there exists an equivalence $\Mod_A^\otimes(\QC(B\GL_d))\simeq \CCC^\otimes$
in $\CAlg(\PR_k)$.
Let $[\Spec A/\GL_d]$ be the quotient stack associated to $A\in \CAlg(\QC(B\GL_d))$.
Then by Proposition~\ref{affinecor},
$\QC^\otimes([\Spec A/\GL_d])\simeq \Mod_A^\otimes(\QC(B\GL_d))\simeq \CCC^\otimes$.
\QED

\begin{Remark}
By the proof, we can take $d$ in the condition (3)
to be the dimension of a wedge-finite object $C$ in the condition (1).
\end{Remark}

\begin{Definition}
When $\CCC^\otimes$ satisfies the conditions in Theorem~\ref{algebraic},
we shall refer to $\CCC^\otimes$ as a {\it fine algebraic $\infty$-category}.
See also Remark~\ref{finiteproduct}.
\end{Definition}

\begin{Theorem}
\label{main1}
Let $\mathcal{C}^\otimes$ be a $k$-linear symmetric monoidal
stable presentable
$\infty$-category. That is, $\CCC^\otimes$ belongs to $\CAlg(\PR_k)$.
The following conditions are equivalent to one another:
\begin{enumerate}
\renewcommand{\labelenumi}{(\theenumi)}

\item $\CCC^{\otimes}$ is a fine $\infty$-category.

\item There exist a quotient stack $X=[\Spec A/G]$
and an equivalence
$\CCC^\otimes\simeq\QC^\otimes(X)$ in $\CAlg(\PR_k)$
where a pro-reductive group $G$ acts on a derived affine scheme $\Spec A$ with $A\in \CAlg_k$.
\end{enumerate}
\end{Theorem}

We deduce Theorem~\ref{main1} from Theorem~\ref{characterization}
and an elaborated version of arguments in Theorem~\ref{algebraic}.
We will need a few preliminaries.

\vspace{2mm}

Let $X$ be a derived stack over the base field $k$.
We say that $X$ is perfect if the following conditions are satisfied:

\begin{itemize}
\item $\QC(X)$ is compactly generated.

\item Compact and dualizable objects in $\QC(X)$
coincide.
\end{itemize}
The notion of perfect stacks is introduced in
\cite[Definition 3.2, Proposition 3.9]{BFN} in a slightly
different setting.
If $G$ is a reductive algebraic group over $k$ (or pro-reductive group),
then $\QC(BG)\simeq \mathcal{D}(BG)$ satisfies these conditions.
In fact,
the set of finite-dimensional irreducible representations of $G$
generates $\QC(BG)$ as a stable presentable $\infty$-category,
and each irreducible representations is compact in $\QC(BG)$.
Hence $\QC(BG)$ is compactly generated. Moreover,
a unit object is compact. It follows that every 
dualizable object is compact.
Thus, $\QC(BG)$ satisfies the above conditions.
We remark that our definition of $\QC(BG)$ (and pullback functors
between them) agrees with
that of \cite{BFN} if $G$
is an algebraic group over $k$
(see Example~\ref{classifyingstackqc} and \cite[Section 3.1]{BFN}).
In the following proposition,
we use \cite[Theorem 1.2 (1)]{BFN} for the product of
classifying stacks of reductive algebraic groups.

\begin{Proposition}[\cite{BFN}]
\label{Moritaproduct}
Suppose that $X=[\Spec A/G]$ and $[\Spec B/H]$
where $A,B\in \CAlg_k$, and $G$ and $H$ are pro-reductive groups over $k$.
Let $p_X^*:\QC^\otimes(X) \to \QC^\otimes(X\times_k Y)$
and $p_Y^*:\QC^\otimes(Y) \to \QC^\otimes(X\times_k Y)$
be the pullback functors of natural projections.
Let $\QC^\otimes(X)\otimes_k\QC^\otimes(Y)$ denote
the coproduct of $\QC^\otimes(X)$ and $\QC^\otimes(Y)$
in $\CAlg(\PR_k)$,
and let
\[
F:\QC^\otimes(X)\otimes_k\QC^\otimes(Y)\to \QC^\otimes(X\times_k Y)
\]
be the symmetric monoidal functor induced by $p_X^*$
and $p_Y^*$. Then $F$ is an equivalence.

\end{Proposition}

\Proof
This assertion follows from the proof of \cite[Theorem 1.2]{BFN};
our notion of derived stacks is slightly different from that of
\cite{BFN}, but the argument is applicable to our setting.
For the reader's convenience
we outline the proof (to fit our situation).
We note that by \cite[3.2.4.7]{HA}
the underlying $\infty$-category
of $\QC^\otimes(X)\otimes_k\QC^\otimes(Y)$
is a tensor product of $\QC(X)$ and $\QC(Y)$
in $\PR_k$.
It is enough to prove the underlying functor
of $F$ is an equivalence of $\infty$-categories.
When $X=BG$ and $Y=BH$ where $G$ and $H$ are reductive algebraic 
groups, an equivalence of the canonical functor $F:\QC(BG)\otimes_k \QC(BH)\to \QC(BG\times_k BH)$ is a special case of \cite[Theorem 1.2 (1)]{BFN}
(note that for these stacks
our notion of quasi-coherent complexes coincides with
that of \cite{BFN}).

Suppose that
 $G$ and $H$ are pro-reductive,
 and $G=\varprojlim G_{\alpha}$
 and $H=\varprojlim H_\beta$
are directed projective limits of reductive
algebraic groups such that projections $G\to G_\alpha$
and $H\to H_\beta$ are surjective.
The case of $X=BG$ and $Y=BH$ (or the general case)
also follows from the
proof of \cite[Theorem 1.2]{BFN}.
Here we give another (ad-hoc) argument.
As observed in Lemma~\ref{replimit} below,
$\QC(BG)$ is a colimit of $\QC(BG_\alpha)$
in $\PR_k$, and
$\QC(BH)$ is a colimit of $\QC(BH_\beta)$
in $\PR_k$.
Since the tensor product $\otimes_k$ preserves colimits
separately in each variable,
we have
$\QC(BG)\otimes_k\QC(BH)\simeq \QC(B(G\times_kH))\simeq \QC(BG\times_kBH)$
using presentations as colimits.
In particular, $\QC^\otimes(BG)\otimes_k\QC^\otimes(BH)\simeq \QC^\otimes(BG\times_kBH)$.

Next we consider the general 
case where $X=[\Spec A/G]$ and $Y=[\Spec B/H]$.
Let $A\boxtimes B$ denote the tensor product of $p_{BG}^*(A)$ and $p_{BH}^*(B)$
as objects in $\CAlg(\QC(BG\times_kBH))$, where $p_{BG}$ and $p_{BH}$ are
natural projections.
Then $A\boxtimes B$ in $\CAlg(\QC(BG\times_kBH))$
gives rise to the quotient stack $[\Spec (A\boxtimes B)/(G\times_kH)]$,
that is equivalent to $[\Spec A/G]\times_k[\Spec B/H]$.
Then 
we have a natural equivalence
\[
\QC([\Spec A/G]\times_k[\Spec B/H]) \simeq \Mod_{A\boxtimes B}(\QC(BG\times_kBH)),
\]
and by \cite[Proposition 4.1 (2)]{BFN} both sides are also equivalent to 
\[
\Mod_{p_{BG}^*(A)}(\QC(BG\times_k BH))\otimes_{\QC(BG\times_kBH)} \Mod_{p_{BH}^*(B)}(\QC(BG\times_k BH)).
\]
In addition, according to \cite[Proposition 4.1 (1)]{BFN} we have
\[
\Mod_{p_{BG}^*(A)}(\QC(BG\times_k BH))\simeq \Mod_{A}(\QC(BG))\otimes_{\QC(BG)}\QC(BG\times_kBH)\]
and
\[
\Mod_{p_{BH}^*(A)}(\QC(BG\times_k BH))\simeq \Mod_{B}(\QC(BH))\otimes_{\QC(BH)}\QC(BG\times_kBH).
\]
Using these equivalences together with $\QC(BG)\otimes_k\QC(BH)\simeq \QC(BG\times_kBH)$, we obtain
\[
\QC([\Spec A/G]\times_k[\Spec B/H])\simeq \Mod_A(\QC(BG))\otimes_k\Mod_{B}(\QC(BH))
\]
where the right-hand side is naturally equivalent to $\QC([\Spec A/G])\otimes_k\QC([\Spec B/H])$.
\QED

\begin{Lemma}
\label{replimit}
Let $G=\varprojlim_{\beta<\alpha}G_{\beta}$ be a limit of pro-reductive groups indexed by
a limit ordinal $\alpha$. Namely,
 $G=G_\alpha$ is a limit of the sequence
\[
\ldots \to G_{\beta+1}\to G_\beta\to \ldots \to G_1\to G_0
\]
as an affine group scheme, where for any $\beta<\alpha$, $G_{\beta}$ is a pro-reductive group over $k$.
Suppose that for any $\gamma<\beta$
the morphism $G_{\beta}\to G_\gamma$ is surjective.
Then the pullback functors induce an equivalence
\[
\varinjlim \Perf^{\otimes}(BG_{\beta})\to \Perf^\otimes(BG)
\]
where the left-hand side is a colimit in $\CAlg(\Cat)$.
Here $\Perf^{\otimes}(BG_{\beta})$ denotes the stable
subcategory of $\QC^{\otimes}(BG_{\beta})$ spanned by dualizable
objects (note that it coincides with $\DDD_c(BG_\beta)$). 

Moreover, the above equivalence is extended to
an equivalence
\[
\varinjlim \QC^{\otimes}(BG_{\beta})\to \QC^\otimes(BG)
\]
in $\CAlg(\PR_k)$ where
the left-hand side is a colimit in $\CAlg(\PR_k)$.
\end{Lemma}

\Proof
Let $G=G_{\alpha}$.
Note first that for $\alpha \ge \beta \ge \gamma$,
the surjective map $G_{\beta}\to G_{\gamma}$
induces a fully faithful pullback functor
$\Perf(BG_{\gamma})\to \Perf(BG_{\beta})$.
In fact, taking account of the semi-simplicity of the representations of
$G_{\beta}$, we see that
any object $W$ in $\Perf(BG_{\beta})$
has the form $V_0[r_0]\oplus\ldots \oplus V_n[r_n]$
where $V_i$ is a finite dimensional irreducible representation of $G_{\beta}$
and $r_i$ is an integer
for any $n\ge i \ge 0$.
Moreover, $\Hom_{\textup{h}(\QC(BG_{\beta}))}(V_i,V_i[r])$
is a division algebra for $r=0$, and it is zero if $r\neq 0$.
Thus we conclude that
$\Perf(BG_{\gamma})\to \Perf(BG_{\beta})$
is fully faithful, and  
its essential image is spanned by those objects
which have the form
$V_0[r_0]\oplus\ldots \oplus V_n[r_n]$
where 
$V_i$ is an irreducible representation of $G_{\beta}$ arising from
the factorization
$G_{\beta}\to G_{\gamma}$,
and $r_i$ is an integer
for any $n\ge i \ge 0$ (keep in mind that an exact functor
between stable $\infty$-categories is an equivalence if and only if
the induced functor between their homotopy categories is an equivalence,
see e.g. \cite{Tan}).
To prove an equivalence
$\varinjlim \Perf^{\otimes}(BG_{\beta})\to \Perf^\otimes(BG)$,
by \cite[3.2.3.1]{HA}
it is enough to show that
the colimit $\varinjlim \Perf(BG_{\beta})$
in $\Cat$ is naturally equivalent to $\Perf(BG)$.
For this, since each $\Perf(BG_{\gamma})\to \Perf(BG_{\beta})$
is fully faithful, it will suffice to observe that
every object $C$ in $\Perf(BG_{\alpha})$ belongs to
$\Perf(BG_{\beta})$ for some $\beta<\alpha$.
Let $A_{\beta}$ denote the ring of functions
on $G_{\beta}$, that is endowed with a structure of a commutative Hopf algebra.
The formulation 
$G_{\alpha}=\varprojlim_{\beta<\alpha}G_{\beta}$
of the limit gives rise to
$A_{\alpha}=\cup_{\beta<\alpha}A_{\beta}$, where
we regard $A_{\beta}$ as a Hopf subalgebra of $A_{\alpha}$.
Let $W\simeq V_0[r_0]\oplus\ldots \oplus V_n[r_n]$
be an object in $\Perf(BG_{\alpha})$
where 
$V_i$ is a finite dimensional irreducible representation of $G_{\alpha}$, and $r_i$ is an integer
for any $n\ge i \ge 0$.
Since each $V_i$ is finite dimensional,
the corresponding coaction $V_i\to V_i\otimes A_{\alpha}$
factors through $V_i\to V_i\otimes H_i$
for a finitely generated commutative
Hopf algebra $H_i\subset A_{\alpha}$.
Let $\{x_1^i,\ldots,x^i_{s_i}\}$ be the set of generators of
$H_i$ as a commutative $k$-algebra.
If we take a sufficiently large $\beta<\alpha$,
$x_j^i$ lies in $A_{\beta}$ for any $i$ and $j$.
Therefore all $H_i$ are contained in $A_{\beta}$
It follows that $W$ belongs to $\Perf(BG_{\beta})$.

Next we prove that $\varinjlim_{\beta<\alpha} \QC^{\otimes}(BG_{\beta})\to \QC^\otimes(BG)
$.
By taking a left Kan extension \cite[4.8.1.14]{HA}
$\Perf^\otimes(BG_\beta)\to \QC^\otimes(BG)$ is extended to
$\Ind(\Perf^\otimes(BG_\beta))\to \QC^\otimes(BG)$
which preserves small colimits.
Observe that $\Ind(\Perf^\otimes(BG_\beta))\simeq \QC^\otimes(BG_\beta)$.
Since objects in $\Perf(BG_\beta)$
are compact in $\QC(BG_\beta)$, it follows from \cite[5.3.4.12]{HTT}
that
the left Kan extension $\Ind(\Perf(BG_\beta))\to \QC(BG_\beta)$
is fully faithful.
Note that $G_\beta$ is a pro-reductive group, and therefore
the abelian category of representations of $G_\beta$ is
semi-simple. As is well-known,
every representation $W$ of $G_{\beta}$
can be described as a filtered colimit $\varinjlim V_z$
of finite dimensional subrepresentations $V_z$.
Thus $\Ind(\Perf(BG_\beta))\to \QC(BG_\beta)$ is essentially surjective.
Using the equivalence
\[
\Map_{\CAlg(\Cat)}^{\textup{ex}}(\Perf^\otimes(BG_{\beta}),\DDD^\otimes)\simeq \Map_{\CAlg(\PR)}(\QC^\otimes(BG_\beta),\DDD^\otimes)
\]
for $\DDD^\otimes\in \CAlg(\PR_{\mathbb{S}})$
and $\alpha \ge \beta$
we deduce that $\QC^\otimes(BG)$ is a filtered colimit
$\varinjlim_{\beta<\alpha}\QC^\otimes(BG_\beta)$ in $\CAlg(\PR)$.
Here the superscript ``$\textup{ex}$'' indicates the full subcategory
spanned by exact functors.
By \cite[4.2.3.5, 3.2.3.1]{HA},
$\QC^\otimes(BG_\beta)$ is a colimit $
\varinjlim_{\beta<\alpha}\QC^\otimes(BG_\beta)$ in $\CAlg(\PR_k)$.
\QED

{\it Proof of Theorem~\ref{main1}}.
We prove that (1) implies (2).
Let $\CCC^\otimes$ be an object in $\CAlg(\PR_k)$.
Let $\{C_{\lambda}\}_{\lambda\in \Lambda}$
be a small set of wedge-finite objects
such that $\CCC^\otimes$ is generated by 
$\{C_{\lambda},C^\vee_\lambda \}_{\lambda\in \Lambda}$.
Choose a bijective map $\Lambda\simeq \alpha$
where $\alpha$ is a cardinal.
We replace $\{C_{\lambda}\}_{\lambda\in \Lambda}$
by $\{C_{\beta}\}_{\beta<\alpha}$.
We will construct a pro-reductive group $G$
and a morphism $F:\Rep^\otimes(G)\simeq \QC^\otimes(BG)\to \CCC^\otimes$
by transfinite induction.

Let $n_\beta$ be the dimension of the wedge-finite object $C_{\beta}$.
By Theorem~\ref{characterization},
there exists
a morphism $F_1:\QC^\otimes(B\GL_{n_{0}})\to \CCC^\otimes$
in $\CAlg(\PR_k)$ which carries the standard representation of
$\GL_{n_{0}}$ (placed in degree zero) to $C_0$.
Set $G_1:=\GL_{n_{0}}$.

Suppose that $G_{\beta}$ and $F_{\beta}:\QC^\otimes(BG_\beta)\to \CCC^\otimes$
have been constructed for $\beta$.
In addition, assume that $G_{\beta}=\varprojlim_{\gamma<\beta}G_\gamma$
 if $\beta$ is a limit ordinal, and $G_{\beta}=G_{\beta-1}\times_k\GL_{n_{\beta-1}}$ if otherwise (by convention $G_0$ is trivial). Moreover, suppose that $G_{\beta}\to G_{\gamma}$
 is surjective for $\gamma<\beta$.
By Theorem~\ref{characterization} we have
$F_{\beta+1}':\QC^\otimes(B\GL_{n_{\beta}})\to \CCC^\otimes$
which carries the standard representation of $\GL_{n_{\beta}}$
to $C_{\beta}$.
Using Proposition~\ref{Moritaproduct} we prove that
\[
\QC^\otimes(BG_{\beta}\times_kB\GL_{n_{\beta}})\simeq \QC^\otimes(BG_{\beta})\otimes_k\QC^\otimes(B\GL_{n_{\beta}}).
\]
Then the ``coproduct'' of
$F_{\beta}$ and $F'_{\beta+1}$ induces
\[
F_{\beta+1}:\QC^\otimes(BG_{\beta}\times_kB\GL_{n_{\beta}})\simeq \QC^\otimes(BG_{\beta})\otimes_k\QC^\otimes(B\GL_{n_{\beta}}) \to \CCC^\otimes.
\]
Here the ``coproduct" is that in $\CAlg(\PR_k)$.
Note that by \cite[5.5.8.11, 5.5.8.12]{HTT}
$B(G_{\beta}\times_k\GL_{n_{\beta}})\simeq BG_{\beta}\times_kB\GL_{n_{\beta}}$.
We define $G_{\beta+1}$ to be $G_{\beta}\times_k\GL_{n_{\beta}}$.
If $p_{\beta}:G_{\beta+1}=G_{\beta}\times_k\GL_{n_{\beta}}\to G_{\beta}$
is the first projection, then
we have a commutative diagram (i.e. 2-cell)
\[
\xymatrix{
\QC^\otimes (BG_{\beta})\ar[r]^{p_{\beta}^*} \ar[rd]_{F_\beta} & \QC^\otimes(BG_{\beta+1}) \ar[d]^{F_{\beta+1}} \\
 & \CCC^\otimes
}
\]
in $\CAlg(\PR_k)$.

Let $\beta$ be a limit ordinal.
Suppose that
a linearly ordered sequence indexed by $\beta$
\[
\cdots \to G_{\gamma+1}\stackrel{p_{\gamma}}{\to} G_{\gamma}\to \cdots \stackrel{p_1}{\to} G_1
\]
of pro-reductive groups
and 
\[
\xymatrix{
\QC^\otimes(BG_{1})\ar[r]^{p_1^*} \ar[dr]_{F_0} & \cdots \ar[r] & \QC^\otimes(BG_{\gamma})\ar[r]^{p_{\gamma}^*} \ar[dl]_{F_\gamma} & \QC^\otimes(BG_{\gamma+1}) \ar[r] \ar[dll]_{F_{\gamma+1}} & \cdots \\
 &  \CCC^\otimes & & &
}
\]
in $\CAlg(\PR_k)_{/\CCC^\otimes}$ have been defined.
Suppose that each $p_\gamma$ is surjective.
Let $G_{\beta}:=\varprojlim_{\gamma<\beta}G_{\gamma}$.
Then by Lemma~\ref{replimit}
$\varinjlim_{\gamma<\beta}\QC^\otimes(BG_\gamma)\simeq \QC^\otimes(BG_\beta)$.
Hence by the universal property of the colimit
and Lemma~\ref{replimit}
the above diagram induces
a morphism $\QC^\otimes(BG_\beta)\to \CCC^\otimes$
in $\CAlg(\PR_k)$.
By transfinite induction we have a pro-reductive group $G:=G_\alpha$
and $F:=F_\alpha:\QC^\otimes (BG)\to \CCC^\otimes$.

Next
we prove that
$F:\QC^\otimes (BG)\to \CCC^\otimes$ satisfies the following conditions:
\begin{itemize}
\item There is a small set of compact and dualizable objects $\{I_{\lambda}\}_{\lambda\in \Lambda}$ which generates $\QC(BG)$ as a stable presentable $\infty$-category.

\item $\{F(I_\lambda)\}_{\lambda\in \Lambda}$
is a set of compact objects in $\CCC$ which generates $\CCC$
as a stable presentable $\infty$-category.
\end{itemize}
If we define $\{I_{\lambda}\}_{\lambda\in \Lambda}$
to be the set of (finite dimensional) irreducible representations of $G$,
then the first condition is satisfied.
To check the second condition, note that there are natural
surjective homomorphisms $G\to G_{\beta+1} =G_\beta\times_k \GL_{n_{\beta}}\to  \GL_{n_\beta}$.
The pullback of the composite
induces an irreducible representation of $G$
from the standard representation of $\GL_{n_\beta}$.
Thus 
$\{C_{\lambda},C^\vee_\lambda \}_{\lambda\in \Lambda}$
is contained in the essential image of $F$.
Hence the second condition is satisfied (notice that
dualizable objects are compact in $\CCC$).
Let $H$ be a right adjoint functor of $F$. As in
the proof of Theorem~\ref{algebraic}, $H(\uni_\CCC)$
belongs to $\CAlg(\QC(BG))\simeq \CAlg(\Rep(G))$.
Now we apply to Proposition~\ref{PMA} to $F$
and obtain an equivalence $\QC^\otimes([\Spec A/G])\simeq \Mod_A^\otimes(\Rep(G))\simeq \CCC^\otimes$ where $[\Spec A/G]$ is the quotient stack associated to
$A\in \CAlg(\Rep(G))$.

Next we prove that (2) implies (1).
As in the proof of Theorem~\ref{algebraic},
if $\{I_{\lambda}\}_{\lambda\in \Lambda}$ is the set of
irreducible representations of $G$, then
$\{A\otimes I_\lambda\}_{\lambda\in \Lambda}$
is the set of compact and dualizable objects
which generates $\Mod_A(\QC(BG))=\QC([\Spec A/G])$
as a stable presentable $\infty$-category.
Every $A\otimes I_{\lambda}$ is wedge-finite.
Finally, the unit of $\QC([\Spec A/G])$ is compact
since the unit in $\QC(BG)$ is compact (use adjoint pair
$\QC(BG)\rightleftarrows \QC([\Spec A/G])$).
\QED

\begin{Remark}
\label{finiteproduct}
Let $\mathcal{C}^\otimes$ be a fine $\infty$-category.
Suppose further that there is a finite set $\{C_1,\ldots,C_r\}$
of wedge-finite objects such that
$\mathcal{C}^\otimes$ is generated by
$\{C_1,\ldots,C_r,C_1^\vee,\ldots, C_r^\vee\}$
as a symmetric monoidal stable presentable $\infty$-category.
The proof of Theorem~\ref{main1} reveals that
in that case there exist a quotient stack of the form
$[\Spec A/(\GL_{n_1}\times \cdots\times \GL_{n_r})]$
and an equivalence $\mathcal{C}^\otimes\simeq \QC^\otimes([\Spec A/(\GL_{n_1}\times \cdots\times \GL_{n_r})])$.
In particular, by Theorem~\ref{algebraic}, $\mathcal{C}^\otimes$ be a fine algebraic $\infty$-category.
\end{Remark}

\subsection{}
For a fine $\infty$-category $\CCC^\otimes$
there are many choices of quotient forms
$[\Spec A/G]$ such that $\CCC^\otimes\simeq \QC^\otimes([\Spec A/G])$.
One pleasant
feature of our construction in the proof of Theorem~\ref{algebraic}
and Theorem~\ref{main1}
is that
given a set of wedge-finite generators we have an
explicit quotient form $[\Spec A/G]$.
For example, as in the proof of Theorem~\ref{algebraic}
and Theorem~\ref{main1}, we can take $G$ to be a product of general
linear groups.
It is useful for many applications.
We will describe $A$ in terms of a given set of generators.

To begin, we consider the case when
a fine $\infty$-category
$\CCC^\otimes$ has a single wedge-finite (compact) generator $C$, i.e.,
the fine algebraic case.
Let $d$ be the dimension of $C$.

Let $\lambda$ be a Young diagram with $n$ boxes.
As in the case of $\textup{Alt}^n$,
we let $\mathbb{S}_{\lambda}C$ be the image of
the associated idempotent map
$C^{\otimes n}\to  C^{\otimes n}$ (in the idempotent complete
homotopy category of $\CCC^\otimes$).
To a Young diagram $\lambda$ with $n$ boxes,
by choosing a lift to a Young tableau,
we associate the Young symmetrizer $c_{\lambda}\in \QQ[\Sigma_n]$
which satisfies $c_\lambda c_\lambda=a_{\lambda} c_{\lambda}$ where $a_\lambda$
is a certain rational number (cf. \cite[Lecture 4]{Ful1}).
This $a_\lambda^{-1} c_\lambda$ gives an idempotent map $C^{\otimes n}\to C^{\otimes n}$
via permutation.
We define $\mathbb{S}_{\lambda}C$ to be $\Ker(1- a_\lambda^{-1} c_\lambda)$.

Let $\mathsf{Hom}_{\CCC}(-,-)$ denote the hom complex which belongs
to $\Mod_k$. Namely, for any $D\in \CCC$, we have
the adjoint pair 
\[
D\otimes s(-):\Mod_k \rightleftarrows \CCC:\mathsf{Hom}_{\CCC}(D,-)
\]
where $s$ is the ``structure'' functor $\Mod_k^\otimes\to \CCC^\otimes$,
and 
the existence of the right adjoint functor $\mathsf{Hom}_{\CCC}(D,-)$
is implied by the adjoint functor theorem and the fact that
$D\otimes s(-)$ preserves small colimits.
We often
think of $\Mod_k$ as the $\infty$-category obtained
from the category of (possibly unbounded) complexes of $k$-vector spaces
by inverting quasi-isomorphisms.
By the highest weight theory,
the set of isomorphism classes of irreducible representations of
$\GL_d$ bijectively corresponds to the set 
\[
\ZZ^{\oplus d}_{\star}:=\{\lambda=(\lambda_1,\ldots,\lambda_d)\in \ZZ^{\oplus d}|\ \lambda_1\ge \lambda_2\ge \cdots \ge \lambda_d\}.
\]
That is, when $\lambda_d \ge 0$,
$\lambda=(\lambda_1,\ldots,\lambda_d)$ determines a partition
of $\lambda_1+\ldots +\lambda_d$, and it corresponds to the irreducible representation
$\mathbb{S}_\lambda K$ where $K$ is the standard representation of
$\GL_d$.
When $\lambda_d<0$, $\lambda^+=(\lambda_1-\lambda_d,\lambda_2-\lambda_d,\ldots,\lambda_d-\lambda_d)$ determined a partition of $(\lambda_1+\ldots +\lambda_d)-d\lambda_d$ (regarded as a Young diagram),
and $\lambda$ corresponds to the irreducible representation
$(\mathbb{S}_{\lambda^+}K)\otimes (\wedge^dK^\vee)^{\otimes (-\lambda_d)}$.
If $\lambda_d<0$, we define $\mathbb{S}_{\lambda}K$ to be
$(\mathbb{S}_{\lambda^+}K)\otimes (\wedge^dK^\vee)^{\otimes (-\lambda_d)}$.
Replacing $K$ by a wedge-finite object $C$ of a fine $\infty$-category we define $\mathbb{S}_{\lambda}C$ for any $\lambda \in \ZZ^{\oplus d}_{\star}$
in a similar way.

\begin{Proposition}
\label{warm}
Let $\CCC^\otimes$ be a fine algebraic $\infty$-category. Suppose that
a fine $\infty$-category $\CCC^\otimes$
 admits a single $d$-dimensional wedge-finite object $C$
such that $\{C,C^\vee\}$ generates $\CCC^\otimes$ as a symmetric
monoidal stable presentable $\infty$-category.
Then in (3) in Theorem~\ref{algebraic},
we can take a derived stack $[\Spec A/\GL_d]$ 
such that
\[
A\simeq \bigoplus_{\lambda\in \ZZ^{\oplus d}_{\star}}\mathsf{Hom}_\CCC(\mathbb{S}_\lambda C,\uni_{\CCC})\otimes \mathbb{S}_\lambda K
\]
in $\Rep(\GL_d)$.
The action of $\GL_d$ on the right-hand side is through $\mathbb{S}_\lambda K$.
\end{Proposition}

\Proof
In the proof of $(3)\Rightarrow (1)$ in Theorem~\ref{algebraic},
using Theorem~\ref{characterization}
we constructed a $k$-linear symmetric monoidal
colimit-preserving functor $F:\Rep^\otimes(\GL_d)\to \CCC^\otimes$
sending the standard representation $K$ to $C$,
which has a (lax symmetric monoidal)
right adjoint $G:\CCC^\otimes\to \Rep^\otimes(\GL_d)$.
Put $A=G(\uni_\CCC)$. We have proved that $\QC^\otimes([\Spec A/\GL_d])\simeq
\CCC^\otimes$ in Theorem~\ref{algebraic}.
Write
\[
P:=\Hom_{\hhh(\Rep(\GL_d))}((\mathbb{S}_\alpha K)[n], \oplus_{\lambda\in \ZZ^{\oplus d}_{\star}}\mathsf{Hom}_\CCC(\mathbb{S}_\lambda C,\uni_{\CCC})\otimes \mathbb{S}_\lambda K).
\]
There are natural isomorphisms of abelian groups
\begin{eqnarray*}
Q:=\oplus_{\lambda \in \ZZ^{\oplus d}_{\star}}\Hom_{\hhh(\Rep(\GL_d))}((\mathbb{S}_{\alpha}K)[n],\mathsf{Hom}_\CCC(\mathbb{S}_\lambda C,\uni_{\CCC})\otimes \mathbb{S}_\lambda K)&\simeq& H_n(\mathsf{Hom}_\CCC(\mathbb{S}_\alpha C,\uni_{\CCC})) \\
&\simeq& \Hom_{\hhh(\CCC)}((\mathbb{S}_\alpha C)[n],\uni_{\CCC}) \\
&\simeq& \Hom_{\hhh(\Rep(\GL_d))}((\mathbb{S}_\alpha K)[n],A) \\
\end{eqnarray*}
where the final isomorphism is implied by
the adjoint pair (notice also that
$F(\mathbb{S}_\alpha K)=\mathbb{S}_{\alpha}C$).
We can observe the first isomorphism as follows:
Note that by the representation theory of general
linear groups,
$\Hom_{\hhh(\Rep(\GL_d))}(\mathbb{S}_\lambda K,(\mathbb{S}_{\mu}K)[m])$
is isomorphic to $k$ (resp. $0$)
if $\lambda=\mu$
and $m=0$ (resp. if otherwise).
It follows that
\begin{eqnarray*}
Q &\simeq& \Hom_{\hhh(\Rep(\GL_d))}((\mathbb{S}_{\alpha}K)[n],\mathsf{Hom}_\CCC(\mathbb{S}_\alpha C,\uni_{\CCC})\otimes \mathbb{S}_\alpha K) \\ &\simeq& \Hom_{\hhh(\Mod_k)}(k[n],\mathsf{Hom}_\CCC(\mathbb{S}_\alpha C,\uni_{\CCC})) \\
&\simeq& H_n( \mathsf{Hom}_\CCC(\mathbb{S}_\alpha C,\uni_{\CCC})).
\end{eqnarray*}
Thus, we obtain the first isomorphism.
The second isomorphism follows from the sequence of adjunctions
$\SSS\stackrel{e}{\rightleftarrows} \Mod_{k}\stackrel{D\otimes s(-)}{\rightleftarrows} \CCC$ where $e:\SSS\to \Mod_k$ is a colimit-preserving functor
which carries a contractible space to $k\in \Mod_k$.
By the compactness of $\mathbb{S}_\alpha K$,
the canonical map $Q\to P$ is an isomorphism,
so that $P\simeq \Hom_{\hhh(\Rep(\GL_d))}((\mathbb{S}_\alpha K)[n],A)$.
Every object $M\in \Rep(\GL_d)$ is a coproduct of objects
$(\mathbb{S}_{\alpha}K)[n]$ with $\alpha \in \ZZ^{\oplus d}_{\star}$
and $n\in \ZZ$.
Consequently, we see that $A\simeq \oplus_{\lambda\in \ZZ^{\oplus d}_{\star}}\mathsf{Hom}_\CCC(\mathbb{S}_\lambda C,\uni_{\CCC})\otimes \mathbb{S}_\lambda K$.
\QED

\begin{Remark}
In general, it is difficult
to describe the commutative algebra
structure of
\[
\bigoplus_{\lambda\in \ZZ^{\oplus d}_{\star}}\mathsf{Hom}_\CCC(\mathbb{S}_\lambda C,\uni_{\CCC})\otimes \mathbb{S}_\lambda K
\]
that inherits from $A$ in an explicit way.
(In fact, $A$ is a so-called $E_\infty$-algebra.)
Nevertheless, if one thinks of $A$ as a commutative algebra
object in the symmetric monoidal homotopy category
$\hhh(\Rep(\GL_d))$, there is an explicit presentation of
the multiplication and the unit.
We will describe it.
Note that any object of $\hhh(\Rep(\GL_d))$
is a coproduct of shifted irreducible
representations of the form $(\mathbb{S}_{\lambda}K)[n]$.
Moreover, $\Hom_{\hhh(\Rep(\GL_d))}((\mathbb{S}_{\lambda}K)[n],(\mathbb{S}_{\mu}K)[m])=k$ if $\lambda=\mu$ and $n=m$, and $\Hom_{\hhh(\Rep(\GL_d))}((\mathbb{S}_{\lambda}K)[n],(\mathbb{S}_{\mu}K)[m])=0$ if otherwise.
Let
\[
p:(\SSSS_{\alpha}K)[m]\to A\simeq \oplus \mathsf{Hom}_\CCC(\mathbb{S}_\lambda C,\uni_{\CCC})\otimes \mathbb{S}_\lambda K \ \ 
\textup{and}\ \ p':(\SSSS_{\beta}K)[n]\to A\simeq \oplus \mathsf{Hom}_\CCC(\mathbb{S}_\lambda C,\uni_{\CCC})\otimes \mathbb{S}_\lambda K
\]
be
morphisms in $\hhh(\Rep(\GL_d))$ (for ease of notation, we omit the indexes).
We will consider the composite
\[
c:(\SSSS_{\alpha}K)[m]\otimes (\SSSS_{\beta}K)[n]\stackrel{p\otimes p'}{\to} A\otimes A\to A
\]
where the right map is the multiplication.
We refer to the composite as the multiplication of $p$ and $p'$.
We remark that the multiplication $A\otimes A\to A$
is given by
\[
A\otimes A=G(\uni_{\CCC})\otimes G(\uni_{\CCC})\to G F(G(\uni_{\CCC})\otimes G(\uni_{\CCC}))\simeq G(F(G(\uni_{\CCC}))\otimes F(G(\uni_{\CCC})))\to G(\uni_{\CCC}\otimes \uni_{\CCC})\simeq G(\uni_{\CCC})=A
\]
where the left arrow is induced by
the unit $\textup{id}\to GF$ of the adjunction, and the right arrow
is induced by the counit $FG\to \textup{id}$.
We give an explicit presentation of the composite.
First, let $q:(\SSSS_{\alpha}C)[m]\to \uni_{\CCC}$ and
$q':(\SSSS_{\beta}C)[n]\to \uni_{\CCC}$ be left adjuncts
of $p$ and $p'$, respectively.
We may regard $q$ and $q'$ as elements of $H_{m}(\mathsf{Hom}_{\CCC}((\SSSS_{\alpha}C),\uni_{\CCC}))$ and $H_{n}(\mathsf{Hom}_{\CCC}((\SSSS_{\alpha}C),\uni_{\CCC}))$
determined by $p$ and $p'$, respectively (see the proof of Proposition~\ref{warm}).
Then
\[
(\SSSS_{\alpha}K)[m]\otimes (\SSSS_{\beta}K)[n]\stackrel{c}{\to}A\simeq  \oplus_{\lambda\in \ZZ^{\oplus d}_{\star}} \mathsf{Hom}_\CCC(\mathbb{S}_\lambda C,\uni_{\CCC})\otimes \mathbb{S}_\lambda K
\]
is a right adjunct of
\[
q\otimes q':(\SSSS_{\alpha}C)[m]\otimes (\SSSS_{\beta}C)[n]\to \uni_{\CCC}\otimes \uni_{\CCC}\simeq \uni_{\CCC}.
\]
We easily check this in a categorical way
by using the adjoint pair
$(F,G)$ and the facts
(i) $F$ is symmetric monoidal, (ii) $G$ is lax symmetric monoidal,
and (iii) the unit $\textup{id}\to GF$ is a symmetric monoidal
natural transformation.
Moreover,
there are
isomorphisms
\[
(\SSSS_{\alpha}K)[m]\otimes (\SSSS_{\beta}K)[n]\simeq 
(\SSSS_{\alpha}K)\otimes (\SSSS_{\beta} K)[m+n]\stackrel{l}{\simeq} \oplus_{\lambda\in \ZZ^{\oplus d}_{\star}}(\SSSS_{\lambda}K)^{\oplus r_{\lambda}}[m+n]
\]
where $r_{\lambda}$ are non-negative
integers determined by the Littlewood-Richardson
rule. We remark that $l$ is a fixed
non-canonical isomorphism. Then since $F$ is symmetric monoidal,
the composite
\[
\oplus(\SSSS_{\lambda}K)^{\oplus r_{\lambda}}[m+n]\simeq
(\SSSS_{\alpha}K)[m]\otimes (\SSSS_{\beta}K)[n]\stackrel{c}{\to}\oplus \mathsf{Hom}_\CCC(\mathbb{S}_\lambda C,\uni_{\CCC})\otimes \mathbb{S}_\lambda K
\]
is a right adjunct of $\oplus(\SSSS_{\lambda}C)^{\oplus r_{\lambda}}[m+n]\simeq (\SSSS_{\alpha}C)[m]\otimes (\SSSS_{\beta}C)[n]\stackrel{q\otimes q'}{\to} \uni_{\CCC}\otimes \uni_{\CCC}\simeq \uni_{\CCC}$ (we omit indexes, and $\oplus(\SSSS_{\lambda}C)^{\oplus r_{\lambda}}[m+n]\simeq (\SSSS_{\alpha}C)[m]\otimes (\SSSS_{\beta}C)[n]$ comes from $l$). Namely, each factor
$(\SSSS_{\gamma}K)[m+n]\hookrightarrow \oplus(\SSSS_{\lambda}K)^{\oplus r_{\lambda}}[m+n] \to \oplus \mathsf{Hom}_\CCC(\mathbb{S}_\lambda C,\uni_{\CCC})\otimes \mathbb{S}_\lambda K$
is determined by the class of the factor $(\SSSS_{\gamma}C)[m+n]\hookrightarrow \oplus(\SSSS_{\lambda}C)^{\oplus r_{\lambda}}[m+n]\stackrel{q\otimes q'}{\to}  \uni_{\CCC}$
in $H_{m+n}(\mathsf{Hom}_\CCC(\mathbb{S}_\gamma C,\uni_{\CCC}))$.
Remember that the unit $k\to A$ is the unit $k\to GF(k)=A$ of the adjunction
where $k$ indicates the one dimensional trivial
representation.
Therefore, the unit $k\to \oplus \mathsf{Hom}_\CCC(\mathbb{S}_\lambda C,\uni_{\CCC})\otimes \mathbb{S}_\lambda K$ is determined by the canonical inclusion $k\to \mathsf{Hom}_{\CCC}(\uni_{\CCC},\uni_{\CCC})\otimes k$ induced by the class of the identity in $H_0(\mathsf{Hom}_{\CCC}(\uni_{\CCC},\uni_{\CCC}))$.

We unwind this algebraic
structure in the simplest case where $\GL_{d}=\mathbb{G}_m$,
i.e., $d=1$. In this case,
\[
A\simeq \oplus_{w\in \ZZ}\mathsf{Hom}_{\CCC}(C^{\otimes w},\uni_{\CCC})\otimes \chi_{w}
\]
where $\chi_w$ is the character of $\mathbb{G}_m$ whose weight is $w$,
and $F(\chi_w)=C^{\otimes w}$.
Let $\theta$ be
an element of $H_m(\mathsf{Hom}_{\CCC}(C^{\otimes a},\uni_{\CCC}))$.
The element $\theta$
amounts to a morphism
$\chi_a[m]\to \oplus_{w\in \ZZ}\mathsf{Hom}_{\CCC}(C^{\otimes w},\uni_{\CCC})\otimes \chi_w$ in $\hhh(\Rep(\mathbb{G}_m))$ which we denote by $f_{\theta}$.
Let $\theta'$ be
another element of $H_n(\mathsf{Hom}_{\CCC}(C^{\otimes b},\uni_{\CCC}))$
and let 
$f_{\theta'}:\chi_b[n]\to \oplus_{w\in \ZZ}\mathsf{Hom}_{\CCC}(C^{\otimes w},\uni_{\CCC})\otimes \chi_w$ in $\hhh(\Rep(\mathbb{G}_m))$ be the corresponding morphism.
Then the multiplication $\chi_{a+b}[m+n]\simeq \chi_{a}[m]\otimes \chi_b[n]\to \oplus_{w\in \ZZ}\mathsf{Hom}_{\CCC}(C^{\otimes w},\uni_{\CCC})\otimes \chi_w$
of $f_{\theta}$ and $f_{\theta'}$ corresponds to
the ``multiplication'' of $\theta$ and $\theta'$,
that is, the image of $\theta\otimes \theta'$ under
the canonical map
\[
H_m(\mathsf{Hom}_{\CCC}(C^{\otimes a},\uni_{\CCC}))\otimes H_n(\mathsf{Hom}_{\CCC}(C^{\otimes b},\uni_{\CCC}))\to H_{m+n}(\mathsf{Hom}_{\CCC}(C^{\otimes a}\otimes C^{\otimes b},\uni_{\CCC}\otimes \uni_{\CCC})).
\]
The unit is $k=\chi_0\to \mathsf{Hom}_{\CCC}(\uni_{\CCC},\uni_{\CCC})\otimes \chi_{0}\hookrightarrow \oplus_{w\in \ZZ}\mathsf{Hom}_{\CCC}(C^{\otimes w},\uni_{\CCC})\otimes \chi_{w}$ which is determined by the identity element
$H_0(\mathsf{Hom}_{\CCC}(\uni_{\CCC},\uni_{\CCC}))$.
\end{Remark}

Next we treat an arbitrary fine $\infty$-category.
We first collect some points from the proof of Theorem~\ref{main1}:
Suppose that $\CCC^\otimes$ is a $k$-linear fine $\infty$-category
and $\{C_{\lambda}\}_{\lambda \in \Lambda}$ is 
a set of wedge-finite objects
such that $\{C_{\lambda},C^\vee_{\lambda}\}_{\lambda \in \Lambda}$
generates $\CCC^\otimes$ as a symmetric monoidal stable presentable $\infty$-category.
Then we have constructed a pro-reductive group $G$ and an adjoint pair 
\[
F:\QC^\otimes(BG)\rightleftarrows \CCC^\otimes:H
\]
where $F$ is a $k$-linear symmetric monoidal (left adjoint)
colimit-preserving functor. We put $A=H(\uni_\CCC)$
and proved $\CCC^\otimes\simeq \QC^\otimes([\Spec A/G])$.
By the construction,
$G$ is a product $\prod_{\lambda \in \Lambda}\GL_{n_{\lambda}}$
where $n_{\lambda}$ is the dimension of $C_{\lambda}$.
Hence $G$ has the form $\varprojlim_{S\in P_{\textup{fin}}(\Lambda)} G_{S}$, where $P_{\textup{fin}}(\Lambda)$
is the set of finite subsets of $\Lambda$, and $G_{S}$ denotes the product
of $\prod_{s\in S}\GL_{n_s}$.
The commutative Hopf algebra $\Gamma(G)$ of $G$ is a union of
Hopf subalgebras of $G_{S}$ with $S\in P_{\textup{fin}}(\Lambda)$.
Hence every finite dimensional representation of $G$ factors through
some quotient $G\to G_{S}$.

\begin{Lemma}
\label{irreform}
Every irreducible representation of $G_S=\GL_{n_1}\times \cdots \times \GL_{n_r}$ is of the form $p^*_1(V_{1})\otimes \cdots \otimes p^*_r(V_r)$
such that $V_i$ is an irreducible representation of $\GL_{n_i}$
and $p_i$ is the natural projection $BG_S\to B\GL_{n_i}$.
The endomorphism algebra
$\End_{\hhh(\QC(BG_S))}(p_1^*(V_{1})\otimes \cdots \otimes p^*_r(V_r))$ is $k$.
\end{Lemma}

\begin{Remark}
\label{infinite}
Consequently, every irreducible
representation of $\prod_{\lambda \in \Lambda}\GL_{n_{\lambda}}$
has the form $\otimes_{s\in S}p_s^*(V_s)$ where $S$ is a finite set of
$\Lambda$, $p_s$ is the natural projection $B\prod_{\lambda \in \Lambda}\GL_{n_{\lambda}}\to B\GL_{n_s}$, and $V_s$ is an irreducible representation of
$\GL_{n_s}$.
\end{Remark}

\begin{Remark}
\label{absimple}
We also remark that if each $V_i$ is an irreducible representation
of $\GL_{n_i}$,
then $V_1\otimes \cdots \otimes V_r$ is an irreducible representation
of $G_S$.
Indeed, by applying \cite[Theorem 1.2 (1)]{BFN} to $B\GL_{n_1}\times \cdots \times B\GL_{n_r}$ (it is possible to apply it to a product of reductive algebraic groups), we see
\[
\End_{\hhh(\QC(BG_S))}(p_1^*(V_{1})\otimes \cdots \otimes p^*_r(V_r))\simeq \End(V_1)\otimes_k\cdots \otimes_k \End(V_r)\simeq k\otimes_k\cdots \otimes_kk\simeq k.
\]
\end{Remark}

{\it Proof of Lemma~\ref{irreform}.}
It is a standard fact but we outline the proof for the reader's convenience.
According to \cite[Proposition 3.24]{BFN} (and its proof)
the set of objects $\{p^*_1(V_1)\otimes \cdots\otimes  p^*_r(V_r) \}$ where each $V_i$ run through
irreducible
representations of $\GL_{n_i}$ is a set of compact objects in $\QC(BG_S)$
which generates $\QC(BG_S)$ as a stable presentable $\infty$-category.
We give a direct proof of this fact:
If we let $W_i$ be a finite dimensional
faithful representation
of $\GL_{n_i}$, then $W:=p^*_1(W_{1})\oplus \cdots \oplus p^*_r(W_r)$
is a faithful representation of $G_S$,
so that the set of compact objects
$\{W^{\otimes n}, (W^{\vee})^{\otimes n}\}_{n\ge0}$
generates $\QC(BG_S)$ as a stable presentable $\infty$-category.
It follows that
$\{p^*_1(V_1)\otimes \cdots\otimes  p^*_r(V_r) \}$
generates $\QC(BG_S)$ as a stable presentable $\infty$-category.
Consequently, every irreducible representation $V$ of $G_S$ (regarded as an object
in $\QC(BG_S)$) is a filtered colimits of objects in
$\{p^*_1(V_1)\otimes \cdots\otimes  p^*_r(V_r)[n]\}_{n\in \ZZ}$.
Since the formulation of cohomology groups is compatible
with filtered colimits,
$V$ is a filtered colimit of objects in $\{p^*_1(V_1)\otimes \cdots\otimes  p^*_r(V_r)\}$
in the abelian category of representations of $G_S$.
Consequently, (by the semi-simplicity and irreducibility of $V$)
we deduce that $V$ is isomorphic to
an object of the form $p^*_1(V_1)\otimes \cdots\otimes  p^*_r(V_r)$.
Remark~\ref{absimple} implies the second assertion.
\QED

Using Lemma~\ref{irreform}, Remark~\ref{infinite}, \ref{absimple}
we deduce the following explicit formula as in Proposition~\ref{warm}:

\begin{Proposition}
\label{warm2}
Let
\[
A_S=\bigoplus_{(\alpha_{\xi})\in \sqcap_{\xi\in S} \ZZ_\star^{\oplus n_\xi}} \mathsf{Hom}_{\CCC}(\otimes_{\xi\in S}\mathbb{S}_{\alpha_{\xi}}C_{\xi}, \uni_\CCC)\otimes (\otimes_{\xi \in S}\mathbb{S}_{\alpha_{\xi}}K_{\xi}).
\]
Here $K_\xi$ is the standard representation of $\GL_{n_{\xi}}$ which we naturally regard as an irreducible representation of $G$.
The set $\ZZ_\star^{\oplus n_\xi}$ parameterizes the isomorphism
classes of irreducible representations of $\GL_{n_\xi}$.
Then there exists an equivalence
\[
A \simeq \varinjlim_{S\in P_{\textup{fin}}(\Lambda)} A_S
\]
in $\Rep(G)$. We regard $P_{\textup{fin}}(\Lambda)$
as a poset by inclusions, and $S\hookrightarrow S'$
induces $A_S\to A_{S'}$ in the obvious way.
\end{Proposition}

\subsection{}
As an immediate application, we
 conclude this Section by explaining how to construct the Tannaka dual:

\begin{Remark}[Tannaka dual]
\label{Tandual}
Let $\CCC^\otimes$ be a fine $\infty$-category
and $F:\CCC^\otimes\to \Mod_k^\otimes$ a morphism in $\CAlg(\PR_k)$.
Let $\{C_{\lambda}\}_{\lambda\in \Lambda}$ be a set of wedge-finite generators
and suppose that $F(C_{\lambda})$ in $\Mod_k$ is concentrated in degree zero
for each $\lambda\in \Lambda$.
Then the stack $[\Spec A/G]$ and the equivalence $\QC^\otimes([\Spec A/G])\simeq \CCC^\otimes$
associated to $\CCC^\otimes$ 
and $\{C_{\lambda}\}_{\lambda\in \Lambda}$
give rise to the composite
$\QC^\otimes(BG)\to \QC^\otimes([\Spec A/G])\simeq \CCC^\otimes\stackrel{F}{\to} \Mod_k^\otimes$ where the first functor is the pullback functor induced
by the natural morphism $[\Spec A/G]\to BG$.
By our construction, this composite is equivalent to the forgetful functor.
Then the right adjoint of this
forgetful functor carries the unit of $\Mod_k$ to the ring of
functions $\Gamma(G)$ (placed in degree zero).
It yields a morphism $p:\Spec k\simeq[\Spec \Gamma(G)/G]\to [\Spec A/G]$.
The based loop stack $\Omega_*[\Spec A/G]:=\Spec k\times_{[\Spec A/G]}\Spec k$
is a derived affine group scheme (cf. \cite[Appendix]{Tan}), i.e., a group
object in $\Aff_k$. This fiber product $\Spec k\times_{[\Spec A/G]}\Spec k$
can be regarded as a $G$-equivariant version of the bar construction
(keep in mind that for an argmented object $B\to k$ in $\CAlg_k$
one can think of the pushout
$k\otimes_Bk$ as the ``standard'' bar construction).
Note also that we have a natural identification $F\simeq p^*$.
By the main result of \cite[Theorem 4.8]{Bar},
this derived affine group scheme $\Omega_*[\Spec A/G]$
represents the automorphism group $\Aut(F)$ of $F$ (see \cite{Bar} for the precise
formulation).
We define the {\it Tannaka dual} of $\CCC^\otimes$ with respect to $F$
to be $\Omega_*[\Spec A/G]$.
\end{Remark}

\section{Symmetric monoidal functors and Correspondences}
\label{correspondence}

As observed in the introduction,
a symmetric monoidal functor $\QC^\otimes(Y)\to \QC^\otimes(X)$
is not necessarily the pullback functor
of a morphism $X\to Y$.
For example, by Theorem~\ref{characterization}
giving a $k$-linear symmetric monoidal functor
\[
\QC^\otimes(B\GL_d)\to \QC^\otimes(\Spec k)
\]
amounts to giving a $d$-dimensional
wedge-finite object in $\QC^\otimes(\Spec k)$.
Let $V[2n]$ be a $d$-dimensional $k$-vector space
placed in (homological) degree $2n$.
Then $V[2n]$ is a $d$-dimensional wedge-finite object,
and it gives rise to 
a symmetric monoidal functor
$\phi_{2n}:\QC^\otimes(B\GL_d)\to \QC^\otimes(\Spec k)$
which carries the standard representation of $\GL_d$
to $V[2n]$. On the other hand,
a morphism $\Spec k\to B\GL_d$ of stacks corresponds to
$\GL_d$-torsor over $\Spec k$, that is, the trivial torsor.
In particular, the pullback functor of $\Spec k\to B\GL_d$
sends the standard representation of $\GL_d$ to
a $k$-vector space placed in degree zero.
If $n\neq 0$, then $\phi_{2n}$ is not the pullback functor.
This means that
morphisms of stacks are not enough for our purpose,
and we need a new geometric notion.

In this Section, for a derived stack $X$
we write $\CAlg(X):=\CAlg(\QC(X))$.

\begin{Definition}
Let $X$ and $Y$ be two derived stacks over a base field $k$
of characteristic zero.
A correspondence from $Y$ to $X$ is an object
$P$ of $\CAlg(Y\times_k X)$
such that
\begin{itemize}
\item $(p_Y)_*(P)\simeq \OO_Y$,

\item the composite of
functors $\Mod_P(\QC(Y\times_kX))\stackrel{\textup{forget}}{\to} \QC(Y\times_kX)
\stackrel{(p_Y)_*}{\to}\QC(Y)$ is conservative.
\end{itemize}
Here $p_Y$ is the projection to $Y$.
Let $\Cor_k(Y,X)$ be the full subcategory
of $(\CAlg(Y\times_kX)^{op})^{\simeq}$ spanned by correspondences from $Y$
to $X$. We shall refer to $\Cor(Y,X)$
as the space (or $\infty$-groupoid) of correspondences from $Y$ to $X$.
\end{Definition}

Correspondences can be regarded as
``twisted morphisms''. The notion of correspondences generalizes
that of morphisms of derived stacks.
Namely, there is a natural functor
from the mapping
space $\Map_{\Sh(\Aff_k)}(Y,X)$
to $\Cor(Y,X)$,
see Remark~\ref{morvscor}.

We define the composition of correspondences.
Let $X$, $Y$, and $Z$ are derived stacks over $k$ and
$p_{YX}:Z\times_kY\times_kX\to Y\times_kX$
the natural projection.
The projections $p_{ZY}$ and $p_{YX}$ are defined in a similar manner.
The projection $p_{YX}$ induces $p_{YX}^*:\CAlg(Y\times_kX)\to \CAlg(Z\times_kY\times_kX)$ induced by the pullback functor.
There is the
pushforward functor
$(p_{YX})_*:\CAlg(Z\times_kY\times_kX)\to \CAlg(Y\times_kX)$
which carries $E$ to $(p_{YX})_*(E)$.

We define the functor
\[
\CAlg(Z\times_kY)\times \CAlg(Y\times_kX) \to \CAlg(Z\times_kX);\ \  (P,Q) \mapsto P\star Q
\]
by the formula
$(P,Q) \mapsto (p_{ZX})_*(p_{ZY}^*(P)\cdot p_{YX}^*(Q))$.
Here $p_{ZY}^*(P)\cdot p_{YX}^*(Q)$ denotes
the coproduct $p_{ZY}^*(P)\otimes p_{YX}^*(Q)$ in
$\CAlg(Z\times_kY\times_kX)$.
As discussed in Remark~\ref{kurukuru},
the composition $\CAlg(Z\times_kY)\times \CAlg(Y\times_kX)\to \CAlg(Z\times_kX)$
induces
\[
\Cor(Z,Y)\times \Cor (Y,X)\to \Cor (Z,X).
\]
If we write $\Delta_X:X\to X \times_kX$ for the diagonal,
then $(\Delta_X)_*(\OO_X)$ is the identity correspondence
of $X$.

The purpose of this Section is to prove the following result:

\begin{Theorem}
\label{recmain}
Let $X=[\Spec A/G]$ and $Y=[\Spec B/H]$ be two quotient stacks
where $A,B\in \CAlg_k$, and
$G$ and $H$ are pro-reductive groups
over $k$. Then we have
\begin{enumerate}
\renewcommand{\labelenumi}{(\roman{enumi})}
\item There is a natural homotopy equivalence
\[
\textup{Cor}(Y,X) \to \Map_{\CAlg(\PR_k)}(\QC^\otimes(X),\QC^\otimes(Y))
\]
which carries $P$ to
$P^*$ defined by
\[
P^*:\QC^\otimes(X)\to \QC^\otimes(Y);\ \ M\mapsto  (p_Y)_*(p_X^*(M)\otimes_{\OO_{Y\times_kX}}P)
\]
where $p_X:Y\times_kX\to X$ and $p_Y:Y\times_kX\to Y$ are natural
projections.

\item 
Let $f:\QC^\otimes(X)\to \QC^\otimes(Y)$
and $g:\QC^\otimes(Y)\to \QC^\otimes(Z)$
be morphisms in $\CAlg(\PR_k)$.
Let $C_f\in \Cor(Y,X)$ and $C_g\in \Cor(Z,Y)$ be
correspondences corresponding to $f$ and $g$ respectively.
Then through the equivalence
$\Cor(Z,X)\simeq \Map_{\CAlg(\PR_k)}(\QC^\otimes(X),\QC^\otimes(Z))$,
the composite $g\circ f$ corresponds to
$C_g\star C_f=(p_{ZX})_*(p_{ZY}^*(C_g)\cdot p_{YX}^*(C_f))$.
\end{enumerate}
\end{Theorem}

\begin{Remark}
For $P \in \CAlg(Y\times_kX)$, the functor
$P^*:\QC(Y)\to \QC(X)$ given by
\[
(p_X)_*(p_Y^*(-)\otimes_{\OO_{Y\times_kX}}P)
\]
 is only a lax symmetric monoidal functor, but
if we provide that $P\in \Cor(Y,X)$, then $P^*$ is a symmetric monoidal
functor.
\end{Remark}

\begin{Remark}
\label{morvscor}
Let $f:Y\to X$ be a morphism of derived stacks.
Let $(\textup{id}_{Y}, f):Y\to Y\times_kX$ be the morphism
determined by the identity
and $f$. Then, $(\textup{id}_{Y}, f)_*(\OO_Y)\in \CAlg(Y\times_kX)$.
It gives rise to a natural functor
\[
\Map_{\Sh(\Aff_k)}(Y,X)\to \Cor(Y,X).
\]
Intuitively, we can think that
this functor carries $f:Y\to X$
to the structure sheaf of ``the graph of $f$''.
\end{Remark}

We need some Lemmata for the proof of Theorem~\ref{recmain}.

The opposite $\infty$-category $\Cor(Y,X)^{op}$ of correspondences
can naturally be identified with
the largest Kan subcomplex in the full subcategory of
 $\CAlg(\QC^\otimes(Y\times_kX))\simeq \CAlg(\QC^\otimes(X)\otimes_k\QC^\otimes(Y))$ (cf. Proposition~\ref{Moritaproduct}).

Let $\CAlg'(\PR_k)_{\QC^\otimes(Y\times_kX)/}$
be the full subcategory of 
$\CAlg(\PR_k)_{\QC^\otimes(Y\times_kX)/}$
spanned by those functors
$\phi:\QC^\otimes(Y\times_kX)\to \CCC^\otimes$
such that $\QC^\otimes(Y)\stackrel{p_Y^*}{\to} \QC^\otimes(X\times_kY)\stackrel{\phi}{\to} \CCC^\otimes$
is an equivalence.

There is a functor
\[
\eta:\Cor(Y,X)^{op}\hookrightarrow\CAlg(Y\times_kX)^\simeq\to \CAlg(\PR_k)_{\QC^\otimes(Y\times_kX)/}^{\simeq}
\]
which carries $P$ to $\pi^*:\QC^\otimes(Y\times_kX)\to \Mod_P^\otimes(\QC^\otimes(Y\times_kX))$.
According to \cite[4.8.5.21]{HA} it is fully faithful.
Moreover, we have:

\begin{Lemma}
\label{rec1}
The functor $\eta$ induces
an equivalence
\[
\Cor(Y,X)^{op}\to \CAlg'(\PR_k)^{\simeq}_{\QC^\otimes(Y\times_kX)/}.
\]
\end{Lemma}

\Proof
We first show that for any $P \in \Cor(Y,X)$,
$\eta(P)$ belongs to $\CAlg'(\PR_k)_{\QC^\otimes(Y\times_kX)/}$.
Namely, we will prove that $((-)\otimes P)\circ (p_Y)^*:\QC^\otimes(Y)\to \Mod_P^{\otimes}(\QC(Y\times_kX))$
is an equivalence.
Since $P \in \Cor(Y,X)$,
the pushforward $\Mod_P(\QC(Y\times_kX))\to \QC(Y\times_kX)\stackrel{(p_Y)_*}{\to} \QC(Y)$
is conservative.
Let $\{V_\lambda\}_{\lambda\in \Lambda}$ is a (small)
set of compact (and dualizable) objects which generates
$\QC(Y)$ as a stable presentable $\infty$-category.
One may take $\{V_\lambda\}_{\lambda\in \Lambda}$
to be the set $\{p^*(U_{i})\}_{i\in I}$
where $\{U_i\}_{i\in I}$ is the set of irreducible representations of
$H$, and $p$ is the projection $Y=[\Spec B/H]\to BH$.
Put $V'_\lambda=p_Y^*(V_\lambda)\otimes_{\OO_{Y\times X}}P\in \Mod_{P}(\QC(Y\times_kX))$.
Observe that
$\{V'_{\lambda}\}_{\lambda\in \Lambda}$ is a set of compact and
dualizable objects which generates
$\Mod_{P}(\QC(Y\times_kX))$ as a stable presentable $\infty$-category.
Since unit objects
in $\QC(Y\times_kX)$ and $\Mod_P(\QC(Y\times X))$ are compact
(in fact, by Theorem~\ref{main1}
$\QC(Y\times_kX)$ is a fine $\infty$-category),
dualizable objects $V_\lambda'$ are also compact in $\Mod_{P}(\QC(Y\times_kX))$.
Using the adjoint pair $((-)\otimes P)\circ (p_Y)^*:\QC(Y)\rightleftarrows \Mod_P(\QC(Y\times_kX))$ and the fact that $\Mod_P(\QC(Y\times_kX))\to \QC(Y)$
is conservative, we see that
the vanishing $\Hom_{\textup{h}(\Mod_P(\QC(Y\times_kX)))}(V'_{\lambda},N[r])=0$
for any $(\lambda,r)\in \Lambda\times \ZZ$ implies that $N\simeq 0$.
By Proposition~\ref{PMA},
$((-)\otimes P)\circ p_Y^*:\QC^\otimes(Y)\to \Mod_P^\otimes(\QC(Y\times_kX))$
is extended to an equivalence $\Mod_{(p_Y)_*(P)}^\otimes(\QC(Y))\simeq \Mod_P^\otimes(\QC(Y\times_kX))$
(the composite $\QC^\otimes(Y)\stackrel{(p_Y)_*(P)\otimes(-)}{\longrightarrow} \Mod_{(p_Y)_*(P)}(\QC(Y))\to \Mod_P^\otimes(\QC(Y\times_kX))$
is equivalent to
$((-)\otimes P)\circ  p_Y^*$). By the equivalence
$(p_Y)_*(P)\simeq \OO_Y$, we see that 
$\QC^\otimes(Y)\simeq \Mod_{(p_Y)_*(P)}^\otimes(\QC(Y))\simeq \Mod_P^\otimes(\QC(Y\times_kX))$.

Conversely, suppose that
$\phi:\QC^\otimes(Y\times_kX)\to \CCC^\otimes$
belongs to $\CAlg'(\PR_k)_{\QC^\otimes(Y\times_kX)/}$,
that is, the composite
$\phi\circ p_Y^*:\QC^\otimes(Y)\to\QC^\otimes(Y\times_kX)\to \CCC^\otimes$ is an equivalence.
Let $\psi:\CCC\to \QC(Y\times_kX)\simeq \QC(Y\times_kX)$
be a (lax symmetric monoidal) right adjoint of $\phi$. Put $A=\psi(\uni_{\CCC})\in \CAlg(\QC(Y\times_kX))$ where $\uni_\CCC$ denotes the
a unit of $\CCC$.
Since $\phi\circ p_Y^*$ is an equivalence
and $\QC^\otimes(Y\times_kX)$ and $\QC^\otimes(Y)$ are
fine $\infty$-categories (cf. Theorem~\ref{main1}),
we can apply Proposition~\ref{PMA}
to deduce that $\phi$
is extended to $\Mod_A^\otimes(\QC(Y\times_kX))\simeq \CCC^\otimes$.
Therefore $P$ lies in $\Cor(Y,X)$, and we have the diagram
\[
\xymatrix{
              \QC^\otimes(Y\times_kX)   \ar[d]_{A\otimes(-)} \ar[rd]^{\phi}    &  \\
   \Mod_{A}^\otimes(\QC(Y\times_kX)) \ar[r]_(0.7){\simeq}  & \CCC^\otimes.
}
\]
Hence our claim follows.
\QED

\begin{Lemma}
\label{rec2}
There is a natural homotopy equivalence
\[
\Map_{\CAlg(\PR_k)}(\QC^\otimes(X),\QC^\otimes(Y)) \to \CAlg'(\PR_k)^\simeq_{\QC^\otimes(Y\times_kX)/}.
\]
\end{Lemma}

\Proof
Let $\langle\textup{id}:\QC^\otimes(Y)\to \QC^\otimes(Y)\rangle$
be the full subcategory of $\CAlg(\PR_k)^\simeq_{\QC^\otimes(Y)/}$
spanned  by those objects
$\QC^\otimes(Y)\to \QC^\otimes(Y)$
which are equivalent to
the identify functor $\QC^\otimes(Y)\to \QC^\otimes(Y)$.
It is obvious that $\langle\textup{id}:\QC^\otimes(Y)\to \QC^\otimes(Y)\rangle$
is equivalent to a contractible space, i.e., $\Delta^0$.
Note that if $\CAlg''(\PR_k)^\simeq_{\QC^\otimes(Y\times_kX)/}$
is the full subcategory of $\CAlg'(\PR_k)^\simeq_{\QC^\otimes(Y\times_kX)/}$
spanned by
those objects $\phi:\QC^\otimes(Y\times_kX)\to \CCC^\otimes$
such that $\CCC^\otimes=\QC^\otimes(Y)$ and $\QC^\otimes(Y)\stackrel{p_Y^*}{\to} \QC^\otimes(Y\times_kX)\stackrel{\phi}{\to} \QC^\otimes(Y)$
is equivalent to the identity.
Then the inclusion 
\[
\CAlg''(\PR_k)^\simeq_{\QC^\otimes(Y\times_kX)/}\hookrightarrow \CAlg'(\PR_k)^\simeq_{\QC^\otimes(Y\times_kX)/}
\]
is a homotopy equivalence.
We have a pullback square
\[
\xymatrix{
\Map_{\CAlg(\PR_k)_{\QC^\otimes(Y)/}}(\QC^\otimes(Y\times_kX),\QC^\otimes(Y)) \ar[d] \ar[r] & \CAlg''(\PR_k)^\simeq_{\QC^\otimes(Y\times_kX)/} \ar[d] \\
\langle\textup{id}:\QC^\otimes(Y)\to \QC^\otimes(Y)\rangle \ar[r] & \CAlg(\PR_k)^\simeq_{\QC^\otimes(Y)/}
}
\]
in $\widehat{\SSS}$, where the right vertical functor
is determined by $p_Y^*:\QC^\otimes(Y)\to \QC^\otimes(Y\times_kX)$.
The essential image of the right vertical functor is 
$\langle\textup{id}:\QC^\otimes(Y)\to \QC^\otimes(Y)\rangle$,
and the bottom horizontal arrow is a fully faithful functor.
Therefore the top horizontal functor is an equivalence.
By Proposition~\ref{Moritaproduct},
$\QC^\otimes(Y\times_kX)\simeq \QC^\otimes(X)\otimes_k\QC^\otimes(Y)$.
Thus, the adjoint pair
\[
\QC^\otimes(Y)\otimes_k(-):\CAlg(\PR_k)\rightleftarrows \CAlg(\Mod_{\QC^\otimes(Y)}(\PR_k))
\simeq \CAlg(\PR_k)_{\QC^\otimes(Y)/}:\textup{forget}
\]
implies
a homotopy equivalence
\[
\Map_{\CAlg(\PR_k)_{\QC^\otimes(Y)/}}(\QC^\otimes(Y\times_kX),\QC^\otimes(Y))\simeq \Map_{\CAlg(\PR_k)}(\QC^\otimes(X),\QC^\otimes(Y)).
\]
Hence our assertion follows.
\QED

{\it Proof of Theorem~\ref{recmain} (i).}
Our claim follows from Lemma~\ref{rec1} and Lemma~\ref{rec2}.
\QED

\begin{Remark}
\label{exconst}
Let $f:\QC^\otimes(X)\to \QC^\otimes(Y)$ be a morphism
in $\CAlg(\PR_k)$.
The corresponding correspondence $C_f$ is constructed as follows:
Let $f_Y:\QC^\otimes(X\times_kY)\simeq \QC^\otimes(X)\otimes_k\QC^\otimes(Y)\to \QC^\otimes(Y)$ be a morphism determined by $f$
and the identity functor
$\QC^\otimes(Y)\to \QC^\otimes(Y)$ (note the universal property of
the coproduct). Let $f_Y'$ be a right adjoint of $f_Y$.
Then as the proof of Lemma~\ref{rec1} reveals,
$C_f$ is equivalent to $f'_Y(\OO_Y)$ (note that $f_Y'$ is a
lax symmetric monoidal functor, and $f'_Y(\OO_Y)$
lies in $\CAlg(Y\times_kX)$).
\end{Remark}

\begin{Remark}
Suppose that $X$ and $Y$ are quasi-projective varieties over $k$.
Then the above argument also works for $X$ and $Y$,
and we have an equivalence
\[
\textup{Cor}(Y,X) \simeq \Map_{\CAlg(\PR_k)}(\QC^\otimes(X),\QC^\otimes(Y)).
\]
It has been proved in \cite{FI} that $\textup{Cor}(Y,X)$
is naturally equivalent to
$\Map_{\Sh(\Aff_k)}(Y,X)$. In particular, every correspondence is a graph of a morphism.
\end{Remark}

{\it Proof of Theorem~\ref{recmain} (ii).}
Identifying $\QC^\otimes(Z)\otimes_k\QC^\otimes(X)$ and $\QC^\otimes(Z)\otimes_k\QC^\otimes(Y)$ with $\QC^\otimes(Z\times_kX)$ and
$\QC^\otimes(Z\times_kY)$ respectively (Proposition~\ref{Moritaproduct}),
we have the diagram
\[
\xymatrix{
 \QC^\otimes(X) \ar[r]^f \ar[d]_{p_X^*} & \QC^\otimes(Y) \ar[d]^{p_Y^*} \ar[rd]^g &  \\
\QC^\otimes(Z)\otimes_k\QC(X) \ar[r]^{\textup{id}_Z\otimes f} & \QC^\otimes(Z)\otimes_k\QC^\otimes(Y) \ar[r]^(0.6){g_Z} & \QC^\otimes(Z) \\
\QC^\otimes(Z) \ar[u]^{p_Z^*} \ar[ur] \ar[urr]_{\textup{id}_Z} &  & 
}
\]
where $g_Z$ is determined by $g$
and $\textup{id}_Z$.
For a left adjoint functor $F$, we write $F'$ for a right adjoint of $F$.
Then by Remark~\ref{exconst}, $C_g\simeq g_Z'(\OO_Z)$.
Note that
$g_Z\circ (\textup{id}_Z\otimes f) \simeq (g\circ f)_Z$
where $(g\circ f)_Z$ is determined by $g\circ f$ and $\textup{id}_Z$.
Thus $(\textup{id}_Z\otimes f)'(C_g)\in \CAlg(X\times_kZ)$ corresponds to $g\circ f$.
It will suffice to prove that
$(\textup{id}_Z\otimes f)'(C_g)$
is equivalent to 
$(p_{ZX})_*(p_{ZY}^*(C_g)\otimes p_{YX}^*(C_f))$.
Unwinding the construction of $f$ obtained from $C_f$,
we see that $f$ is the composite
\[
\QC^\otimes(X)\stackrel{p_X^*}{\to} \QC^\otimes(Y\times_kX)\stackrel{C_f\otimes(-)}{\to} \Mod_{C_f}^\otimes(\QC(Y\times_kX))\stackrel{\sim}{\leftarrow} \QC^\otimes(Y).
\]
Therefore,
the right adjoint of $\textup{id}_Z\otimes f$
is
the composite
\begin{eqnarray*}
\QC(Z)\otimes\QC(Y)\stackrel{\sim}{\to}   \QC(Z)\otimes\Mod_{C_f}(\QC(Y\times_kX)) &\simeq& \Mod_{p^*_{YX}(C_f)}(\QC(Z\times_k Y\times_kX))  \\
&\stackrel{\textup{forget}}{\to} & \QC(Z\times_kY\times_kX) \\ &\stackrel{(p_{ZX})_*}{\to}& \QC(Z\times_kX).
\end{eqnarray*}
The image of $C_g$ under the composite is $(p_{ZX})_*(p_{ZY}^*(C_g)\otimes p_{YX}^*(C_f))$, as desired.
\QED

\begin{Remark}
\label{kurukuru}
Theorem~\ref{recmain} (i) and (ii) implies that
if $U\in \Cor(Y,X)$ and  $V\in \Cor(Z,Y)$,
$V\star U$ lies in $\Cor(Z,X)$.
One can also prove it by verifying the definition directly.
\end{Remark}

\section{Fine $\infty$-categories and Examples}
\label{EX}

In this Section, we give some examples and applications.
For this purpose, we start with some usable results.

\subsection{Elementary properties}

\begin{Proposition}
\label{useful}
Let $\CCC^\otimes$ be a symmetric monoidal idempotent complete
additive category (the tensor product is additive
separately in each variable).
Suppose that the endomorphism algebra of a unit of $\CCC$ is
a field $K$ of characteristic zero (hence $\CCC$ is $K$-linear).
Let $C$ be a nonzero dualizable object in $\CCC$ and
suppose that the $(n+1)$-fold wedge-product $\wedge^{n+1}C$ is a zero
object.
Suppose that $n$ is the minimal natural number
such that $\wedge^{n+1}C$ is a zero
object.
Then $\wedge^nC$ is invertible, i.e., $(\wedge^nC)\otimes (\wedge^nC)^\vee
\simeq (\wedge^nC)^\vee\otimes (\wedge^nC)$ is a unit for some object $(\wedge^nC)^\vee$.
In particular, $C$ is a $n$-dimensional wedge-finite object.
\end{Proposition}

\Proof
Let $\chi(C)$ be the trace defined as an element of $K:=\Hom_{\CCC}(\uni_\CCC,\uni_\CCC)$ given by
\[
\uni_\CCC\to C^\vee\otimes C\stackrel{\textup{flip}}{\simeq} C\otimes C^\vee\to \uni_\CCC
\]
where the left map is the coevaluation and the right map is the evaluation.
Taking account of $\wedge^{n+1}C\simeq 0$ we see
by \cite[Lemma 4.16, Corollary 4.20]{Ivo}
that $\chi(C)=n\in \ZZ\subset K$.
Since $\chi(\wedge^nC)=\frac{1}{n!}\chi(C)(\chi(C)-1)\cdots (\chi(C)-n+1)$,
we have $\chi(\wedge^n C)=1$.
Combining $\chi(\wedge^n C)=1$
and \cite[Corollary 3.16, Corollary 4.20]{Ivo} gives
$\wedge^2(\wedge^n(C))\simeq 0$.
Then according to \cite[8.2.9]{Kim2} and our hypothesis that 
the endomorphism algebra of the unit has no non-trivial idempotents,
$\wedge^nC$ is invertible.
\QED

\begin{Remark}
In Proposition~\ref{useful}, if one drops the assumption on the
endomorphism algebra of the unit, then the assertion does not hold.
Namely, one can not deduce that $C$ is wedge-finite
from the condition
that $C$ is dualizable and $(n+1)$-fold wedge-product $\wedge^{n+1}C$ is zero
for some $n$.
Let $X=\Spec A \sqcup \Spec B$ is a non-connected usual affine scheme
and let $L$ be an $\OO_X$-module which is an invertible sheaf on
$\Spec A$ and is zero on $\Spec B$. Then $L$ is dualizable in the symmetric monoidal category of $\OO_X$-modules and $\wedge^2L\simeq 0$, but it
is not an invertible object in the symmetric monoidal category of $\OO_X$-modules.
\end{Remark}

\begin{Proposition}
\label{affinefine}
Let $\XX$ be a derived stack $\CAlg_k\to \widehat{\mathcal{S}}$ such that
$\QC^\otimes(\XX)$ is a fine $\infty$-category.
Let $\YY$ be another sheaf and $f:\YY\to \XX$ an affine morphism,
i.e., for any $\Spec A\to \XX$ the fiber product
$\YY\times_{\XX}\Spec A$ is affine.
Then $\QC^\otimes(\YY)$ is a fine $\infty$-category.
\end{Proposition}

\Proof
Let $\{V_{\lambda}\}_{\lambda\in \Lambda}$
be a set of wedge-finite objects such that
$\{V_\lambda,V_\lambda^\vee\}_{\lambda\in \Lambda}$
generates $\QC^\otimes(\XX)$ as a symmetric monoidal
stable presentable $\infty$-category.
Note that each wedge-finite object $p^*(V_\lambda)$
is compact.
Indeed, unwinding the definition of $\QC(\XX)$ and $\QC(\YY)$
and using the base change formula for affine morphisms,
we may assume that $\XX$ and $\mathcal{Y}$ are affine.
The fact that $f_*$ preserves all small colimits implies that
\[
\Map_{\QC(\YY)}(f^*(V_\lambda),\varinjlim_i M_i)\simeq \Map_{\QC(\XX)}(V_\lambda,f_*(\varinjlim_i M_i))\simeq \varinjlim_i\Map_{\QC(\XX)}(V_\lambda, f_*(M_i))
\]
for any filtered colimit $\varinjlim_i M_i$.
It also 
follows that the unit of $\QC^\otimes(\YY)$ is compact.
In addition, since $f_*$ is conservative, we deduce
from Remark~\ref{homlevelcpt}
that the set $\{f^*(V_\lambda),f^*(V_\lambda)^\vee\}_{\lambda\in \Lambda}$
of compact objects generates $\QC(\YY)$ as a symmetric monoidal stable
presentable $\infty$-category.
Hence $\QC^\otimes(\YY)$ is fine.
\QED

\begin{Proposition}
Let $\CCC^\otimes$ and $\DDD^\otimes$ be two fine $\infty$-categories.
Then $\CCC^\otimes\otimes_k \DDD^\otimes$ is also fine.
\end{Proposition}

\Proof
Combine Theorem~\ref{main1} and Proposition~\ref{Moritaproduct}.
\QED

\subsection{Classical Tannakian categories}
We discuss a relationship with (classical) neutral Tannakian categories.
Let $G$ be an algebraic group over a field $k$
of characteristic zero.
Let $\QC^\otimes(BG)$ be the $k$-linear stable
presentable $\infty$-category of quasi-coherent complexes over $BG$.

Let us observe that $\QC^\otimes(BG)$ is a fine $\infty$-category.
The underlying $\infty$-category
$\QC(BG)$ is compactly generated, and compact and dualizable objects
coincide (see for example \cite[Corollary 3.22]{BFN}: $\QC^\otimes(BG)$ in this paper
agrees with that of {\it loc. cit.}).
We 
take a closed immersion $G\hookrightarrow \GL_r$ that makes $G$
a subgroup scheme of $\GL_r$.
Furthermore, by \cite[Lemma 3.1]{Tot}
we can choose $G\hookrightarrow \GL_r$
so that the quotient $\GL_r/G$ is quasi-affine over $k$.
The morphism $p:BG\to B\textup{GL}_r$
induced by $G\hookrightarrow \GL_r$
is quasi-affine since $GL_r/G$ is a usual quasi-affine scheme
(in particular, the structure sheaf is very ample).
Let $V$ be the standard representation of
$\GL_r$.
Then by the standard use of the adjoint pair $(p^*,p_*)$
(see the proof of \cite[Proposition 3.21]{BFN}),
the set
$\{p^*(V),p^*(V)^\vee\}$
generates $\QC(BG)$ as a symmetric monoidal stable presentable $\infty$-category.
Note that $p^*(V)$
is compact and dualizable.
Recall that $V$ is wedge-finite and so is $p^*(V)$.
Therefore we conclude:

\begin{Proposition}
\label{alggp}
If $G$ is an algebraic group over $k$, then
$\QC^\otimes(BG)$ is a fine algebraic $\infty$-category.
\end{Proposition}

\begin{Corollary}
\label{nonred}
Let $[\Spec A/G]$ be a derived quotient stack,
where an (possibly  non-reductive) algebraic group $G$ (over $k$) acts on $\Spec A$ with $A\in \CAlg_k$.
Then $\QC^\otimes([\Spec A/G])$ is a fine algebraic $\infty$-category.
\end{Corollary}

\Proof It follows from
Proposition~\ref{affinefine} and~\ref{alggp}.
\QED

\begin{Remark}
Let $\CCC^\otimes$ be a $k$-linear symmetric monoidal stable
presentable $\infty$-category.
By Corollary~\ref{nonred}, the conditions
in Theorem~\ref{algebraic} are also equivalent 
to the condition: $\CCC^\otimes$ is equivalent to 
$\QC^\otimes([\Spec A/G])$ for some $[\Spec A/G]$ such that an (possibly  non-reductive) algebraic group $G$ (over $k$) acts on $\Spec A$ with $A\in \CAlg_k$.

\end{Remark}

\begin{Remark}
\label{alggpR}
By Theorem~\ref{main1}, if $G$ is a pro-reductive group, $\QC^\otimes(BG)$
is a fine $\infty$-category.
For an arbitrary pro-algebraic group $G$ over $k$,
$\QC^\otimes(BG)$ is not necessarily fine (since the unit is not compact
when $G$ has infinite cohomological dimension).
For our purpose
a correct generalization of $\QC^\otimes(BG)$ to arbitrary pro-algebraic
groups is
given by the Ind-category $\Ind^\otimes(\textup{Coh}(BG))$,
where $\textup{Coh}(BG)$ is the stable subcategory of $\QC(BG)$
spanned by dualizable objects.
Namely, it is the symmetric monoidal
compactly generated stable $\infty$-category
of {\it Ind-coherent complexes} on $BG$.
Note that for a pro-algebraic group $G$,
a finite dimensional representation is
a wedge-finite object in $\Ind^\otimes(\textup{Coh}(BG))$.
Thus, $\Ind^\otimes(\textup{Coh}(BG))$
is a fine $\infty$-category because the set of finite dimensional representations
of $G$ generates $\Ind^\otimes(\textup{Coh}(BG))$
as a stable presentable $\infty$-category, and objects
in $\textup{Coh}(BG)$ are compact
in $\Ind^\otimes(\textup{Coh}(BG))$.
\end{Remark}

\subsection{Stable $\infty$-category of mixed motives, fine $\infty$-categories and Kimura finiteness}
\label{exmot}
We study a relationship between
 fine $\infty$-categories, the symmetric monoidal
stable
$\infty$-category of mixed motives, and Kimura finiteness of Chow
motives. We also consider the $\infty$-category of noncommutative motives.

We begin by briefly recalling its background; the 
reason being that we would like to
regard the category of {\it mixed motives}
as a fine $\infty$-category.
One of the main themes of motives is a motivic Galois theory which generalizes
the classical Galois theory of fields (see e.g. \cite{A}).
The Galois theory for Artin motives corresponds to the classical Galois
theory (cf. \cite{A}, \cite[Section 8]{Bar}). The Galois group of motives (motivic Galois group)
should encode the structure of periods of motives.
A conjectural abelian category of mixed motives is expected to
be a Tannakian category. Furthermore,
it has been conjectured by Beilinson and Deligne
that {\it the} abelian category of mixed motives should be the 
heart of a conjectural so-called motivic $t$-structure in the triangulated category
of mixed motives $DM$ that was constructed by Hanamura, Levine and Voevodsky.  
But the existence of a motivic $t$-structure is inaccessible as of this moment
(except the case of mixed Tate motives). 


With this in mind, we study an $\infty$-categorical enhancement of $DM$
by means of fine $\infty$-categories.
That is, the above conjectural line and derived Tannakian viewpoint
suggest the following picture.
The $\infty$-categorical enhancement of $DM$ should constitute
a fine $\infty$-category, i.e., our $\infty$-categorical analogue of Tannakian
category.
A motivic Galois group should appear as
a group object arising from a pointed derived stack corresponding to
the fine $\infty$-category equipped with a realization functor associated with a Weil cohomology theory.
In \cite{Tan},
we constructed derived automorphism group schemes of realization functors of mixed motives associated with a Weil
cohomology theory (i.e. motivic Galois groups) by means of {\it tannakization},
and proved a consistency with the above 
traditional line (but, the general construction
of derived automorphism group schemes in \cite{Tan} is somewhat abstract, whereas foundational properties are proved).

{\it Stable $\infty$-category of mixed motives.}
Now let us consider the $\QQ$-linear
symmetric monoidal (stable) presentable $\infty$-category
$\DM^\otimes$ of (Voevodsky's) mixed motives over a perfect field $S=\Spec K$,
which is treated in \cite{Tan}, \cite{Bar}, \cite{PM}, \cite{Rob}
(see these papers for further details).
Here, we use the symmetric
monoidal model category $\textup{DM}^\otimes$
studied in \cite[Example 7.15]{CD1}
and let $\DM^\otimes$ be the stable presentably symmetric monoidal
$\infty$-category (i.e., an object
of $\CAlg(\PR_{\mathbb{S}})$)
obtained from (the full subcategory of cofibrant objects in)
$\textup{DM}^\otimes$ by inverting
weak equivalences.
For a smooth variety $X$, i.e., a smooth scheme separated
of finite type over $S$,
there is a motive $M(X)$ of $X$ that is an object of $\DM$.
We will work with $\QQ$-coefficients.
Namely, $\textup{DM}$
consists of symmetric
spectrum objects, with respect to a Tate twist, in the category
of chain complexes of Nisnevich sheaves of $\QQ$-vector spaces
(with transfers) on the category of finite correspondences
over $S$
(see \cite{CD1}, \cite{MVW}).
As a result, $\DM^\otimes$ is a
$\QQ$-linear symmetric monoidal presentable $\infty$-category (see \cite[Section 5]{Tan}
for more details).
We can consider a direct generalization (of this subsection)
to relative mixed motives over a smooth variety $S$, but
for simplicity, we consider the case when $S$ is the Zariski spectrum of
a perfect field.

{\it Chow
motives}. There is a symmetric monoidal
$\QQ$-linear (ordinary) category $CHM^\otimes$ of the (homological)
Chow motives (cf. \cite{Scho}, see also \cite[4.1]{PM} for homological convention).
In $CHM$, every object is dualizable.
For a projective smooth variety $X$ over $K$, there exists a Chow
motive $h(X)$ in $CHM$.
Moreover,
there is
a
symmetric monoidal $\QQ$-linear fully faithful functor $CHM\to \textup{h}(\DM)$
which carries $h(X)$ to $M(X)$ (cf. \cite[20.2]{MVW}).
Hence Chow motives can be regarded as objects in $\DM$.

{\it Kimura finiteness of Chow motives}.
The work of Kimura \cite{Kim} and others places
Kimura finiteness at the heart of recent developments of motivic theory.
Let us recall this notion.
An object $M$ in the underlying category
$CHM$ is evenly finite dimensional
(resp. oddly finite dimensional)
if there is a non-negative integer $n$
such that $\wedge^nM=0$ (resp. $\Sym^nM=0$).
Here $\Sym^nM$ denotes the symmetric product $\Ker(1-\frac{1}{n!}\Sigma_{\sigma\in \Sigma_n}\sigma)$
where $\frac{1}{n!}\Sigma_{\sigma\in \Sigma_n}\sigma:M^{\otimes n}\to M^{\otimes n}$ is the symmetrizer.
An object $M$ in $CHM$ is Kimura finite dimensional
if there exists a decomposition
$M\simeq M^+\oplus M^-$ such that $M^+$ is evenly finite dimensional
and $M^-$ is oddly finite dimensional.
Similarly, we say that an object $M$ in $\DM$ is Kimura finite dimensional
if 
there exists a decomposition
$M\simeq M^+\oplus M^-$ such that $M^+$ is evenly finite dimensional
and $M^-$ is oddly finite dimensional in the homotopy category.

\begin{Lemma}
\label{Kimstep}
If $M^+$ is an evenly (resp. $M^-$ is an oddly)
finite dimensional object 
in
$\DM$,
then
$M^+[2m]$ (resp. $M^-[2m+1]$)
is wedge-finite for any $m\in \ZZ$.
In particular, if $M$ is a Kimura finite dimensional Chow motif
such that $M\simeq M^+\oplus M^-$
where $M^+$ is evenly finite dimensional and $M^-$
is oddly finite dimensional,
then $M^+[2m]\oplus M^-[2n+1]$ is wedge-finite
for any $m,n\in \ZZ$.
\end{Lemma}

\Proof
The endomorphism algebra of a unit of $\textup{h}(\DM)$
is $\QQ$. Indeed, the endomorphism algebra of the unit $M(\Spec K)$
is
$\Hom_{\hhh(\DM)}(M(\Spec K),M(\Spec K))\simeq \textup{CH}^0(\Spec K\times_{\Spec K}\Spec K)\otimes_{\ZZ}\QQ \simeq \textup{CH}^0(\Spec K)\otimes_{\ZZ}\QQ\simeq\ZZ\otimes_{\ZZ} \QQ\simeq \QQ$,
where $CH^i(-)$ denotes the Chow group (cf. \cite[14.5.6]{MVW}).
Thanks to Proposition~\ref{useful}, $M^+$ is wedge-finite.
By the Koszul sign rule (cf. \cite[8A.2]{MVW}),
$M^+[1]$ is oddly finite dimensional.
Similarly, if $M^-$ is oddly finite dimensional,
then $M^-[1]$ is evenly finite dimensional,
and so it is wedge-finite. Now our assertion is clear.
\QED

Let $KF$ be a small set of objects in $\DM$ that consists of
Kimura finite dimensional Chow motives.
(We remark that if $M$ is Kimura finite dimensional,
then the dual object $M^\vee$ is Kimura finite dimensional.)
Let $T$ be a subset of $KF$.
Namely, any element of $T$ is Kimura finite dimensional.
Let $\DM^\otimes \langle T \rangle$
be a symmetric monoidal stable presentable full subcategory 
of $\DM$
generated by $\{M, M^\vee\}_{M\in T}$ as a symmetric monoidal stable presentable $\infty$-category.
That is, it is the smallest stable subcategory
which contains the unit and $\{M, M^\vee\}_{M\in T}$ and is closed under small
coproducts and tensor products.
(We note that a dualizable object in $\DM^\otimes \langle T \rangle$
is not necessarily Kimura finite.)
Known examples of Kimura finite objects include Chow motives $h(X)$ of abelian varieties (and more generally abelian schemes), some algebraic surfaces
(rational surfaces, K3 surfaces of certain types, Godeaux surfaces..), Fano 3-folds,
Tate objects $\QQ(n)$, Artin motives, and forth.
We then have

\begin{Theorem}
\label{motuncond}
The $\QQ$-linear symmetric monoidal
presentable $\infty$-category $\DM^\otimes \langle T \rangle$
is a fine $\infty$-category. Namely, there exists a derived stack
$[\Spec A/G]$, where $G$ is a pro-reductive group over $\QQ$,
and an equivalence
\[
\DM^\otimes \langle T \rangle \simeq \QC^\otimes([\Spec A/G]).
\]
If $T$ is a finite set, then $\DM^\otimes \langle T \rangle$
is a fine algebraic $\infty$-category.
\end{Theorem}

\Proof
Note first that dualizable and compact objects coincide
in $\DM$ (see \cite[Theorem 2.7.10]{CD2}).
In addition, if $X$ is a smooth projective variety, $M(X)$
is dualizable.
Lemma~\ref{Kimstep} implies that $\DM^\otimes \langle T \rangle$ admits
a small set of wedge-finite objects which generates
$\DM^\otimes \langle T \rangle$ as a symmetric monoidal stable presentable $\infty$-category
(consider $M^+[2m]\oplus M^-[2n+1]$).
Hence $\DM^\otimes \langle T \rangle$ is a fine $\infty$-category.
The existence of $[\Spec A/G]$ follows from Theorem~\ref{main1}.
The final assertion follows from Remark~\ref{finiteproduct}.
\QED

\begin{Remark}
It is important to notice that the existence of a motivic $t$-structure
of $\DM^\otimes \langle KF \rangle$ is still unknown and mysterious, but Theorem~\ref{algebraic} and~\ref{main1} are applicable.
For known cases of Kimura finiteness, Theorem~\ref{motuncond}
provides an unconditional application, which is a far-reaching
generalization of the mixed Tate case.
A precursor to the above Theorem in the case of mixed Tate motives
is a theorem of Spitzweck (see \cite{Spi}).
The symmetric monoidal $\infty$-category $\DM^\otimes \langle KF \rangle$ unconditionally contains the important class of
mixed motives; mixed motives generated by abelian schemes (as a
symmetric monoidal presentable $\infty$-category).

The statement of the above form seems to be somewhat abstract.
But, thanks to
Proposition~\ref{warm} and~\ref{warm2}
we have an explicit presentation of the underlying complex of
$A$  by means of motivic complexes,
Weyl construction and the product of general linear groups.
We note that this presentation depends on the choice of a set of
wedge-fine generators $\{C_{\lambda}\}_{\lambda\in \Lambda}$
that appears in Definition~\ref{fine}.
For various applications (see the next Remark),
it would be nice to have
$\{C_{\lambda}\}_{\lambda\in \Lambda}$
such that each $R(C_{\lambda})$ belongs to the heart of
the standard $t$-structure of $\Mod_{k}$
(i.e., the concentrated in degree zero) where
$R:\DM^\otimes \to \Mod_{k}^\otimes$ is a realization functor (e.g.,
\'etale,  Betti, de Rham realizations).
In all known Kimura finite cases at the writing of this paper, fortunately,
one can take such sets of wedge-finite generators.
\end{Remark}

\begin{Remark}
\label{motGalcon}
Theorem~\ref{motuncond} and variants can
be applied to explicit constructions and studies
of motivic Galois groups of $\DM^\otimes \langle T \rangle$
by means of the construction of based loop spaces
(equivariant bar construction) of $[\Spec A/G]$.
Let
\[
R:\QC^\otimes([\Spec A/G])\simeq \DM^\otimes \langle T \rangle \to \Mod_{k}^\otimes
\]
be a realization functor associated to mixed Weil (co)homology with coefficients in $k$ (see e.g. \cite{Tan}).
Suppose that each $R(C_{\lambda})$ belongs to the heart of
the standard $t$-structure of $\Mod_{k}$ for a
prescribed set $\{ C_\lambda\}$
of wedge-finite objects.
Then as discussed in Remark~\ref{Tandual} it gives rise to a morphism (``geometric point'')
\[
p:\Spec k\to [\Spec A/G]
\]
and the realization functor
$R$ can be identified with the pullback functor $p^*$.
From this, we have the based loop space $\Omega_*[\Spec A/G]=\Spec k\times_{[\Spec A/G]}\Spec k$ that is a derived affine group scheme; similar constructions yield the Betti-de Rham comparison torsor,
and motivic Galois group representing the automorphism group of the realization functor (see \cite{Bar}).
(This construction can be
generalized to the context of realization of relative mixed motives.)
The interested reader is referred to \cite{PM} and \cite{Bar}
for detailed study and further applications to mixed motives.
\end{Remark}

It is natural to expect

\begin{Conjecture}
\label{myconj}
The $\QQ$-linear symmetric monoidal
stable presentable $\infty$-category $\DM^\otimes$ 
is a fine $\infty$-category.
\end{Conjecture}

Recall the following well-known conjecture:

\begin{Conjecture}[Kimura, O'Sullivan]
\label{Kimuraconj}
Every object in $CHM$
is Kimura finite dimensional.
\end{Conjecture}

The conjecture of Kimura and O'Sullivan
does not imply the existence of a motivic $t$-structure on $DM$,
but we have the following nice implication:

\begin{Proposition}
Conjecture~\ref{Kimuraconj} implies Conjecture~\ref{myconj}.
\end{Proposition}

\Proof
Let $CHM\hookrightarrow \textup{h}(\DM)$
be the canonical fully faithful functor.
The essential image of this functor
generates $\DM$ as a stable presentable $\infty$-category (cf. \cite[2.7.10]{CD2}). Moreover, the unit is comapct (see Theorem~\ref{motuncond}).
By Conjecture~\ref{Kimuraconj}, every Chow motive $M$ has
a decomposition $M\simeq M^+\oplus M^-$
such that $M^+$ is evenly finite dimensional and $M^-$
is oddly finite dimensional.
By Lemma~\ref{Kimstep}, both $M^+$ and $M^-[1]$ are wedge-finite.
Hence our claim follows.
\QED

{\it Noncommutative motives.}
Next we consider
noncommutative motives.
Here
we use the theory developed by Blumberg, Gepner
and Tabuada \cite{BGT}, \cite{BGT2}.
Let $\mathcal{M}_{\textup{add}}^\otimes$ be the stable
presentably
symmetric monoidal
$\infty$-category of noncommutative motives
(we use the same notation as in
\cite[Section 6]{BGT}, \cite[Section 5]{BGT2}).
Let $\SP^\otimes$ be
the stable presentably symmetric monoidal
$\infty$-category of spectra.
Let $s:\SP^\otimes\to \mathcal{M}_{\textup{add}}^\otimes$
be a symmetric monoidal colimit-preserving
functor, which is unique up to a contractible space of choices.
We define $\mathcal{M}_{\textup{add},\QQ}^\otimes$
to be $\Mod_{s(H\QQ)}^\otimes(\mathcal{M}^\otimes_{\textup{add}})$,
that is, the ``$\QQ$-linearization'' of $\mathcal{M}^\otimes_{\textup{add}}$
($H\QQ$ is the Eilenberg-MacLane spectrum).
There is a canonical symmetric monoidal colimit-preserving
functor $\Mod^\otimes_{\QQ}\to \mathcal{M}_{\textup{add},\QQ}^\otimes$
that is induced by $s$.
We remark
that by the Morita theory (e.g. \cite[Proposition 4.1 (1)]{BFN}),
one can naturally identify $\mathcal{M}_{\textup{add},\QQ}^\otimes$
with $\Mod^\otimes_{\QQ}\otimes_{\SP^\otimes} \mathcal{M}^\otimes_{\textup{add}}$. The Kimura finiteness in $\mathcal{M}_{\textup{add},\QQ}^\otimes$ is
defined in a similar way: An object $M$
of the homotopy category
$\hhh(\mathcal{M}_{\textup{add},\QQ})$ is said to be Kimura finite dimensional
if
there is a decomposition
$M\simeq M^+\oplus M^-$ such that $\wedge^nM^+\simeq 0$ and
$\textup{Sym}^nM^-\simeq 0$ for a sufficiently large natural number $n>0$.

\begin{Proposition}
\label{noncommutativemot}
Let $T=\{M_i\}_{i\in I}$ be a small set of dualizable objects 
of $\mathcal{M}_{\textup{add},\QQ}$ such that each $M_i$
is Kimura finite dimensional.
Let
$\mathcal{M}_{\textup{add},\QQ}\langle T \rangle$
be the stable subcategory of $\mathcal{M}_{\textup{add},\QQ}$
generated by $T'=\{M_i,M_i^\vee\}_{i\in I}$
as a symmetric monoidal stable
presentable $\infty$-category.
Namely, $\mathcal{M}_{\textup{add},\QQ}^\otimes \langle T \rangle$
is the smallest stable subcategory of
$\mathcal{M}_{\textup{add},\QQ}^\otimes$
which contains the unit object and $T'$, and is closed under
small coproducts
and tensor products.
Then $\mathcal{M}_{\textup{add},\QQ}^\otimes \langle T \rangle$ is a fine $\infty$-category.
\end{Proposition}

\Proof
We will prove that $\mathcal{M}_{\textup{add},\QQ}^\otimes \langle T \rangle$
satisfies two conditions (i), (ii) in Definition~\ref{fine}.
Since a unit object $\mathbf{1}$
in $\mathcal{M}_{\textup{add}}$ is compact
and the forgetful functor $v:\mathcal{M}_{\textup{add},\QQ}\to \mathcal{M}_{\textup{add}}$ preserves filtered colimits, it follows that
a unit object $\mathbf{1}_{\QQ}:=\mathbf{1}\otimes s(H\QQ)$ in $\mathcal{M}_{\textup{add},\QQ}$
is also compact. (The compactness of $\mathbf{1}$ in $\mathcal{M}_{\textup{add}}$ follows from the construction of $\mathcal{M}_{\textup{add}}$
and \cite[5.5.7.3]{HTT}:
The compactness of $\mathbf{1}$ is clearly stated and proved in
Proposition 9.22 of the version 3 of \cite{BGT} on arXiv).
Thus, the unit object $\mathbf{1}_{\QQ}$
is compact in $\mathcal{M}_{\textup{add},\QQ} \langle T \rangle$.
For any $M_i\in T$,
there is a decomposition
$M_i^+\oplus M_i^-$ such that $\wedge^n M_i^+\simeq 0$
and $\wedge^n (M_i^-[1])\simeq (\Sym^nM_i^-)[n] \simeq 0$
for a sufficiently large natural number $n>0$.
We will show that both $M_i^+$ and $M_i^-[1]$ are wedge-finite.
In view of Proposition~\ref{useful},
it is enough to show that the endomorphism algebra
of the unit object $\mathbf{1}_{\QQ}$ in the homotopy category
$\hhh(\mathcal{M}_{\textup{add},\QQ})$ is $\QQ$.
By a presentation of the mapping spectrum in terms of $K$-theory
\cite[Theorem 7.13]{BGT}, $\Hom_{\hhh(\mathcal{M}_{\textup{add}})}(\mathbf{1}[n],\mathbf{1}) =\pi_n(K(\mathbb{S}))=K_n(\mathbb{S})$ for $n\in \ZZ$ ($K_n(\mathbb{S})=0$ if $n<0$).
Here $K(-)$ is the connective K-theory spectrum and $\mathbb{S}$
is the sphere spectrum.
By \cite[5.19]{BGT2} (or a simple application of Proposition~\ref{PMA} to
$s:\SP^\otimes \to \mathcal{M}^\otimes_{\textup{add}}$),
there
is a (symmetric monoidal)
fully faithful functor
$\Mod_{K(\mathbb{S})}\hookrightarrow
 \mathcal{M}_{\textup{add}}$
whose essential image is the smallest stable subcategory
which contains the unit object and is closed under small coproducts.
Consequently, we have equivalences
\[
\Map_{\mathcal{M}_{\textup{add},\QQ}}(\mathbf{1}_{\QQ},\mathbf{1}_{\QQ})\simeq 
\Map_{\mathcal{M}_{\textup{add}}}(\mathbf{1},v(\mathbf{1}_\QQ))\simeq \Map_{\Mod_{K(\mathbb{S})}}(K(\mathbb{S}),K(\mathbb{S})\otimes_{\mathbb{S}}H\QQ)\simeq \Map_{\SP}(\mathbb{S},K(\mathbb{S})\otimes_{\mathbb{S}}H\QQ).
\]
In particular, $\Hom_{\hhh(\mathcal{M}_{\textup{add},\QQ})}(\mathbf{1}_\QQ,\mathbf{1}_\QQ)=\pi_0(\Map_{\mathcal{M}_{\textup{add},\QQ}}(\mathbf{1}_{\QQ},\mathbf{1}_{\QQ}))\simeq \pi_0(K(\mathbb{S})\otimes_{\mathbb{S}}H\QQ)$.
Since $K(\mathbb{S})$ and $H\QQ$ are connective
and $K_0(\mathbb{S})\simeq K_0(\ZZ)\simeq \ZZ$, we see that
$\pi_0(K(\mathbb{S})\otimes_{\mathbb{S}}H\QQ)\simeq \pi_0(K(\mathbb{S}))\otimes_{\ZZ}\QQ\simeq \QQ$.
\QED

\subsection{Quasi-coherent complexes on an algebraic variety}

We will apply our duality theorem to the derived $\infty$-category
of quasi-coherent sheaves on a quasi-projective variety.
Let $X$ be a quasi-projective scheme over a field $k$.
Note that $X$ admits a Zariski covering $\sqcup_{1\le i \le n}\Spec A_i\to X$
and its Cech nerve gives rise to a groupoid object $X_\bullet:\NNNN(\Delta)^{op}\to \Aff_k$.

Let $\QC^\otimes(X)$ be the $k$-linear symmetric monoidal
$\infty$-category of quasi-coherent complexes on $X$,
that is, $\QC^\otimes(X):=\varprojlim\QC^\otimes(X_\bullet([n]))$.
Let $\DDD_{qc}(X)$ be the derived $\infty$-category of (ordinary)
$\OO_X$-modules whose cohomology is quasi-coherent on $X$
(cf. \cite[1.3.5.8]{HA}).
We then remark that there is an equivalence
$\QC(X)\simeq \DDD_{qc}(X)$ (indeed, by \cite[2.1.8, 2.3.1]{DAG8}
there is an equivalence $\QC(X)^+\simeq \DDD^+_{qc}(X)$
between the full subcategories of left bounded objects
with respect to the ``standard'' $t$-structures, and thus
the left completeness of $\DDD_{qc}(X)$ and $\QC(X)$
\cite[B1]{HNR}, \cite[2.3.18]{DAG8}
implies $\QC(X)\simeq \DDD_{qc}(X)$).

\begin{Theorem}
\label{dSerre}
Suppose that $k$ is of characteristic zero.
The $k$-linear symmetric monoidal presentable $\infty$-category
$\QC^\otimes(X)$ is a fine $\infty$-category,
so that there exist a derived stack $[\Spec A/\mathbb{G}_m]$ and an
equivalence
\[
\QC^\otimes(X)\simeq \QC^\otimes([\Spec A/\mathbb{G}_m])
\]
where $\mathbb{G}_m=\GL_1$.
Moreover, there is an equivalence
$A\simeq \oplus_{r\in \ZZ}\mathsf{Hom}_{\QC(X)}(\OO_X,\LLL^{\otimes r})\otimes \chi_r$ in $\QC(B\mathbb{G}_m)$ where $\chi_r$ is the character of weight $r$
of $\mathbb{G}_m$, and $\mathcal{L}$ is a very ample invertible sheaf.
\end{Theorem}

\Proof
Note first that $\QC(X)$
is compactly generated,
and dualizable and compact objects coincide (cf. \cite{BFN}).
Moreover, if $\LLL$ is a very ample invertible sheaf on $X$
a single compact object
$\oplus_{0\ge i \ge -d}\LLL^{\otimes i}$ for some $d\ge0$
generates $\QC(X)$ as a stable presentable $\infty$-category (see \cite[Theorem 4]{Orl},\cite[Lemma 3.2.2]{Van}).
It follows that $\{\LLL, \LLL^\vee\}$ generates
$\QC^\otimes(X)$
as a symmetric monoidal stable presentable $\infty$-category.
Note that $\LLL^{\vee}$ is wedge-finite and $1$-dimensional.
Let $A\otimes\chi_r$ denote the image
of $\chi_r$
under the natural pullback functor
$\QC(B\mathbb{G}_m)\to \QC([\Spec A/\mathbb{G}_m])$.
Then by Theorem~\ref{algebraic} we obtain
a derived stack $[\Spec A/\mathbb{G}_m]$ and
an equivalence
$\QC^\otimes([\Spec A/\mathbb{G}_m]) \stackrel{\sim}{\to} \QC^\otimes(X)$
in $\CAlg(\PR_k)$
which carries $A\otimes \chi_r$ to $\LLL^{\otimes (-r)}$.
Therefore by Proposition~\ref{warm}
$A\simeq \oplus_{r\in \ZZ}\mathsf{Hom}_{\QC(X)}(\OO_X,\LLL^{\otimes r})\otimes \chi_r$
in $\QC(B\mathbb{G}_{m})$, where $\mathsf{Hom}_{\QC(X)}(-,-)$ denote
the hom complex.
The truncation is given by $\pi_0(A)\simeq \oplus_{r\in \ZZ} H^0(X,\LLL^{\otimes r})\otimes \chi_r$.
\QED

Recall Serre's theorem which identifies the category
of coherent sheaves on a projective variety $X$
with the category of quasi-finitely generated graded modules of $\oplus_{r\in \ZZ} H^0(X,\LLL^{\otimes r})$ modulo torsion sheaves (see e.g. \cite[Ex. 5.8]{Har}).
We think of Theorem~\ref{dSerre} as a {\it derived analogue of Serre's
theorem}. In spite of
the equivalence $\QC^\otimes(X)\simeq \QC^\otimes([\Spec A/\mathbb{G}_m])$,
$[\Spec A/\mathbb{G}_m]$ is not equivalent to $X$ in general.

\subsection{Coherent complexes on a topological space and Rational homotopy theory}
\label{exRHT}

We will discuss the $\infty$-category of Ind-coherent complexes on a topological space from a viewpoint of rational homotopy theory. We work with
coefficients in a fixed base field $k$ of characteristic zero.

In his foundational work \cite{SGA1} Grothendieck developed the theory of Galois categories.
The category $\textup{Cov}(S)$
of finite topological covers of a topological space $S$
is a Galois category. A base point $s$ of $S$
determines a symmetric monoidal functor $f:\textup{Cov}(S)\to \textup{Fin}$
to the category of finite sets (with respect to
cartesian monoidal structures), which carries a cover $\phi:X\to S$
to $\phi^{-1}(s)$. The automorphism group
of $f$ is equivalent to the pro-finite completion $\hat{\pi}_1(S,s)$ of 
the fundamental group $\pi_1(S,s)$.
Moreover, $\hat{\pi}_1(S,s)$ continuously acts on the fiber
$\phi^{-1}(s)$.
It gives rise to a categorical equivalence between
$\textup{Cov}(S)$ and the category of finite sets endowed with continuous $\hat{\pi}_1(S,s)$-actions.
We will describe a generalization of this story to the context of
rational homotopy theory
by dint of fine $\infty$-categories.

Let $S$ be a connected topological space which we regard as an object in
$\mathcal{S}$.
We can think of $S$ as a constant sheaf $\Aff^{op}_k\to \mathcal{S}$
taking the value $S$.
Let $\QC^\otimes(S)$ denote the $k$-linear symmetric monoidal presentable $\infty$-category
of quasi-coherent complexes on $S$ (cf. Section~\ref{presection}).
If $S$ is a contractible space, $\QC^\otimes(S)$ is equivalent to
$\Mod_k^\otimes$.
For an arbitrary (small) topological space $S$,
$\QC^\otimes(S)$ is the limit $\varprojlim_{S}\Mod^\otimes_k$
of a constant diagram of $\Mod_k^\otimes$
indexed by the space $S$.
The underlying $\infty$-category $\QC(S)$ is nothing but (equivalent to)
the function complex $\Fun(S,\Mod_k)$.
We will use the $\infty$-category of Ind-coherent complexes
of $S$ instead of quasi-coherent complexes
since dualizable objects on $S$ are not necessarily compact objects.
Let us define the full subcategory of bounded coherent complexes.
Let $\Mod_{k,\ge 0}$ (resp. $\Mod_{k,\le0}$) be the
full subcategory of $\Mod_k$ that consists of objects $C$
such that $H_i(C)=0$ for $i<0$ (resp. $i>0$).
The pair $(\Mod_{k,\ge 0}, \Mod_{k,\le0})$ together with
the truncation functors $\tau_{\ge0 }:\Mod_k\to \Mod_{k,\ge 0}$,
$\tau_{\le0 }:\Mod_k\to \Mod_{k,\le 0}$
determines a $t$-structure on the stable $\infty$-categories $\Mod_k$.
The pair $(\Fun(S,\Mod_{k,\ge 0}),\Fun(S,\Mod_{k,\le0}))$ determines
a $t$-structure on $\Fun(S,\Mod_k)$.
To see this, we first note that
there are adjoint pairs
\[
\iota_{\ge0}^S:\Fun(S,\Mod_{k,\ge0}) \rightleftarrows \Fun(S,\Mod_k):\tau_{\ge0 }^S\ \ \textup{and}\ \ 
\tau_{\le0 }^S:\Fun(S,\Mod_k)\rightleftarrows \Fun(S,\Mod_{k,\le0}):\iota^S_{\le0}
\]
induced by the compositions with $\tau_{\ge0}$ and $\tau_{\le0}$.
The functors $\iota_{\ge0}^S$ and $\iota_{\le0}^S$ are canonical
fully faithful functors.
Let $s:\Delta^0 \to S$ be a point on $S$. The pullback
functor
$s^*:\QC(S)=\Fun(S,\Mod_k)\to \Fun(\Delta^0,\Mod_k)\simeq \Mod_k$
is induced by the composition with $s:\Delta^0 \to S$.
Note that the pullback functor $\Fun(S,\Mod_k)\to \Mod_k$
is conservative because $S$ is connected.
Namely, a morphism $f:M\to N$ in $\Fun(S,\Mod_k)$
is a zero map (i.e. null-homotopic) if and only if
$s^*(f)$ is a zero map in $\Mod_k$. Thus, if $M$ belongs to
$\Mod_{k,\ge0}$ and $N$ belongs to $\Mod_{k,\le0}$,
then any morphism $f:M\to N[-1]$ is a zero map
because $s^*(f)$ is a zero map.
Write $\tau^S_{\le -1}:=[-1]\circ \tau^S_{\le0}\circ [1]$.
Using adjoint pairs and the fact that
$s^*$ is an exact conservative functor,
we see that for any object $M$ of $\Fun(S,\Mod_k)$
there is a canonical fiber sequence
\[
\tau_{\ge0}^S(M)\to M\to \tau_{\le-1}^S(M)
\]
where we view $\tau_{\ge0}^S(M)$ and $\tau_{\le-1}^S(M)$
as objects of $\Fun(S,\Mod_{k})$.
Thus, 
$(\Fun(S,\Mod_{k,\ge 0}),\Fun(S,\Mod_{k,\le0}))$ determines
a $t$-structure on $\Fun(S,\Mod_k)$.
Let $\Coh(S)$ be the full subcategory of $\QC(S)$
spanned by objects $C$ such that $C$ is bounded with respect to
the $t$-structure and $H_i(C)$ is finite dimensional for every $i\in \ZZ$
(after pulling back to the heart of $\Mod_k$).
Namely, if $s:\Delta^0\to S$ is a base point,
then $s^*C$ is represented by a bounded (chain) complex such that
$H_i(s^*C)$ is finite dimensional for $i\in \ZZ$.
(Notice that $\Coh(S)$ can be defined to be the full subcategory of {\it dualizable objects}, i.e., perfect complexes.)
The symmetric monoidal $\infty$-category $\Coh^\otimes(S)$ will play a role analogous to
$\textup{Cov}(S)$.
We denote by $\ICoh^\otimes (S):=\textup{Ind}(\Coh^\otimes(S))$ the symmetric monoidal 
stable presentable $\infty$-category of Ind-objects and refer to
it
as the symmetric monoidal $\infty$-category of Ind-coherent complexes
on $S$.
In this subsection, we observe that $\ICoh^\otimes(S)$ is a fine
$\infty$-category (see Proposition~\ref{RRHT1}).
Moreover, we explain how one can obtain the (higher) rational homotopy groups
and the pro-algebraic completion of the fundamental group of $S$
from the associated
derived stack under a certain finiteness condition
(see Theorem~\ref{homotopy}).

\begin{Proposition}
\label{RRHT1}
Let $S$ be a connected space.
Then $\ICoh^\otimes(S)$ is a fine $\infty$-category.
\end{Proposition}

Let $G:=\pi_1(S,s)$ be the fundamental group of $S$ at a fixed base
point $s\in S$.
Let $\textup{B}G$ denote the fundamental groupoid of $S$ and let
$f:S\to \textup{B}G$ be the natural projection.

{\it Proof of Proposition~\ref{RRHT1}.}
Observe first that the heart of $\Coh(S)$ with respect to the $t$-structure
is naturally equivalent to
$\Fun(S,\Vect^{f}_k)$, where $\Vect^{f}_k$ is (the nerve of) the category
of finite dimensional $k$-vector spaces regarded as the complexes
placed in degree zero.
Every functor $S\to \Vect^{f}_k$ factors as $S\to \textup{B}G\to \Vect^f_k$
in a unique way.  More precisely, we have a natural categorical equivalence
$\Fun(S,\Vect^f_k)\simeq \Fun(\textup{B}G,\Vect^f_k)$.
Note that a functor $BG\to \Vect^f_k$ amounts to an action of the group
$G$ on a finite dimensional vector space. Thus if $G_{alg}$ denotes
the pro-algebraic completion of $G$, then by the universal property of 
the completion, $\Fun(\textup{B}G,\Vect^f_k)$ is equivalent to
the category $\textup{Vect}^f(G_{alg})$ of finite dimensional representations of $G_{alg}$
as symmetric monoidal categories.
Similarly, the heart of $\Coh(\textup{B}G)$ is equivalent to
the category of finite dimensional representations of $G_{alg}$,
and the pullback functor $f^*:\Coh(\textup{B}G)\to \Coh(S)$ induces
an identity of the heart when both hearts are identified with
$\textup{Vect}^f(G_{alg})$.
We will prove that the set $P$ of simple objects of $\textup{Vect}^f(G_{alg})$
regarded as objects in the heart of $\Coh(S)$ is a set of wedge-finite generators.
Note that every object of $\Coh(S)$
is compact in $\ICoh(S)$, and every object of the heart of $\Coh(S)$ is wedge-finite
(since wedge-finiteness can be verified after the base change along
$s:\Delta^0\to S$).
Therefore, it is enough to show that $\Coh(S)$ is contained in the smallest stable subcategory $\CCC$ which contains $P$.
To see this, recall that every object of $\textup{Vect}^f(G_{alg})$
has a filtration of finite length whose graded quotients are simple objects.
Hence we find that the heart of $\Coh(S)$ is contained in $\CCC$.
We then proceed by induction on the length with respect to $t$-structure.
Suppose that objects $D$ such that $H_i(D)=0$ for $i<0$ and $i>r$ belong
to $\CCC$. 
Let $C$ be an object in $\Coh(S)$.
Assume that $H_i(C)=0$ for $i<0$ and $i>r+1$.
Then using the $t$-structure
we have the distinguished triangle
\[
\tau_{\ge0}(C[-1])[1] \to C \to H_0(C) \to \tau_{\ge0}(C[-1])[2].
\]
By the assumption $\tau_{\ge0}(C[-1])[1]$ belongs to $\CCC$.
As observed above $H_0(C)$ belongs to $\CCC$.
Thus $C$ lies in $\CCC$. It follows that arbitrary shifts of $C$ lie in $\CCC$, as desired.
\QED

Next applying Theorem~\ref{main1} to $\ICoh^\otimes(S)$,
we will define a derived stack.
Let $P$ be the set of simple objects that belongs to the heart
in $\Coh(S)$ (cf. the proof of Proposition~\ref{RRHT1}).
An object in $P$ can be thought of as a simple local system on $S$.
We fix an order on the set of wedge-finite generators
$P$.
By Proposition~\ref{RRHT1},
$\ICoh^\otimes(S)$ is a fine $\infty$-category.
Invoking Theorem~\ref{main1},
we associate a derived quotient stack $\XX:=[\Spec A/H]$
to $\ICoh^\otimes(S)$ and the ordered set $P$.
Here $A$ is a commutative differential graded algebra and
$H$ is a pro-reductive group.
There is a canonical equivalence $\QC^\otimes(\XX)\simeq \ICoh^\otimes(S)$
as objects in $\CAlg(\Pr^{\textup{L}}_k)$.

\begin{Remark}
There are other choices of $P$.
For example, we can take $P$ to be
the set of all (finite dimensional)
semi-simple objects. When $G$ is a finite group,
we can take
$P$ to be the set of a single faithful finite dimensional
representation
$\{V\}$ of $G$.

From the viewpoint of the reconstruction problem, the initial categorical data should
be a pair
\[
(\ICoh^\otimes(S), s^*:\ICoh^\otimes(S)\to \ICoh^\otimes(\Delta^0)\simeq\Mod_k^{\otimes})
\]
where we think of $\ICoh^\otimes(S)$ and $s^*$ as an object in $\CAlg(\PR_k)$
and a morphism respectively. The set $P$ of simple local systems can be
obtained from the pair as follows. Let $\mathcal{H}_k$ be the
full subcategory of $\Mod_k$ spanned by those objects $E$ such that
$\Hom_{\textup{h}(\Mod_k)}(\uni_k,E)$ is finite dimensional, and
$\Hom_{\textup{h}(\Mod_k)}(\uni_k,E[n])=0$ for $n\neq 0$ where $\uni_k$ is a unit
of $\Mod_k^\otimes$. Put $\mathcal{F}:=(s^*)^{-1}(\mathcal{H}_k)$.
Then $P$ is the set of simple objects in the homotopy category of $\mathcal{F}$ which is a $k$-linear category.
\end{Remark}

Let $\Aut(s^*):\Aff_k^{op}\to \Grp(\widehat{\SSS})$ be the automorphism group
functor which carries $\Spec R$ to the ``space of automorphisms'' of the composite of symmetric monoidal functors
\[
\ICoh^\otimes(S)\stackrel{s^*}{\to} \ICoh^\otimes(*)\simeq \QC^\otimes(*)\simeq \Mod_k^\otimes\stackrel{\otimes_kR}{\longrightarrow} \Mod_R^\otimes
\]
 (see \cite[Section 3]{Tan} for the precise definition).
 Here $s^*$ denotes the pullback along the point $s:\Delta^0=*\to S$,
and $\Grp(\widehat{\SSS})$ denotes the $\infty$-category of group objects
in $\widehat{\SSS}$.
By the main result of \cite{Bar} and the equivalence $\ICoh^\otimes(S)\simeq \QC^\otimes(\XX)$,
$\Aut(s^*)$ is represented by the based loop stack $\Spec k\times_{\XX}\Spec k=\Omega_*\XX$, that is a derived affine group scheme (that is, a group object
in $\Aff_k$, see \cite[Appendix]{Tan}).
Here let us recall how to get base points on stacks from $s^*$ (Remark~\ref{Tandual}).
The point of $\XX=[\Spec A/H]$ is given by $w:[\Spec \Gamma(H)/H]\simeq \Spec k\to [\Spec A/H]$ where we identify the ring of functions $\Gamma(H)$
as the image of the unit of $\Mod_k$ under the right adjoint of the composite
(i.e., the forgetful functor)
$\QC^\otimes(BH)\to \ICoh^\otimes(S)\stackrel{s^*}{\to} \QC^\otimes(*)\simeq \Mod_k^\otimes$.
Moreover, $w^*:\QC^\otimes([\Spec A/H])\to \QC(\Spec k)\simeq \Mod_k^\otimes$ can naturally be identified with $s^*:\ICoh^\otimes(S)\to \Mod_k^\otimes$.
Here, we record the following result:

\begin{Proposition}
The automorphism group functor $\Aut(s^*):\Aff_k^{op}\to \Grp(\widehat{\SSS})$
is representable by the derived affine group scheme
$\Omega_*\XX$ over $k$.
\end{Proposition}

For a derived stack $\YY:\CAlg_k\to \widehat{\mathcal{S}}$ equipped with a base point $y:\Spec k\to \YY$,
we denote by
$\pi_i(\YY,y)$ the sheafification of
the composite 
\[
\CAlg_k^{\textup{dis}}\hookrightarrow \CAlg_k\stackrel{(\YY,y)}{\to} \widehat{\mathcal{S}}_* \stackrel{\pi_i(-)}{\to} \Grp
\]
with respect to flat (fpqc) topology where $\Grp$ is the category of (large)
groups, $\widehat{\mathcal{S}}_*:=\widehat{\mathcal{S}}_{\Delta^0/}$,
and $\pi_i(-)$ is the $i$-th homotopy group with respect to the base point.
We write $\CAlg_k^{\textup{dis}}$ for the full subcategory of $\CAlg_k$ 
spanned by discrete objects, i.e., those objects $C$ such that $H_i(C)=0$
for $i\neq 0$. It is equivalent to the nerve of ordinary commutative $k$-algebras.

\begin{Theorem}
\label{homotopy}
Suppose that $\pi_1(S,s)$ is algebraically good (see below for this notion),
and a universal cover of $S$ has the homotopy type of a
finite CW complex.
Then we have
\[
\pi_i(\XX,w)=
\left\{
\begin{array}{rl}
 \pi_i(S,s)_{uni} &  \mbox{for $i>1$} \\
 \pi_1(S,s)_{alg} & \mbox{for $i=1$}
 \end{array}
 \right.
\]
where $\pi_i(S,s)_{uni}$ is the pro-unipotent completion of $\pi_i(S,s)$.
We remark that for $i>1$
the unipotent algebraic group $\pi_i(S,s)_{uni}$ is isomorphic to
the additive group $\mathbb{G}_a^{\times r_i}$
where $r_i$ is the rank of $\pi_i(S,s)$.
\end{Theorem}

Let $G\to G_{alg}(k)$ be the canonical homomorphism from
the discrete group $G$ to the group of the $k$-valued points of $G_{alg}$.
It gives rise to a morphism $\textup{B}G \to BG_{alg}$, regarding
$\textup{B}G$ as the constant functor, and we have the pullback
functor $\QC^\otimes(BG_{alg})\to \QC^\otimes(\textup{B}G)$ and its restriction
$\Coh^\otimes(BG_{alg})\to \Coh^\otimes(\textup{B}G)$. Here $\Coh^\otimes(BG_{alg})$
is defined in a similar way, i.e., it consists of
bounded complexes with finite dimensional (co)homology.
Consequently, $\Coh^\otimes(BG_{alg})$ coincides with the full subcategory
spanned by dualizable (but not necessarily compact) objects.
Let $G_{red}$ be the maximal pro-reductive quotient of $G_{alg}$, that is,
the pro-reductive completion of $G$. 
The full subcategory $\Coh(BG_{red})$ of $\QC(BG_{red})$ is defined
in a similar way.
In his important work \cite{aff} where the theory of affine stacks
and schematizations of spaces are developed, To\"en
introduced the notion of algebraically goodness, which we will use.
We shall recall this notion.
Let $H^i(G_{alg},-)$ (resp. $H^i(G,-)$) is the $i$-th derived functor of invariants $\textup{Vect}(G_{alg})\to \Vect_k$, $V\mapsto V^{G_{alg}}$
(resp. $\textup{Vect}(G)\to \Vect_k$, $V\mapsto V^{G}$), where $\textup{Vect}(G)$
 is the category of possibly infinite dimensional representations of the discrete group
$G$, and $\Vect_k$ is the category of $k$-vector spaces.
The natural homomorphism $G\to G_{alg}(k)$ induces the
natural transformation $H^i(G_{alg},-)\to H^i(G,-)$.
The group $G$ is said to be algebraically good when
the natural map $H^i(G_{alg},V)\to H^i(G,V)$ is an isomorphism
for every finite dimensional representation $V$ of $G_{alg}$ and every $i\in \ZZ$. Known examples of algebraically good groups
include finite groups, finitely generated free group, finitely
generated abelian groups, fundamental groups of Riemann surfaces and so on.
The proof of the next Lemma is routine and is left to the reader.

\begin{Lemma}
\label{good}
Suppose that $G$ is algebraically good.
The natural functor
$\Coh(BG_{alg})\to \Coh(\textup{B}G)$ is a categorical equivalence.
\end{Lemma}

Unwinding Theorem~\ref{main1} and its proof including the inductive
construction,
we have the following additional properties of $[\Spec A/H]$
and the equivalence $\QC^\otimes([\Spec A/H])\simeq \ICoh^\otimes(S)$.
We can construct (i) a homomorphism $G_{red}\to H$, (ii) a symmetric monoidal $k$-linear functor
$\QC^\otimes(BH)\to \QC^\otimes([\Spec A/H])\simeq \ICoh^\otimes(S)$ which factors as
\[
\QC(BH)\stackrel{t^*}{\to} \QC(BG_{red}) \to \Ind(\Coh(BG))\stackrel{f^*}{\to} \ICoh(S),
\]
where we abuse notation by denoting by $f^*$ the left Kan extension $\Ind(f^*)$ of the restriction
$f^*:\Coh(BG)\to \Coh(S)$ (cf. \cite[5.3.5.10]{HTT}),
and $t$ is the induced morphism $BG_{red}\to BH$.

\begin{Remark}
\label{remap}
There are some more remarks on the properties.
The commutative algebra object $A$ is the image of the unit $\uni_S$ of $\ICoh(S)$
under the lax symmetric monoidal right adjoint
functor of $\QC(BH)\to \ICoh(S)$.
The homomorphism $G_{red}\to H$ is a closed immersion.
In particular, the induced morphism $BG_{red}\to BH$ is an affine morphism
since $H$ and $G_{red}$ are pro-reductive.
To see this, remember that
the essential image of compacts objects in $\QC(BH)$ under
the constructed functor $t^*:\QC(BH)\to \QC(BG_{red})$
forms a set of compact generators of $\QC(BG_{red})$ (we have constructed such a functor).
It follows that the right adjoint functor
$t_*:\QC(BG_{red})\to \QC(BH)$ is conservative.
Suppose that the kernel $G_{red}'$ of $G_{red}\to H$ is non-trivial.
Take a non-trivial irreducible representation $V$ of $G_{red}'$, and
let $h_*V\in \QC(BG_{red})$ be the pushforward along the natural affine morphism $h:BG_{red}'\to BG_{red}$, that is not zero. But since the composite $BG_{red}'\to BH$ factors as $BG_{red}'\to \Spec k\to BH$, we have $t_*h_*V\simeq 0$.
It gives rise to a contradiction. We conclude that $G_{red}\to H$ is a closed
immersion.
\end{Remark}

In the rest of this Section, $G$ is assumed to be algebraically good
and a universal cover of $S$ has the homotopy type of a
finite CW complex.

Let $\pi:U\to S$ be a universal cover of $S$.
Then we have a pullback diagram in $\mathcal{S}$,
\[
\xymatrix{
U \ar[r]^\eta \ar[d]_{\pi} & \ast \ar[d]^q \\
S \ar[r]_f & \textup{B}G
}
\]
where $\ast$ denotes the contractible space.
We also fix a base point $s'$ of $U$ lying over $s$.
The natural morphism $f:S\to \textup{B}G$ induces the adjunction
$f^*:\QC(\textup{B}G)  \rightleftarrows \QC(S):f_*$.
Note that $S\times_{\textup{B}G}\ast\simeq U$ is a simply connected finite CW complex.
Therefore, by \cite[Lemma 3.4, 3.17]{BF},
$\eta_*$ is conservative and preserves small colimits,
and there is an equivalence $\QC^\otimes(U)\simeq \Mod_C^\otimes$
where $C$ is $\eta_*(\uni_U)$
with $\uni_U$ a unit object in $\QC^\otimes(U)$.
The pushforward functor
$\eta_*:\QC(U)\simeq \Mod_C \to \Mod_k\simeq \QC(*)$ 
can be identified with the forgetful functor. Note that
$C$ is equivalent to
the singular cochain complex of $U$ which belongs to
$\Coh(\ast)$ (keep in mind that $U$ is of finite type).
We remark that the definition of $\QC(U)$ in \cite{BFN}, \cite{BF}
is equivalent to
$\Mod(U,\Mod_k)$, that is, our definition of $\QC(U)$.
Next observe that the canonical base change morphism
$q^*f_*(\uni_S) \to \eta_*\pi^*(\uni_S)\simeq C$ is an equivalence
where $\uni_S$ denotes a unit object of $\QC^\otimes(S)$.
To see this, consider the natural equivalences
$\QC(BG)\simeq \lim_{\ast/BG}\QC(\ast)$
and $\QC(S)\simeq \lim_{\ast/BG}\QC(\ast\times_{BG}S)$ where $*/BG$ is the full subcategory of the overcategory $\SSS_{/BG}$
which consists of morphisms $* \to BG$ from  the contractible space,
and
$\lim_{\ast/BG}\QC(\ast)$ is a limit of the constant diagram indexed by
$*/BG$.
By abuse of notation
we denote by $\{E_{\alpha}\}_{*/BG}$ an object
of $\lim_{\ast/BG}\QC(\ast\times_{BG}S)\simeq \QC(S)$
such that $E_\alpha$ is the pullback of $\{E_{\alpha}\}_{*/BG}$
to $\QC(*\times_{BG}S)$ along the base change
$\alpha':\ast\times_{BG}S\to S$ of $\alpha:\ast\to BG$.
Let $\eta_{\alpha}:\ast\times_{BG}S\to \ast$
be the base change of $S\to BG$ along $\alpha:\ast\to BG$.
Let $f_+:\lim_{\ast/BG}\QC(*\times_{BG}S)\to \lim_{\ast/BG}\QC(\ast)$
be a functor informally given by $\{E_{\alpha}\}_{*/BG} \mapsto \{(\eta_{\alpha})_*E_{\alpha}\}_{*/ BG}$. Here
we regard $\{(\eta_{\alpha})_*E_{\alpha}\}_{*/ BG}$
as an object of $\QC(BG)\simeq \lim_{\ast/BG}\QC(\ast)$.
More precisely, 
as in the construction of
the quotient stack $[\Spec A/G]$ from $A\in \CAlg(\QC(BG))$ in Section~\ref{QCC},
one can construct the functor $f_+$ by using coCartesian fibrations and
the relative adjoint functor theorem.
The functor $f_+$ is a right adjoint of $f^*$.
Indeed, for any $N\in \QC(BG)$ and any $M\in \QC(S)$,
there are natural equivalences
\begin{eqnarray*}
\Map_{\QC(S)}(N,f_+M) &\simeq& \lim_{*/BG}\Map_{\QC(\ast\times_{BG}S)}(\alpha^*N, \alpha^*f_+M) \\
  &\simeq& \lim_{*/BG}\Map_{\QC(\ast\times_{BG}S)}(\alpha^*N, (\eta_{\alpha})_*(\alpha')^*M) \\
 &\simeq& \lim_{*/BG}\Map_{\QC(\ast\times_{BG}S)}(\eta_{\alpha}^*\alpha^*N, (\alpha')^*M) \\
 &\simeq& \lim_{*/BG}\Map_{\QC(\ast\times_{BG}S)}((\alpha')^*f^*N, (\alpha')^*M) \\
&\simeq& \Map_{\QC(S)}(f^*N,M).
\end{eqnarray*}
Thus, we see that $f_+$ is a right adjoint $f_*$ of $f^*$.
Consequently, we see that
$f_*(\uni_S)$ lies in $\Coh(\textup{B}G)$.
The restriction $f^*:\Coh(\textup{B}G)\rightleftarrows \Coh(S):f_*$
is an adjunction.
Take left Kan extensions
$\Ind(f^*):\ICoh(\textup{B}G)\to \ICoh(S)$ and $\Ind(f_*):\ICoh(S)\to \ICoh(\textup{B}G)$ of
$\Coh(\textup{B}G)\to \Coh(S)\subset \ICoh(S)$ and $\Coh(S)\to \Coh(\textup{B}G)\subset \ICoh(\textup{B}G)$ respectively (cf. \cite[5.3.5.10]{HTT}).
It gives rise to an adjunction $f^*:\ICoh(\textup{B}G)\rightleftarrows \ICoh(S):f_*$ (we abuse notation by writing $f^*$ and $f_*$ for them).
The natural morphism $g:BG_{alg}\to BG_{red}$ induced by the quotient map $G_{alg}\to G_{red}$ determines the pullback functor $g^*:\Coh(BG_{red})\to \Coh(BG_{alg})$ and its left Kan extension $g^*:\QC(BG_{red})\simeq \Ind(\Coh(BG_{red}))\to \Ind(\Coh(BG_{alg}))$ (we abuse notation). By adjoint functor theorem,
we have a right adjoint functor which we denote by $g_*$.
Therefore there is a sequence of adjunctions
\[
\xymatrix{
\QC(BG_{red}) \simeq \ICoh(B G_{red}) \ar@<2pt>[r]^{g^*} & \ICoh(BG_{alg}) \simeq \ICoh(\textup{B}G)  \ar@<2pt>[l]^{g_*} \ar@<2pt>[r]^(0.7){f^*}  & \ar@<2pt>[l]^(0.3){f_*} \ICoh(S).
}
\]
The middle equivalence follows from Lemma~\ref{good}.
The right adjoint functors $f_*$ and $g_*$ are lax symmetric monoidal
functors, and hence they carry commutative algebra objects to commutative algebra objects.
We regard $\uni_S$ as an object of $\CAlg(\ICoh(S))$
and put $C=f_*\uni_S\in \CAlg(\Coh(BG_{alg}))\simeq \CAlg(\Coh(\textup{B}G))$
(strictly speaking, we abuse notation. The pullback of $C$ to $\Spec k$
is $\eta_*\uni_U$).
Let $B\in \CAlg(\QC(BG_{red}))$ be $g_*f_*\uni_S$. (Observe and keep
in mind that since $C=f_*\uni_S$ lies in $\Coh(BG_{alg})$, $B$ coincides with the image of $C$ under the right adjoint functor $\QC(BG_{alg})\to \QC(BG_{red})$
of the pullback functor
$\QC(BG_{red})\to \QC(BG_{alg})$.)
Recall the adjoint pair $t^*:\QC(BH) \rightleftarrows \QC(BG_{red}):t_*$.
The functor $t_*$ sends $B$ to $A$ in $\CAlg(\QC(BH))$.

\vspace{2mm}

We fix our convention for stacks in the rest of this subsection.
Note that in Theorem~\ref{homotopy} we treat $R$-valued points
of $\XX=[\Spec A/H]$ only for $R\in \CAlg_k^{\textup{dis}}$.
Thus, by the restriction we will consider $[\Spec A/H]$ to be 
a sheaf on $\CAlg_k^{\textup{dis}}$.
Moreover,
flat (fpqc) sheaves $\pi_i(\XX,w)$ $(i\ge0)$ depend only
on a flat hypercomplete sheafification of $\XX$, that is,
a hypercomplete sheafification with respect to the flat topology
on $\CAlg_k^{\textup{dis}}$ (cf. a theorem of Dugger-Hollander-Isaksen \cite{DHI}).
Therefore, we will take $[\Spec A/H]$ to be a quotient (i.e. a geometric realization)
as a flat hypercomplete sheaf, that is,
a hypercomplete sheaf with respect to the flat topology
on $\CAlg_k^{\textup{dis}}$.
Similarly, we take other stacks such as $[\Spec C/G_{alg}]$, $BG_{alg}$
to be flat hypercomplete sheaves on $\CAlg_k^{\textup{dis}}$
(or one may work with flat topology from the beginning of this subsection).
(This definition of stacks agrees with that of \cite{aff}.)

\begin{Lemma}
\label{isomisom}
There exist natural equivalences of stacks
\[
[\Spec C/G_{alg}]\simeq [\Spec B/G_{red}]\simeq [\Spec A/H].
\]
\end{Lemma}

\Proof
We claim that there are equivalences of stacks
regarded as restricted functors $\Fun(\CAlg_k^{\textup{dis}},\widehat{\mathcal{S}})$.
We will prove the first equivalence.
It will suffice to show that the fiber product
$\Spec k\times_{BG_{red}}[\Spec C/G_{alg}]$ is represented by $\Spec B$
where $\Spec k\to BG_{red}$ is the natural projection.
Notice first that the fiber product $BG_{alg}\times_{BG_{red}}\Spec k$
is equivalent to $[G_{red}/G_{alg}]\simeq BG_{uni}$ where
$G_{uni}$
is the unipotent radical of $G_{alg}$.
Hence there is a fiber sequence of
\[
\Spec C\to \Spec k\times_{BG_{red}}[\Spec C/G_{alg}]\simeq [\Spec C/G_{uni}] \to BG_{uni}
\]
in $\Fun(\CAlg_k^{\textup{dis}},\widehat{\mathcal{S}})$.
By virtue of \cite[2.4.1]{aff}, $BG_{uni}$
can be represented by $\Spec E$ (in $\Fun(\CAlg_k^{\textup{dis}},\widehat{\mathcal{S}})$)
such that $E$ is a
coconnective object in $\CAlg_{k}$, i.e., $H_i(E)=0$
for $i>0$, and $H_0(E)=k$;
in \cite{aff} cosimplicial $k$-algebras
are used but the base field $k$ is of characteristic zero,
and so
one can use coconnective commutative differential graded $k$-algebras (cf. \cite[Section 6.4]{Fre}).
Moreover, since 
$\pi_i([\Spec C/G_{uni}],p)$ is represented by
a pro-unipotent group
for any $i\ge1$ and any base point $p$ by the fiber sequence and \cite[2.4.5]{aff}, we deduce
from \cite[2.4.1]{aff} that $[\Spec C/G_{uni}]\simeq \Spec F$ such that $F$ is a
coconnective object in $\CAlg_k$, that is a limit $\lim_{[n]} C\otimes_k\Gamma(G_{uni})^{\otimes n}$ of the cosimplicial diagram given by the Cech nerve associated to $\Spec C\to [\Spec C/G_{uni}]$
(the affinization \cite[Section 2.2]{aff} commutes with colimits).
Note that $F \simeq \lim_{[n]} C\otimes_k\Gamma(G_{uni})^{\otimes n}$
is the image of $C\in \CAlg(\QC(BG_{uni}))$
under $\QC(BG_{uni})\to \QC(\Spec k)$.
Since $C$ is homologically bounded above (as an object $\Mod_k$),
we can apply the base change formula for quasi-coherent complexes
to the pullback diagram
\[
\xymatrix{
BG_{uni} \ar[r] \ar[d] & BG_{alg} \ar[d]^t \\
\Spec k \ar[r] & BG_{red}
}
\]
and conclude that $F\simeq B$ (note also that $B$ is coconnective).
Therefore, it yields a natural equivalence $[\Spec B/G_{red}]\simeq [\Spec C/G_{alg}]$.
Since we have observed that $BG_{red}\to BH$ is an affine morphism
with the fiber $H/G_{red}\simeq BG_{red}\times_{BH}\Spec k$
over the natural morphism $\Spec k\to BH$ (see Remark~\ref{remap}), 
the proof of $[\Spec B/G_{red}]\simeq [\Spec A/H]$ is similar and easier.
\QED

{\it Proof of Theorem~\ref{homotopy}.}
As in the case of $[\Spec A/H]$, we obtain base points $u:\Spec k\to [\Spec C/G_{alg}]$ and $v:\Spec k\to [\Spec B/G_{red}]$ in a similar way. These points commute with the above equivalences
$[\Spec C/G_{alg}]\simeq [\Spec B/G_{red}]\simeq [\Spec A/H]$ up to an equivalence.
Remember that
$C$ is given by the chain complex computing the singular cohomology of $U$.
Namely, $C\simeq k^{U}$ in $\CAlg_k$
($U$ is considered to be an object of $\SSS$,
and the presentable $\infty$-category $\CAlg_k$ is cotensored over $\SSS$).
Note that $\QC^\otimes(U)\simeq \Mod^\otimes_{C}$ (here $C$ forgets the $G_{alg}$-action). If $D\simeq C\otimes_k \Gamma(G_{alg})$ denotes the image of the unit of $\Mod^\otimes_{C}$ under the pushforward
along the composite $U\to S\to BG_{alg}$, we have a natural equivalence
of stacks $\Spec C\simeq [\Spec D/G_{alg}]$.
Thus the point $u:\Spec k\to [\Spec C/G_{alg}]$ factors as
$\Spec k\simeq [\Spec \Gamma(G_{alg})/G_{alg}]\stackrel{u'}{\to} \Spec C\simeq [\Spec D/G_{alg}]\to [\Spec C/G_{alg}]$.
Applying \cite[2.3.3, 2.5.3]{aff} to $\Spec C\simeq \Spec k^{U}$,
we deduce that $\pi_{i}(\Spec C,u') \simeq \pi_{i}(U,s')_{uni}$.
Combining the fiber sequence $\Spec C\to [\Spec C/G_{alg}]\to BG_{alg}$, the vanishing of the higher homotopy group of $BG_{alg}$ and Lemma~\ref{isomisom},
we have isomorphisms
\[
\pi_{i}([\Spec A/H],w)\simeq \pi_{i}([\Spec C/G_{alg}],u) \simeq \pi_{i}(\Spec C,u') \simeq \pi_{i}(U,s')_{uni}
\]
for $i>1$.
For $i=1$,
$\pi_1([\Spec A/H],w)\simeq\pi_1([\Spec C/G_{alg}],u)\simeq G_{alg}=\pi_1(S,s)_{alg}.$
This proves the theorem.
\QED

\begin{Remark}
\label{classicRHT}
There are several formalisms of rational homotopy types
and rationalizations of
non-nilpotent topological space: for example, fibrewise rationalizations,
schematizations, pro-algebraic homotopy types,
Tannakian differential graded categories,
see
\cite{BK}, \cite{BS}, \cite{GHT}, \cite{aff}, \cite{P}, \cite{Mori}.
The author does not know which formalism is the most adequate one
(perhaps it depends on the purposes).
If we consider the fine $\infty$-category $\ICoh^\otimes(S)$
to be a categorical invariant of $S$, then
this viewpoint is similar to the idea by Moriya \cite{Mori}
that the Tannakian differential
graded category $\operatorname{T}_{\textup{PL}}(S)$
(in the sense of {\it loc. cit.})
associated with a topological space $S$ is
a model of a rationalization of $S$.
The main difference is that our construction of $\ICoh^\otimes(S)$
may be viewed as a generalization of $\QC(-)$ or $\Coh(-)$
in algebraic geometry, while the construction in \cite{Mori}
is a generalization of
Sullivan's construction which is a ``de Rham-theoretic" approach.
On the other hand, in light of Theorem~\ref{main1}, we might think that the associated stack
$\mathcal{X}:\CAlg_k^{\textup{dis}}\to \widehat{\SSS}$ is a candidate for a model of
a rationalization of $S$.
It seems quite likely that for an arbitrary topological
space,
the associated stack agrees with
the schematization in the sense of \cite{aff}.
\end{Remark}

\begin{Remark}
It is interesting to compare this subsection with a Tannakian reconstruction
of schemes and Deligne-Mumford stacks discussed in \cite{FI}.
In {\it loc. cit.}, emphasizing  ``derived Tannakian viewpoint''
we give a reconstruction of schemes and Deligne-Mumford stacks
$X$
from $\QC^\otimes(X)$ (without reference to any $t$-structure).
Our approach to rational homotopy theory
in this subsection gives a unified picture.
\end{Remark}

{\it Acknowledgements.}
The author would like to
thank Prof. K. Fujiwara, B. Kahn, K. Kimura, S. Kimura, S. Mochizuki, S. Moriya, S. Saito
and T. Yamazaki
for valuable conversations and comments related to the subject of this paper.
The author would like to express his gratitude to
anonymous referees for careful readings and constructive comments
which substantially improved this paper.
The author is partially supported by Grant-in-Aid for Scientific Research,
Japan Society for the Promotion of Science.

\end{document}